\documentclass[a4paper,reqno]{amsart}

\usepackage{amsfonts}
\usepackage{amssymb}
\usepackage{amsthm}
\usepackage{mathrsfs}
\usepackage[english]{babel}
\usepackage[left=2.5cm,right=2.5cm,top=3.5cm,bottom=3cm,headsep=0.7cm]{geometry}

\usepackage{pgf,tikz}
\usetikzlibrary{arrows}
\usepackage[scriptsize,bf]{caption}
\usepackage{floatrow}
\usepackage{enumerate}
\usepackage{enumitem}
\usepackage{scalerel}
\usepackage{hyperref}

 \usepackage[latin1]{inputenc}

\usepackage[titletoc]{appendix}

\theoremstyle{plain}
\begingroup
\theoremstyle{plain}
\newtheorem{theorem}{Theorem}[section]

\newtheorem{proposition}[theorem]{Proposition}
\newtheorem{lemma}[theorem]{Lemma}
\theoremstyle{definition}
\newtheorem{definition}[theorem]{Definition}
\theoremstyle{remark}
\newtheorem{remark}[theorem]{Remark}

\endgroup

\theoremstyle{definition}
\theoremstyle{remark}

\mathchardef\emptyset="001F

\numberwithin{equation}{section}

\definecolor{ddorange}{rgb}{1,0.5,0}
\definecolor{ddcyan}{rgb}{0,0.2,1.0}

\definecolor{corrections}{RGB}{169,22,159}

\newcommand{\RR}{\mathbb{R}}
\newcommand{\NN}{\mathbb{N}}

\newcommand{\D}{\mathcal{D}}
\newcommand{\A}{\mathscr{A}}

\newcommand{\E}{\mathcal{E}}

\newcommand{\M}{\mathcal{M}}

\renewcommand{\H}{\mathcal{H}}
\renewcommand{\L}{\mathcal{L}}

\newcommand{\R}{\mathcal{R}}

\newcommand{\sm}{\setminus}
\newcommand{\compact}{\Subset}
\newcommand{\de}{\partial}

\renewcommand{\d}[1]{\, \mathrm{d} #1}

\renewcommand{\epsilon}{\varepsilon}

\renewcommand{\tilde}{\widetilde}
\newcommand{\weak}{\rightharpoonup}
\newcommand{\wstar}{\stackrel{*}\rightharpoonup}

\newcommand{\mres}{\mathbin{\vrule height 1.6ex depth 0pt width 0.13ex\vrule height 0.13ex depth 0pt width 1.3ex}}

\newcommand{\x}{{\times}}
\newcommand{\integral}[3]{\int \limits_{#1} \! #2 #3}

\newcommand{\ol}{\overline}
\newcommand{\ul}{\underline}
\newcommand{\seo}{s_\epsilon^\circ}
\newcommand{\teo}{t_\epsilon^\circ}

\newcommand{\Veo}{V_\epsilon^\circ}
\newcommand{\aeo}{\alpha_\epsilon^\circ}
\newcommand{\ueo}{u_\epsilon^\circ}

\renewcommand{\div}{\mathrm{div}}

\newcommand{\bunderline}[1]{\underline{#1\mkern-4mu}\mkern4mu }
\newcommand{\V}{\bunderline{V}}

\DeclareMathOperator*{\esssup}{ess\,sup}
\DeclareMathOperator*{\essvar}{ess\,Var}
\DeclareMathOperator*{\argmin}{argmin}

\makeatletter
\newcommand{\myitem}[1][]{
  \protected@edef\@currentlabel{#1}%
\item[#1]
}
\makeatother

\title[Fatigue effects in elastic materials]{Fatigue effects in elastic materials with variational damage models: A vanishing viscosity approach}
\author{Roberto Alessi, Vito Crismale, Gianluca Orlando}

\address[Roberto Alessi]{Department of Mathematics, SAPIENZA Università di Roma, Piazzale Aldo Moro 5\\
Roma, 00185, Italy}
\email[Roberto Alessi]{roberto.alessi@uniroma1.it}

\address[Vito Crismale]{CMAP, \'Ecole Polytechnique, UMR CNRS 7641, 91128 Palaiseau Cedex, France}
\email[Vito Crismale]{vito.crismale@polytechnique.edu}

\address[Gianluca Orlando]{Zentrum Mathematik, Technische Universit\"at M\"unchen, Boltzmannstrasse 3\\
Garching bei M\"unchen, 85747, Germany}
\email[Gianluca Orlando]{orlando@ma.tum.de}

\begin{document}
\begin{abstract}
\noindent We study the existence of quasistatic evolutions for a family of gradient damage models which take into account fatigue, that is the process of weakening in a material due to repeated applied loads. The main feature of these models is the fact that damage is favoured in regions where the cumulation of the elastic strain (or other relevant variables, depending on the model) is higher. To prove the existence of a quasistatic evolution, we follow a vanishing viscosity approach based on two steps: we first let the time-step $\tau$ of the time-discretisation and later the viscosity parameter $\epsilon$ go to zero. As $\tau \to 0$, we find $\epsilon$-approximate viscous evolutions; then, as $\epsilon \to 0$, we find a rescaled approximate evolution satisfying an energy-dissipation balance.
\end{abstract}
\maketitle
{\small

\bigskip

\bigskip
\bigskip
\maketitle

\setcounter{tocdepth}{1}
\tableofcontents

\section{Introduction}
In Material Science, fatigue refers to the process which leads to the weakening of a material due to repeated applied loads, which individually would be too small to cause the direct failure of the material itself.
Macroscopic fatigue fractures appear as a consequence of the interaction of many and complicated material phenomena occurring at the micro-scale, such as, for instance,  plastic slip systems and coalescence of micro-voids,~\cite{Suresh1998,SchijveB2009,Skibicki2014}.
Fatigue failure is extremely dangerous, since it often occurs without forewarning resulting in devastating events, and is responsible for up to the 90\% of all mechanical failures~\cite{stephens}.
The main reason is that it is very difficult, in real situations, to identify the fatigue degradation state of a material. Therefore, its prediction still represents an open challenge for modeling and simulation at the cutting edge of mechanics.

Fatigue favours the occurrence of damage and fracture in different types of materials, both brittle and ductile. When the stress level is high enough to induce plastic deformations, the material is usually subjected to a so-called low-cycle fatigue regime; instead, high-cycle fatigue occurs if the stress is below the yield stress such that the strains are primarily elastic. Models where fatigue effects are induced by the cumulation of plastic deformations have been recently studied in~\cite{AMV,AMV2,AMMV,AAGVD} and~\cite{Cri16,CriLaz16,CLO}.

In this paper we study a phenomenological material model where damage is the only inelastic phenomenon and the fatigue weakening of the material is a consequence of repeated cycles of elastic deformations. Our work is inspired by the recent paper~\cite{Ale-Vid-DL}, where the authors propose a similar model in the one-dimensional setting and to which the reader is invited to refer to for further mechanical details.

As usual, damage is expressed in terms of a scalar variable which affects the elastic response of the material and may be interpreted as the local percentage of sound interatomical bonds. In contrast to many previous damage models~\cite{FraGar,MieRou,BMR09,ThoMie,Tho13,KRZ13a,KRZ15,KRZ18}, in this paper the dissipation depends not only on the damage variable itself, but also on the history of the evolution. Indeed, damage is favoured in regions where a suitable history variable has a higher value. This history variable is defined pointwise in the body as the cumulation in time of a given function $\zeta$ that may be the strain, or the stress, or the energy density, according to the model. As a consequence, the material may undergo a damage process even if the variable $\zeta$ remains small during the evolution.

We are here interested in proving the existence of quasistatic evolutions for this model in a two-dimensional antiplane shear setting, following a vanishing viscosity approach. To present in detail our result,  before expressing the strong formulation of the model in terms of differential inclusions,  we introduce the time-incremental minimisation problem corresponding to a time discretisation  $t^i_k := i \frac{T}{k} = i\tau_k$ for the unknowns $\alpha \colon \Omega \to [0,1]$ (the damage variable) and $u \colon \Omega \to \RR$ (the displacement) assuming that the previous states $(\alpha^j_k, u^j_k)_{j=0}^{i-1}$ are known:
\[
(\alpha^i_k, u^i_k) \in \argmin_{\alpha \leq \alpha^{i-1}_k} \Big\{\frac{1}{2}\integral{\Omega}{\mu(\alpha) \big|\nabla u\big|^2}{\d x} + \frac{1}{2} \integral{\Omega}{\big| \nabla \alpha \big|^2}{\d x} + \integral{\Omega}{f(V_k^{i-1}) (\alpha_k^{i-1} - \alpha)}{\d x} + \frac{\epsilon}{2\tau_k} \|\alpha - \alpha^{i-1}_k\|^2_{L^2} \Big\}  \, .
\]
The functional minimised above consists of three parts: the internal energy 
\[  \E(\alpha,u) = \frac{1}{2}\integral{\Omega}{\mu(\alpha) \big|\nabla u\big|^2}{\d x} + \frac{1}{2} \integral{\Omega}{\big| \nabla \alpha \big|^2}{\d x} \, ,
\] 
 given by the sum of the elastic energy and the damage regularisation term; the energy dissipated from the previous state 
\begin{equation} \label{eq:death of R}
\integral{\Omega}{f(V_k^{i-1}) (\alpha_k^{i-1} - \alpha)}{\d x}  = \R(\alpha - \alpha^{i-1}_k;V^{i-1}_k) \,,  \qquad \text{with} \quad
\R(\beta;V) = \begin{cases} \displaystyle
{-} \! \integral{\Omega}{f(V) \beta}{\d x} & \text{if } \beta \leq 0 \, , \\
+ \infty & \text{otherwise;}
\end{cases}
\end{equation}
and the viscosity term, depending on a small parameter $\epsilon$. The elastic response is affected by the factor $\mu(\alpha) > 0$, where $\mu$ in nondecreasing in $\alpha$, according to the fact that $\alpha = 1$ represents a sound material and $\alpha = 0$ a completely damaged one. (Notice that the constraint $\alpha \leq \alpha^{i-1}_k$ enforces the irreversibility of the damage process.) The $L^2$ norm of~$\nabla \alpha$ is the usual regularising term in gradient damage models (see the aforementioned works and~\cite{DMOrlToa16,Cri17,CriOrl18} for coupling with plasticity). The dissipation term characterises the present model in comparison to other damage models, since the fatigue term $f(V^{i-1}_k)$ weights the  damage increment. For every~$j$, the history variable $V^j_k$ is defined by
\[
V^j_k := \sum_{h = 1}^j \big| \zeta^h_k - \zeta^{h-1}_k\big| \, ,
\]
where $\zeta_k^h$ represents the elastic strain, or the stress, or the density of the elastic energy at time $t^h_k$. Notice that $V^j_k = \int_0^{t^j_k} |\dot \zeta_k(s)| \d s$, where $\zeta_k(s)$ is the piecewise affine interpolation of $\zeta^h_k$. The function~$f$ is nonincreasing, so that in the minimisation it is more convenient to take $\alpha$ lower where the cumulation $V^{i-1}_k$ is larger. The viscosity term prevents $\alpha^i_k$ to be too far (in $L^2$) from the previous damage state $\alpha^{i-1}_k$. 

The approach that we follow consists of two main steps: as, e.g., in \cite{MieRou,BMR09,ThoMie,Tho13,KRZ13a,KRZ15}, we let first the time-step of the discretisation~$\tau_k$ and later the viscosity parameter $\epsilon$ tend to 0. More precisely, the starting point is to define for every $k$ the discrete-time evolution $(\alpha_{\epsilon,k}(t),u_{\epsilon,k}(t))$ as the piecewise affine interpolation of $(\alpha^i_k,u^i_k)$ and to derive {\em a priori} estimates (cf.\ Proposition~\ref{prop:enhanced estimates}) which guarantee that $\|\alpha_{\epsilon,k}\|_{H^1(0,T;H^1(\Omega))}$, $\|u_{\epsilon,k}\|_{H^1(0,T;W^{1,p}(\Omega))}$ are bounded uniformly with respect to $k$ (not with respect to $\epsilon$) and $\|\alpha_{\epsilon,k}\|_{W^{1,1}(0,T;H^1(\Omega))}$, $\|u_{\epsilon,k}\|_{W^{1,1}(0,T;W^{1,p}(\Omega))}$ are bounded uniformly with respect to $k$  and $\epsilon$, for some $p>2$.
We exploit the {\em a priori} estimates $H^1$ in time to pass to the limit as $k \to +\infty$: for every $\epsilon$ we obtain an {\em $\epsilon$-approximate viscous evolution} $(\alpha_\epsilon(t), u_\epsilon(t))$ characterised by an equilibrium condition in~$u_\epsilon(t)$, a unilateral stability condition in~$\alpha_\epsilon(t)$ (Karush-Kuhn-Tucker inequality),  and an energy-dissipation balance (cf.\ Definition~\ref{def:appreveps}).  This evolution may be expressed in terms of the differential inclusions (cf.\ {\rm (ev1)$_\epsilon$}--{\rm (ev2)$_\epsilon$} in Definition~\ref{def:appreveps} and {\rm (ev3')$_\epsilon$} in Lemma~\ref{le:2202181112})
\begin{align*}
\de_u \E(\alpha_\epsilon(t),u_\epsilon(t)) &= 0 \quad \text{in } H^{-1}_{\de_D \Omega}(\Omega)\, , \\
\de_\alpha \R(\dot \alpha_\epsilon(t); V_\epsilon(t)) + \epsilon \, \dot \alpha_\epsilon(t) + \de_\alpha \E(\alpha_\epsilon(t),u_\epsilon(t)) & \ni 0 \quad \text{in } (H^1(\Omega))' ,
\end{align*}
for a.e.\ $t \in (0,T)$, where $V_\epsilon(t)$ is the history variable associated to the evolution $(\alpha_\epsilon, u_\epsilon)$, $H^{-1}_{\de_D \Omega}(\Omega)$ is the dual of $\{v \in H^1(\Omega) \colon  v = 0 \text{ on } \de_D \Omega\}$, and $\de_\alpha \R(\beta; V)$ is the (convex analysis) subdifferential of $\R(\, \cdot\, ; V)$, i.e.\ $\xi \in \de_\alpha \R(\ol \beta; V)$  if and only if $\R(\ol \beta; V) + \langle \xi, \beta - \ol \beta \rangle \leq \R(\beta;V)$ for every $\beta \in H^1(\Omega)$. (For the expression of $\de_u \E$ and $\de_\alpha \E$ we refer to Lemma~\ref{lemma:deuE}.)

The {\em a priori} estimates $W^{1,1}$ in time allow us to reparametrise the $\epsilon$-approximate viscous evolutions and to obtain a family of equi-Lipschitz evolutions $(\alpha^\circ_\epsilon(s), u^\circ_\epsilon(s))$ in a slower time scale $s$. At this stage we let $\epsilon \to 0$ and obtain an evolution $(\alpha^\circ(s),u^\circ(s))$ together with a reparametrisation function $t^\circ(s)$ that permits the passage from the slow to the original fast time scale $t$. In Theorem~\ref{teo:evoRepar} we prove that $(\alpha^\circ, u^\circ)$ still satisfies an equilibrium condition in~$u^\circ(s)$, a unilateral stability condition in~$\alpha^\circ(s)$  (Karush-Kuhn-Tucker inequality),  and an energy-dissipation balance. However, the dissipation in the energy balance weights the rate of damage with a function $\tilde{f}^\circ(s)\leq f(V^\circ(s))$, where $V^\circ(s)$ is the history variable associated to the evolution $(\alpha^\circ, u^\circ)$.  In terms of differential inclusions, this reads as (cf.\ {\rm (ev1)}--{\rm (ev2)} in Theorem~\ref{teo:evoRepar} and~\eqref{0106182337} in Remark~\ref{0106182044})
\begin{align*}
\de_u \E(\alpha^\circ(s),u^\circ(s)) &= 0 \quad \text{in } H^{-1}_{\de_D \Omega}(\Omega)\, , \\
\de_\alpha \R(\dot \alpha^\circ(s); \tilde f^\circ(s)) + \de_\alpha \E(\alpha^\circ(s),u^\circ(s)) & \ni 0 \quad \text{in } (H^1(\Omega))' ,
\end{align*}
for a.e.\ $s \in (0,S) \sm U^\circ$, where $\R(\, \cdot \, ; \tilde f^\circ(s))$ is defined as in~\eqref{eq:death of R} with $\tilde f^\circ(s)$ in place of $f(V^\circ(s))$.  The set $U^\circ$ corresponds to jump instants (of the evolution in the fast time scale) reparametrised in the slow time scale: therein the limit evolution is governed by a variational inequality of viscous type, representing a fast unstable propagation in the original time scale.  An interesting issue, that we were not able to solve, is to determine whether there are explicit examples where this inequality is strict and $\tilde{f}^\circ(s)$ is actually the correct weight to consider in the energy-dissipation balance.

In the mathematical treatment of the present model some technical difficulties arise. Here we discuss the main issues in the {\em a priori} estimates and in the limits as $k \to +\infty$ and $\epsilon \to 0$.

The proof of the {\em a priori} estimates rests upon the manipulation of the Discrete  Karush-Kuhn-Tucker  conditions~\eqref{eq:eul1} and~\eqref{eq:eul2} evaluated at two subsequent times $t^{i-1}_k$ and $t^i_k$, respectively, as e.g.\ in~\cite{NocSavVer00,KRZ13a,MieZel14,KRZ15,CriLaz16,KRZ18}. The resulting estimate~\eqref{eq:08111609} contains in the right-hand side also discrete-time derivatives at time $t^{i-1}_k$, in contrast to the aforementioned works, where only discrete-time derivatives at time $t^i_k$ appear. These additional terms are due to the presence of the fatigue weight $f(V^{i-1}_k)$ in the dissipation for the $i$-th incremental minimisation problem and prevent the immediate application of the discrete Gronwall estimate used in the previous works. We refine the usual technique to overcome this issue in \eqref{eq:08111721}--\eqref{eq:before discrete Gronwall}.

The main difficulty in deriving the properties of the $\epsilon$-approximate viscous evolutions $(\alpha_\epsilon(t), u_\epsilon(t))$ consists in passing to the limit as $k \to +\infty$ in the dissipation term containing the fatigue weight $f(V_{\epsilon,k}(t))$. The {\em a priori} estimate on $\|u_{\epsilon,k}\|_{H^1(0,T;W^{1,p}(\Omega))}$ only guarantees that $\nabla \dot u_{\epsilon,k} \weak \nabla \dot u_\epsilon$ weakly in $L^2(0,T;L^p(\Omega;\RR^2))$, and this convergence is not sufficient to deduce the convergence of $V_{\epsilon,k}$ to $V_{\epsilon}$, even in the paradigmatic case where~$\zeta$ is the elastic strain, namely when the history variable is $V(t) = \int_0^t |\nabla \dot u(s)| \d s$. To circumvent this problem we first let~$f(V_{\epsilon,k}(t))$ converge to some $\tilde f_\epsilon(t)$ weakly* in $L^\infty(\Omega)$ for every $t$ by an Helly-type theorem (cf.\ Lemma~\ref{le:compCum}), to get an evolution $(\alpha_\epsilon(t), u_\epsilon(t))$ satisfying the  Karush-Kuhn-Tucker  inequality, and the energy-dissipation balance with~$\tilde f_\epsilon(t)$ in place of $f(V_\epsilon(t))$ (cf.\ Propositions~\ref{prop:stweak} and~\ref{prop:balance weak}). At this stage, we exploit the convergence of all the terms of the discrete-time energy-dissipation balance to the corresponding ones in the continuous-time energy-dissipation balance. This improves the convergence of $\dot \alpha_{\epsilon,k}$ to~$\dot \alpha_\epsilon$ (Proposition~\ref{prop:strconveps}), allowing us to deduce that $\nabla \dot u_{\epsilon,k} \to \nabla \dot u_\epsilon$ strongly in $L^2(0,T;L^p(\Omega;\RR^2))$ and thus that $\tilde f_\epsilon(t) = f(V_\epsilon(t))$. Eventually, we obtain the existence of an $\epsilon$-approximate evolution.

The scenario when $\epsilon \to 0$ is radically different. Indeed, here the energy-dissipation balance does not help to improve the weak convergence $\nabla \dot u^\circ_{\epsilon} \weak \nabla \dot u^\circ$ for the rescaled evolutions $(\alpha^\circ_\epsilon(s), u^\circ_\epsilon(s))$, due to the rate-independence of the system as $\epsilon \to 0$. As a consequence, the limit evolution is formulated with $\tilde f^\circ(s)$, the weak$^*$-$L^\infty$ limit of the fatigue weight reparametrisations $f(V^\circ_\varepsilon(s))$, in place of $f(V^\circ(s))$. This motivates why we pass to the limit in two steps, rather than directly taking a simultaneous limit $\tau_k/\epsilon_k \to 0$, $k \to +\infty$, as in the framework developed in~\cite{MinSav18, MieRosSav16} and followed in \cite{KRZ18}.

\section{Assumptions on the model}\label{sec:Ass}

\subsection*{Vector-valued functions.}

In this paragraph we let $X$ be a Banach space. We will often consider the Bochner integral of measurable functions $v \colon [0,T] \to X$. For the definition of this notion of integral and its main properties we refer to \cite[Appendix]{Bre} or to the textbook~\cite{Dun-Sch}. The Lebesgue space $L^p(0,T;X)$ is defined accordingly. We recall that, if $p \in [1,\infty)$ and $X$ is separable, the dual of $L^p(0,T;X)$ is $L^{p'}(0,T;X')$, where $\frac{1}{p} + \frac{1}{p'} = 1$ and $X'$ is the dual of $X$.

For the definition and the main properties of absolute continuous functions $AC([0,T];X)$ and Sobolev functions $W^{1,p}(0,T;X)$, the reader is referred to~\cite[Appendix]{Bre}. We recall here the Aubin-Lions Lemma~\cite{Aub63,Sim87} about the compactness property enjoyed by $W^{1,p}(0,T;X)$. Let $Y$ be a Banach space compactly embedded in $X$, and let $1 \leq p,q \leq \infty$. Then the space $W=\{v \in L^p(0,T;Y) \ : \ \dot v \in L^q(0,T;X) \}$ is: 1) compact in $L^p(0,T;X)$ if $p < \infty$; 2) compact in $C([0,T];X)$ if $p=\infty$ and $q > 1$.

In this paper, the Banach space $X$ will be either a Lebesgue space $L^q(U;\RR^m)$ or a Sobolev space $W^{1,q}(U)$, where $U$ is an open set of $\RR^n$. Given an element $v \in L^p(0,T;L^q(U;\RR^m))$, $p, q \in [1,\infty)$, we identify it with the function $v \colon [0,T] \x \Omega \to \RR^m$ defined by $v(t;x) := \big( v(t) \big) (x)$.

The norms $\| \cdot \|_{L^p}$ and $\| \cdot \|_{W^{1,p}}$ without any further notation will always denote the $L^p$-norm and the $W^{1,p}$-norm with respect to the space variable $x$, respectively.

\subsection*{The reference configuration.} Throughout the paper, $\Omega$ is a bounded, Lipschitz, open set in~$\RR^2$ representing the cross-section of a cylindrical body in the reference configuration. The deformation $v\colon \Omega {\times} \RR \to \Omega {\times} \RR$ takes the form $v(x_1,x_2,x_3) = (x_1,x_2,x_3 \, {+} \, u(x_1,x_2))$, where $u \colon \Omega \to \RR$ is the vertical displacement.  In this antiplane shear framework, the two dimensional setting is the physical relevant one. This assumption gives the compact embedding $H^1(\Omega)$ in $L^p(\Omega)$ for every $p \in [1, \infty)$, which we employ in the {\em a priori} estimates. 

We assume that $\de \Omega = \overline{\de_D \Omega} \cup \overline{\de_N \Omega}$, where $\de_D \Omega$ and $\de_N \Omega$ are relatively open sets in $\de \Omega$ with $\de_D \Omega \cap \de_N \Omega = \emptyset$ and $\H^1(\de_D \Omega) > 0$. A Dirichlet  boundary datum will be prescribed on the set $\de_D\Omega$.

In order to apply the integrability result~\cite{Gro} to our problem (see Remark~\ref{rmk:integrability} below), we assume that $\Omega \cup \de_N \Omega$ is regular in the sense of \cite[Definition 2]{Gro}. (Notice that in dimension 2 this regularity assumption on $\Omega \cup \de_N \Omega$ is satisfied, e.g., when the relative boundary $\de (\de_N \Gamma)$ in $\de \Omega$ consists  of  a finite number of points.)

It is convenient to introduce the notation $W^{-1,p}_{\de_D \Omega}(\Omega)$ for the dual of the space $\{ u \in W^{1,p'}(\Omega) \ : \ u = 0 \text{ on } \de_D \Omega\}$, where $\frac{1}{p} + \frac{1}{p'} = 1$.

\subsection*{The total energy.} Following~\cite{fremond}, the damage state of the body is represented by an internal variable $\alpha \colon \Omega \to [0,1]$. The value $\alpha = 1$ corresponds to a sound state, whereas $\alpha = 0$ corresponds to the maximum possible damage. As usual in gradient damage models \cite{PhaMar}, the system in analysis comprises a regularizing term $\| \nabla \alpha \|^2_{L^2}$. In particular, the damage variable $\alpha$ belongs to the Sobolev space $H^1(\Omega)$.

For every $\alpha \in H^1(\Omega)$ and $u \in H^1(\Omega)$, the {\em stored elastic energy} is defined by
\[
\frac{1}{2} \! \integral{\Omega}{\mu(\alpha) \big| \nabla u \big|^2}{\d x} \, .
\]
We make the following assumptions on the dependence of the shear modulus $\mu$ on the damage variable $\alpha$:
\begin{gather}
\mu \colon \RR \to [0, +\infty) \text{ is a } C^{1,1}(\RR), \text{ nondecreasing function with } \mu(0) > 0\, , \nonumber \\
\mu(\beta) = \mu(0) \text{ for } \beta \leq 0, \quad \mu(\beta) = \mu(1) \text{ for } \beta \geq 1 \, . \label{eq:mu constant}
\end{gather}
The regularity assumption on $\mu$ is needed in the proof of Proposition~\ref{prop:enhanced estimates} (see \eqref{eq:05121748}). The condition \eqref{eq:mu constant} on $\mu$ forces $\alpha$ to take values in $[0,1]$ in the evolution (see Remark~\ref{rem:3005181826}).

The {\em total energy} corresponding to a damage state $\alpha$ and to a displacement $u$ is
\begin{equation} \label{eq:def of E}
\E(\alpha, u) := \frac{1}{2} \! \integral{\Omega}{\mu(\alpha) |\nabla u|^2}{\d x}  + \frac{1}{2} \! \integral{\Omega}{| \nabla \alpha |^2}{\d x} \, .
\end{equation}
Notice that the constant $\frac{1}{2}$ in the gradient damage regularisation term does not play a role in the mathematical treatment and may be replaced by any positive constant.

We compute here the derivatives of the total energy. Note that an integrability strictly higher than 2 is required on~$\nabla u$ to guarantee the differentiability of the energy with respect to $\alpha$.

\begin{lemma} \label{lemma:deuE}
The following statements hold true:
\begin{itemize}
\item[i)] Let $u \in W^{1,p}(\Omega)$, with $p > 2$. Then the functional $\alpha \in H^1(\Omega) \mapsto \E(\alpha, u)$ is differentiable and
\begin{equation} \label{eq:derivative of E}
\langle \de_\alpha \E(\alpha, u), \beta \rangle = \frac{1}{2} \!\integral{\Omega}{\mu'(\alpha) \big| \nabla u \big|^2 \beta}{ \d x} + \integral{\Omega}{\nabla \alpha \cdot \nabla \beta}{\d x} \, ,
\end{equation}
for every $\alpha, \beta \in H^1(\Omega)$.
\item[ii)] Let $\alpha \in H^1(\Omega)$. Then the functional $u \in H^1(\Omega) \mapsto \E(\alpha,u)$ is differentiable 
and
\[
\langle \de_u \E(\alpha, u), v \rangle = \integral{\Omega}{\mu(\alpha) \nabla u \cdot \nabla v}{ \d x}\, ,
\]
for every $v \in H^1(\Omega)$.
\end{itemize}
\end{lemma}
\begin{proof}
We only prove {\em i)}, the proof of {\em ii)} being trivial. The derivative of $\frac{1}{2}\| \nabla \alpha \|^2_{L^2}$ simply gives the second integral in~\eqref{eq:derivative of E}. As for the differentiability of $\integral{}{\mu(\alpha)|\nabla u|^2}{\d x}$, let us fix $\alpha, \beta \in H^1(\Omega)$, and $\delta > 0$. By Young's inequality we have
\[
\Big| \frac{\mu(\alpha+\delta\beta) - \mu(\alpha)}{\delta}\big| \nabla u \big|^2 \Big| \leq \|\mu'\|_{L^\infty} |\beta| \ \big|\nabla u\big|^2  \leq C \Big[ |\beta|^q + \big| \nabla u \big|^p \Big] \, ,
\]
where $q = \frac{p}{p-2} < \infty$. Thanks to the embedding $H^1(\Omega) \compact L^q(\Omega)$, we can apply the Dominated Convergence Theorem to deduce that the functional $\alpha \in H^1(\Omega) \mapsto \E(\alpha, u)$ is G\^ateaux-differentiable and its G\^ateaux-differential is expressed by \eqref{eq:derivative of E}. Moreover, since $u\in W^{1,p}(\Omega)$, with $p>2$, and $H^1(\Omega) \compact
 L^r(\Omega)$, for any $r \in [1,\infty)$, it is immediate that the functionals in i) and ii) are Fr\'echet-differentiable.
\end{proof}

\subsection*{Fatigue and damage dissipation.} The damage dissipation is affected by the cumulation of a suitable variable of the system during the history of the evolution. This variable may be for instance the elastic strain, the stress, or the density of the elastic energy, according to the material model. In the general case, we consider a function depending on the damage variable $\alpha$ and on the elastic strain $\nabla u$: we take, for given evolutions $\alpha \in AC([0,T]; L^q(\Omega;[0,1]))$, $u \in AC([0,T]; W^{1,p}(\Omega))$, with $p>2$, $\frac{1}{q}+\frac{1}{p}<\frac{1}{2}$, the function
\begin{equation}\label{3005181820}
\zeta(t):= g(\alpha(t)) \nabla u(t)\,,
\end{equation}
where $g\in C^{1,1}([0,1])$. (In the following we will guarantee that the damage variable takes values in $[0,1]$, see Remark~\ref{rem:3005181826}; one could also assume $g\in C^{1,1}(\RR)$ and constant in $(-\infty,0]$ and $[1,\infty)$ as done for $\mu$, the difference is that the terms involving $g$ are constant in the incremental minimisation, see \eqref{eq:incpb}.) For instance, if $g(\alpha) \equiv 1$, then $\zeta$ is simply the elastic strain; if $g(\alpha) = \mu(\alpha)$, then $\zeta$ is the stress.

By our assumption on the evolutions $\alpha$, $u$, we have that $\zeta \in AC([0,T]; L^2(\Omega;\RR^2))$, so we consider the corresponding cumulation
\begin{equation}\label{3005181833}
V^\zeta(t;x)\equiv V(t;x):= \int_0^t \big| \dot{\zeta}(s;x)| \d s\,, \quad x\in \Omega\,,
\end{equation}
defined as the Bochner integral in $L^2(\Omega)$.

In \eqref{3005181833} the notation $\equiv$ represents the fact that we do not write in the following the dependence of the cumulated variable from $\zeta$. We shall also use the notation $V_k, V_\epsilon$, etc.\ for the cumulated variable corresponding to $\zeta_k$, $\zeta_\varepsilon$, etc.\, given by \eqref{3005181820} for $\alpha_k$ and $u_k$, $\alpha_\varepsilon$ and $u_\varepsilon$, etc., respectively, specifying the correspondence in each case.

We notice that one could consider other possible choices for the variable $\zeta$, for which the results of this paper still hold. For instance, one could take $\zeta(t)= g(\alpha) \, |\nabla u|^\theta$, with $\theta \in [1,p)$ , so $\zeta \in AC([0,T]; L^{p/\theta}(\Omega))$ (see also the observations in Proposition~\ref{prop:enhanced estimates} and Lemmas~\ref{le:2202181003} and \ref{lemma:desired convergence}). This covers, e.g., the case where $\zeta$ is the density of the elastic energy, i.e., when $g(\alpha) = \mu(\alpha)$ and $\theta = 2$.

 We denote by $H^1_-(\Omega)$ the functions $\beta \in H^1(\Omega)$ with $\beta\leq 0$ a.e.\ in $\Omega$.  For every measurable function $V \colon \Omega \to [0, +\infty)$, playing the role of the cumulation of $\zeta$, and for every $\beta \in H^1_-(\Omega)$, representing the damage rate, we define the corresponding dissipation potential~by
\begin{equation} \label{eq:def of R}
\R(\beta; V) := - \integral{\Omega}{f(V) \beta}{\d x} \, ,
\end{equation}
where
\[
f \colon [0, + \infty) \to [0,+\infty)  \text{ is a Lipschitz,  nonincreasing function with } f(0) > 0 \, .
\]
The regularity assumptions on $f$, $g$ are used in the proof of Proposition~\ref{prop:enhanced estimates} (see \eqref{eq:05121748}), and in Lemmas~\ref{le:2202181003} and~\ref{lemma:desired convergence}.

\begin{figure}[H]
\begin{tikzpicture}
\begin{scope}[scale=1.5]
\draw[->] (-0.3,0) -- (2.9,0);
\draw[->] (0,-0.3) -- (0,1.9);
\draw[line width=1pt] (0,1.5)  to[out=0,in=180-10] (2.6,0.2);
\draw (2.5,-0.3) node {$V$};
\draw (2,1) node {$f(V)$};
\end{scope}
\end{tikzpicture}
\caption{Graph of the function $f$ appearing in the dissipation potential. The higher the value of $V$, the smaller the weight $f(V)$ in the damage dissipation. Recall that $V$ plays the role of the cumulation of the variable $\zeta$.}
\end{figure}
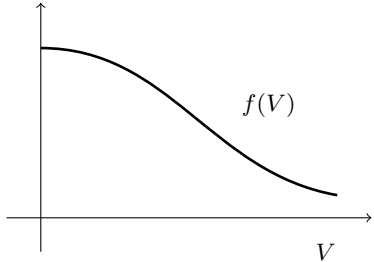

According to the general theory of Rate-Independent systems \cite{MieRouLibro}, $\R$ naturally induces the following dissipation between two damage states  $\alpha_1, \alpha_2 \in H^1(\Omega)$ with $0 \leq \alpha_1 \leq \alpha_2 \leq 1$  a.e.\ in $\Omega$
\begin{equation} \label{eq:def of D}
\D(\alpha_1, \alpha_2; V) :=  \R(\alpha_2 - \alpha_1; V) \, .
\end{equation}

\begin{remark}\label{1207181250}
The dissipation potential $\R$ that we choose here slightly differs from the one proposed in the model of~\cite{Ale-Vid-DL}. In that paper, the dissipation potential features an additional term depending on the gradient of the damage variable. More precisely, using the notation of our paper, a choice more coherent with~\cite{Ale-Vid-DL} would be $\R(\dot \alpha, \nabla \alpha, \nabla \dot \alpha; V) = \int_\Omega f(V)( - \dot \alpha + \nabla \alpha \cdot \nabla \dot \alpha) \d x$. We explain here two reasons that lead us to the decision of not including the term $\nabla \alpha \cdot \nabla \dot \alpha$ in the dissipation potential.

The first reason is a mathematical one. Note that a generic evolution $\alpha(t)$ may not satisfy the inequality $- \dot \alpha + \nabla \alpha \cdot \nabla \dot \alpha \geq 0$; the validity of this condition is however crucial for a physically consistent notion of dissipation potential. Our approach to the problem does not guarantee the {\em a priori} fulfilment of this condition.

The second reason is a modelling one. The model proposed in~\cite{Ale-Vid-DL} is an approach to fatigue fracture via a phase-field model. In a classical phase-field model (without fatigue), the energy dissipated by a fracture is approximated by an energy of the form $\int_\Omega (1-\alpha) + \frac{1}{2} |\nabla \alpha|^2 \d x $, and in that case the term $\int_\Omega \frac{1}{2} |\nabla \alpha|^2 \d x $ should be interpreted as part of the dissipation. This explains why in~\cite{Ale-Vid-DL} the rate of $\frac{1}{2} |\nabla \alpha|^2$ appears in the definition of $\R(\dot \alpha, \nabla \alpha, \nabla \dot \alpha; V)$ and the fatigue weight $f(V)$ also affects this term. In contrast, our aim is to study damage models, whereas the approximation of fracture via damage is not in the scope of this paper. For this reason (as already done in other papers about damage models~\cite{MieRou,BMR09,ThoMie,Tho13,KRZ13a,KRZ15,KRZ18}) we interpret $\int_\Omega \frac{1}{2} |\nabla \alpha|^2 \d x $ as part of the internal energy of the system. In particular, the rate of $\frac{1}{2} |\nabla \alpha|^2$ does not appear in the definition of the dissipation potential.
\end{remark}

\subsection*{Boundary conditions and initial data}

For every  $\overline{\alpha} \in H^1(\Omega)$ with $0 \leq \overline{\alpha} \leq 1$ a.e.\ in $\Omega$  and for every $w \in H^1(\Omega)$, the set of admissible pairs $(\alpha, u)$ with respect to the damage variable $\overline{\alpha}$ and the boundary datum $w$ is defined by:
\[
\A(\overline{\alpha}, w) := \{ (\alpha, u) \in H^1(\Omega) \x H^1(\Omega) \ : \ 0 \leq \alpha \leq \overline{\alpha} \text{ a.e.\ in } \Omega, \ u = w \text{ on } \de_D \Omega \} \, .
\]

The quasistatic evolution will be driven by a boundary datum satisfying
\begin{equation} \label{ass:w}
w \in H^1\big(0,T; W^{1,\tilde p}(\Omega)\big) \, ,
\end{equation}
where $\tilde p > 2$ is a suitable exponent that  is chosen according to Lemma~\ref{lemma:increment of u}. The $\tilde{p}$ integrability  of $\nabla w$ is needed to control the increments of the displacement with those of the damage variable (cf.\ Lemma~\ref{lemma:increment of u}). The $H^1$ regularity in time of the boundary datum is needed for the proof of the {\em a priori} bounds in Proposition~\ref{prop:enhanced estimates} (see \eqref{eq:22031225}).

We prescribe initial conditions $\alpha_0 \in H^1(\Omega)$ and  $u_0 \in  W^{1,\tilde{p}}(\Omega)$  at time $t=0$. We assume, consistently with~\eqref{3005181833}, that the initial cumulation $V_0 = 0$ for notation simplicity. Taking a generic initial cumulation $V_0 \in L^2(\Omega)$ with $V_0 \geq 0$ a.e.\ in $\Omega$ does not entail any mathematical difficulty. (Note that, in that case, definition~\eqref{3005181833} should be modified accordingly by adding the initial cumulation $V_0$.)

 We require
\begin{equation}\label{2705181010}
\de_\alpha \E(\alpha_0, u_0) \in L^2(\Omega)\,.
\end{equation}
 Notice that one could also assume that $\alpha_0$ and $u_0$ are stable with $V_{-1} = 0$, so that the Euler conditions in Lemma~\ref{lemma:euler} hold for $i=0$ too.
The assumption \eqref{2705181010}  is slightly more general, since, for instance, the initial condition $\alpha_0 = 0, u_0 = 0$ is always admissible, no matter whether it is stable or not.

\section{Incremental minimum problems}

\subsection*{Construction of discrete-time evolutions} We fix a sequence of subdivisions $(t^i_k)_{i=0}^{k}$ of the interval $[0,T]$, where $t^i_k := \frac{i}{k}T$ are equispaced nodes. We denote the step of the time discretisation by $\tau_k = \frac{1}{k}$. For notational simplicity, we omit the dependence of $\tau_k$ on $k$ and we use the symbol $\tau$. Moreover, we fix $\epsilon > 0$.

We define the discretisation of the boundary datum $w$ by $w^i_k := w(t^i_k)$, $i =0,\dots, k$.

Let $\alpha^0_k := \alpha_0$, $u^0_k := u_0$, $\zeta_k^0:=g(\alpha_0) \nabla u_0$, and $V^0_k := V_0 =0$. Assuming that we know $\alpha^{i-1}_k$ and $V^{i-1}_k$, we define $(\alpha^i_k, u^i_k)$ as a solution to the incremental minimisation problem (cf.\ \eqref{eq:def of E}, \eqref{eq:def of R},  \eqref{eq:def of D} for the definition of $\E$ and $\D$)
\begin{equation} \label{eq:incpb}
\min \big\{\E(\alpha, u) + \D(\alpha, \alpha^{i-1}_k; V^{i-1}_k) + \tfrac{\epsilon}{2\tau} \|\alpha - \alpha^{i-1}_k\|^2_{L^2} \ : \ (\alpha, u) \in \A(\alpha_k^{i-1}, w_k^i) \big\}
\end{equation}
and we set $\zeta_k^i:= g(\alpha_k^{i}) \nabla u^{i}_k$ and
\[
V^i_k := V^{i-1}_k + \big|\zeta^i_k - \zeta^{i-1}_k\big| = \sum_{j = 1}^i \big| \zeta^j_k - \zeta^{j-1}_k\big| \, .
\]
The existence of a solution to~\eqref{eq:incpb} is obtained by employing the direct method of the Calculus of Variations.

\begin{remark}\label{rem:3005181826}
It is immediate to see that $\alpha^i_k$ is a solution to the problem
\begin{equation} \label{eq:pb alpha 1}
\min \big\{\E(\alpha, u^i_k) + \D(\alpha, \alpha^{i-1}_k; V^{i-1}_k) + \tfrac{\epsilon}{2\tau} \|\alpha - \alpha^{i-1}_k\|^2_{L^2} \ : \  \alpha \in H^1(\Omega), \ 0 \leq \alpha \leq \alpha^{i-1}_k  \leq 1  \big\} \, ,
\end{equation}
where $u = u^i_k$ is fixed. Notice that $\alpha^i_k$ is also a solution to the problem
\begin{equation} \label{eq:pb alpha 2}
\min \big\{\E(\alpha, u^i_k) + \D(\alpha, \alpha^{i-1}_k; V^{i-1}_k) + \tfrac{\epsilon}{2\tau} \|\alpha - \alpha^{i-1}_k\|^2_{L^2} \ : \ \alpha \in H^1(\Omega), \ \alpha \leq \alpha^{i-1}_k \big\} \, ,
\end{equation}
where also competitors $\alpha$ with negative values are taken into account.  Indeed, let us fix a competitor for the problem \eqref{eq:pb alpha 2}, namely $\alpha \in H^1(\Omega; \RR)$ with $\alpha \leq \alpha^{i-1}_k$ and let us set $\alpha^+ := \max \{ \alpha, 0 \}$.
 We employ the fact that $\alpha^+$ is a competitor for \eqref{eq:pb alpha 1}, the assumption \eqref{eq:mu constant}, and the fact that $\alpha^{i-1}_k \geq 0$ to obtain
\[
\begin{split}
\E(\alpha^i_k, u^i_k) & + \D(\alpha^i_k, \alpha^{i-1}_k; V^{i-1}_k) + \tfrac{\epsilon}{2\tau} \|\alpha^i_k - \alpha^{i-1}_k\|^2_{L^2} \\
& \leq \E(\alpha^+, u^i_k) + \D(\alpha^+, \alpha^{i-1}_k; V^{i-1}_k) + \tfrac{\epsilon}{2\tau} \|\alpha^+ - \alpha^{i-1}_k\|^2_{L^2} \\
 & \leq \E(\alpha, u^i_k) + \D(\alpha, \alpha^{i-1}_k; V^{i-1}_k) + \tfrac{\epsilon}{2\tau} \|\alpha - \alpha^{i-1}_k\|^2_{L^2} \, .
\end{split}
\]
This proves the  equivalence between \eqref{eq:pb alpha 1} and \eqref{eq:pb alpha 2}.
\end{remark}

We define the upper and lower piecewise constant interpolations by
\begin{align*}
\ol t_k(t) &:= t^i_k\,, &  \ol \alpha_k(t) &:= \alpha^i_k\,, & \ol u_k(t) &:= u^i_k \,, & \ol \zeta_k(t) &:= \zeta^i_k \,, & \ol w_k(t)&:=w_k^i\,,  \\
\ul t_k(t) &:= t_k^{i-1}  \, ,  & \ul \alpha_k(t) &:= \alpha^{i-1}_k\,,  & \ul u_k(t) &:= u^{i-1}_k\,, & \ul \zeta_k(t) &:= \zeta_k^{i-1} \,, & \ul w_k(t)&:=w_k^{i-1}\, , \\
\end{align*}
and
\begin{equation*}
\V_k(t) := V^{i-1}_k \quad \text{for } t \in (t^{i-1}_k, t^i_k] \, ,
\end{equation*}
for $i = 1, \dots, k$ and $\ol \alpha_k(0)  := \alpha_0$, $\ol u_k(0) := u_0$, $\V_k(0) := V_0 =0$, while $\ul t_k(T):= T$, $\ul \alpha_k(T)  := \alpha_k^k$, $\ul u_k(T) := u_k^k$.
Moreover, we consider the piecewise affine interpolations defined by
\begin{align*}
\alpha_k(t) &:= \alpha^{i-1}_k + (t-t^{i-1}_k) \dot \alpha^i_k \, , \\
u_k(t) &:= u^{i-1}_k + (t-t^{i-1}_k) \dot u^i_k \, , \\
\zeta_k(t) &:= \zeta_k^{i-1} + (t-t^{i-1}_k) \dot \zeta^i_k \, , \quad \text{for } t \in [t^{i-1}_k, t^i_k] \, ,
\end{align*}
for $i = 1,\dots, k$, where
\[
\dot \alpha^i_k := \tfrac{\alpha^i_k - \alpha^{i-1}_k}{\tau} \, , \quad \dot u^i_k := \tfrac{u^i_k - u^{i-1}_k}{\tau}\,,  \quad \dot \zeta^i_k := \tfrac{\zeta^i_k - \zeta^{i-1}_k}{\tau}\,,
\]
and define $w_k$ as the affine interpolation in time of $w$.
We set also
\begin{equation}\label{3005181946}
V_k(t):= \V_k(t)+ \frac{t- \ul t_k(t)}{\tau} \big| \zeta_k(t)- \zeta_k(\ul t_k(t)) \big|\,.
\end{equation}
It is not difficult to verify that Proposition~\ref{prop:essvar} yields
\begin{equation}\label{3105180902}
V_k(t)= \int_0^t \big| \dot{\zeta}_k(s) \big| \d s
\end{equation}
in the sense of Bochner integral in $L^2(\Omega)$.

Note that in the above definitions we dropped the dependence on $\epsilon$ for notation simplicity.

\subsection*{{\em A priori} bounds on discrete-time evolutions} We start the analysis of the discrete evolutions by deducing higher integrability properties of the strain. Following the idea of previous papers (see, e.g., \cite{KRZ13a}), we apply a result proved  in~\cite[Theorem~1]{Gro} (see also \cite[Theorem 1.1]{Her-Mey-Wac} for an extension to the case of elliptic systems with the symmetric gradient in place of $\nabla u$)  regarding the integrability of solutions to elliptic systems with measurable coefficients and with mixed boundary conditions.

\begin{remark}\label{rmk:integrability}
 By \cite[Theorem~1]{Gro},  there exist a constant $C > 0$ and $ \tilde p > 2$ depending on $\|\mu\|_{L^\infty}$ such that the following property is satisfied: for every $\alpha \in H^1(\Omega)$, for every $p \in [2,\tilde p \, ]$, and for every $\ell \in W^{-1,p}_{\de_D \Omega}(\Omega)$, the weak solution $v \in W^{1,p}(\Omega)$ to the problem
\[
\left\{
\begin{aligned}
\div\big( \mu(\alpha) \nabla v \big) &= \ell && \quad  \text{in } \Omega \, , \\
 v & = 0  && \quad \text{on } \de_D \Omega
\end{aligned}
  \right.
\]
satisfies
\[
\| v \|_{W^{1,p}} \leq C \| \ell \|_{W^{-1,p}_{\de_D \Omega}} \, .
\]
\end{remark}

In the following lemma we apply  the regularity given by Remark~\ref{rmk:integrability} to deduce higher integrability of $\nabla u_k(t)$ and to control the increments of the displacement $u$ with the increments of the damage variable $\alpha$.

\begin{lemma}[Higher integrability of the strain] \label{lemma:increment of u}
There exist $\tilde p > 2$ (depending only on $\|\mu\|_{L^\infty}$) and a constant $C > 0$ (depending only on $\|\mu\|_{L^\infty}$, $\|\mu'\|_{L^\infty}$, and $\|w\|_{L^\infty(0,T;W^{1,\tilde p}(\Omega))}$) such that
\begin{subequations}\label{3105181044}
\begin{align}
\| \ol u_k(t) \|_{W^{1,\tilde p}} + \| \ul u_k(t) \|_{W^{1,\tilde p}} + \| u_k(t) \|_{W^{1,\tilde p}} \leq C  \, ,  &\quad \text{for } t \in [0,T] \label{eq:estimate of u} \\
\| \dot u_k(t) \|_{W^{1,p}}  \leq C \Big[ \|\dot \alpha_k(t)\|_{L^q} + \|\dot w_k(t) \|_{W^{1,\tilde p}} \Big] \, , &\quad \text{for } t \in (0,T) \sm \{t^1_k, \dots, t^{k-1}_k\}  \, , \label{eq:estimate of increment}
\end{align}
\end{subequations}
for every $p \in [2,\tilde p \, )$, where $q = \frac{p \tilde p}{\tilde p-p}$.
\end{lemma}
\begin{proof}
Let $\tilde p > 2$ be the exponent given in Remark~\ref{rmk:integrability}. To prove~\eqref{eq:estimate of u}, let us fix $t \in (t^{i-1}_k,t^i_k]$ for $i \in \{1, \dots, k\}$ (notice that the inequality is trivial for $t=0$). By \eqref{eq:incpb}, the function $u^i_k$ minimises $\E(\alpha^i_k, u)$ among all $u \in H^1(\Omega)$ with $u = w^i_k$ on $\de_D \Omega$. Therefore $u^i_k$ is a weak solution to the problem
\begin{equation} \label{eq:pb of uik}
\left\{
\begin{aligned}
\div\big( \mu(\alpha^i_k) \nabla u^i_k \big) &= 0 && \quad  \text{in } \Omega \, , \\
 u^i_k & = w^i_k  && \quad \text{on } \de_D \Omega \, .
\end{aligned}
  \right.
\end{equation}
By Remark~\ref{rmk:integrability}, we have that
\[
\|u^i_k - w^i_k \|_{W^{1,\tilde p}} \leq C \big\| \div\big( \mu(\alpha^i_k) \nabla w^i_k \big) \big\|_{W^{-1,\tilde p}_{\de_D\Omega}} \leq C \|\mu\|_{L^\infty} \|w^i_k\|_{W^{1,\tilde p}},
\]
which implies \eqref{eq:estimate of u} (recall the definition of $\ol u_k$, $\ul u_k$, $u_k$ in terms of the family of $u_k^i$).

To prove~\eqref{eq:estimate of increment}, let us fix $p \in [2, \tilde p \, )$ and $t \in (t^{i-1}_k,t^i_k)$ for $t \in \{1,\dots,k\}$. By \eqref{eq:pb of uik} for $i$ and $i-1$ we get that the function $v := u^i_k - u^{i-1}_k - w^i_k + w^{i-1}_k$ is a weak solution to the problem
\begin{equation} \label{eq:pb of dot uk}
\left\{
\begin{aligned}
\div\big( \mu(\alpha^{i-1}_k) \nabla v \big) &= \ell && \quad  \text{in } \Omega \, , \\
 v & = 0  && \quad \text{on } \de_D \Omega \, ,
\end{aligned}
  \right.
\end{equation}
where $\ell := \div\big( (\mu(\alpha^{i-1}_k) - \mu(\alpha^i_k) ) \nabla u^i_k \big) - \div\big( \mu(\alpha^{i-1}_k) (\nabla w^i_k - \nabla w^{i-1}_k) \big)$. Notice that $\ell \in W^{-1,p}_{\de_D\Omega}(\Omega)$ by \eqref{eq:estimate of u}. By Remark~\ref{rmk:integrability} and by H\"older's inequality we deduce that
\[
\begin{split}
\|v\|_{W^{1,p}} & \leq C \|\ell\|_{W^{-1,p}_{\de_D\Omega}} \leq C \Big[ \| (\mu(\alpha^{i-1}_k) - \mu(\alpha^i_k) ) \nabla u^i_k  \|_{L^p}  + \| \mu(\alpha^{i-1}_k) (\nabla w^i_k - \nabla w^{i-1}_k) \|_{L^p} \Big] \\
& \leq C \Big[ \|\mu'\|_{L^\infty} \|\alpha^i_k - \alpha^{i-1}_k\|_{L^q} \| \nabla u^i_k \|_{L^{\tilde p}} + \|\mu\|_{L^\infty} \|w^i_k - w^{i-1}_k\|_{W^{1,\tilde p}} \Big] \, ,
\end{split}
\]
since $\frac{1}{q} =  \frac{1}{p} - \frac{1}{\tilde p}$. By \eqref{eq:estimate of u} and dividing by $\tau$ we conclude that
\[
\|\dot u^i_k \|_{W^{1,p}} \leq C \Big[ \| \dot \alpha^i_k \|_{L^q} + \| \dot w^i_k \|_{W^{1,\tilde p}} \Big] \, ,
\]
hence the thesis.
\end{proof}

We are now in a position to derive the Euler conditions satisfied by the damage variable in the discrete evolutions. These conditions are also called Discrete  Karush-Kuhn-Tucker  conditions, since we have a constraint of unidirectionality on the damage variable. They are a fundamental ingredient to deduce the {\em a priori} bounds in Proposition~\ref{prop:enhanced estimates}.

\begin{lemma}[Euler conditions] \label{lemma:euler}
For every $t \in (0,T) \sm \{t^1_k, \dots, t^{k-1}_k\}$ we have
\begin{equation} \label{eq:eul1}
\langle \de_\alpha \E(\ol \alpha_k(t), \ol u_k(t)), \beta \rangle + \R(\beta; \V_k(t)) + \epsilon \langle \dot \alpha_k(t) , \beta \rangle_{L^2} \geq 0 \, 
\end{equation}
for every $\beta \in H^1(\Omega)$ such that $\beta \leq 0$ a.e.\ in $\Omega$. Moreover
\begin{equation} \label{eq:eul2}
\langle \de_\alpha \E(\ol \alpha_k(t), \ol u_k(t)), \dot \alpha_k(t) \rangle + \R(\dot \alpha_k(t); \V_k(t)) + \epsilon \| \dot \alpha_k(t)\|^2_{L^2} = 0 \, .
\end{equation}

\end{lemma}
\begin{proof}
Let us fix $t \in (t^{i-1}_k, t^i_k)$ for some $i \in \{1,\dots, k\}$. Let $\beta \in H^1(\Omega)$ with $\beta \leq 0$ a.e.\ in $\Omega$ and let $\delta > 0$. Since $\alpha^i_k$ solves \eqref{eq:pb alpha 2} and $\alpha^i_k + \delta \beta \leq \alpha^{i-1}_k$, we get
\[
\begin{split}
0 & \leq \E(\alpha^i_k + \delta \beta, u^i_k) + \D(\alpha^i_k + \delta \beta , \alpha^{i-1}_k; V^{i-1}_k) + \tfrac{\epsilon}{2 \tau} \|\alpha^i_k + \delta \beta - \alpha^{i-1}_k \|^2_{L^2} + \\
& \quad - \E(\alpha^i_k, u^i_k) - \D(\alpha^i_k, \alpha^{i-1}_k; V^{i-1}_k) - \tfrac{\epsilon}{2 \tau} \|\alpha^i_k - \alpha^{i-1}_k \|^2_{L^2}\, .
\end{split}
\]
Dividing by $\delta$ and letting $\delta \to 0^+$, by \eqref{eq:derivative of E} we get
\[
\frac{1}{2} \!\integral{\Omega}{\mu'(\alpha^i_k) \big| \nabla u^i_k \big|^2 \beta}{ \d x} + \integral{\Omega}{\nabla \alpha^i_k \cdot \nabla \beta}{\d x} - \integral{\Omega}{f(V^{i-1}_k) \, \beta}{\d x} + \epsilon \integral{\Omega}{ \dot \alpha^i_k   \, \beta}{\d x} \geq 0 \, .
\]
This concludes the proof of \eqref{eq:eul1}.

To prove~\eqref{eq:eul2}, notice that $\alpha^i_k - \delta \dot \alpha^i_k \leq  \alpha^{i-1}_k$ for $0< \delta < \tau$. Since $\alpha^i_k$ solves \eqref{eq:pb alpha 2} we get that
\[
\begin{split}
0 & \leq \E(\alpha^i_k - \delta \dot \alpha^i_k, u^i_k) + \D(\alpha^i_k - \delta \dot \alpha^i_k , \alpha^{i-1}_k; V^{i-1}_k) + \tfrac{\epsilon}{2 \tau} \|\alpha^i_k - \delta \dot \alpha^i_k - \alpha^{i-1}_k \|^2_{L^2} + \\
& \quad - \E(\alpha^i_k, u^i_k) - \D(\alpha^i_k, \alpha^{i-1}_k; V^{i-1}_k) - \tfrac{\epsilon}{2 \tau} \|\alpha^i_k - \alpha^{i-1}_k \|^2_{L^2} \, .
\end{split}
\]
Dividing by $\delta$ and letting $\delta \to 0^+$, by~\eqref{eq:derivative of E} this implies~\eqref{eq:eul2}.
\end{proof}
The following proposition ensures that the evolution of $\alpha$ and $u$ is $H^1$ in time uniformly in $k$ for fixed $\varepsilon$, and $AC$ in time uniformly in $k$ and $\varepsilon$, with values in the target spaces $H^1(\Omega)$ and $W^{1,p}(\Omega)$.
\begin{proposition}[{\em A priori} bounds] \label{prop:enhanced estimates}
Let $\tilde{p}$ be as in Lemma~\ref{lemma:increment of u}. There exists a positive constant $C$ independent of $\epsilon$, $k$, and $t$ such that for every $\epsilon > 0$, $k \in \NN$, $t \in (0,T) \sm \{t^1_k, \dots, t^{k-1}_k\}$, $p < \tilde{p}$, it holds that
\begin{align}
\epsilon \| \dot \alpha_k(t) \|_{L^2} & \leq C \exp \Big( C \tfrac{\ol \tau_k(t)}{\epsilon}  \Big)  \, , \label{eq:estimate 1} \\
\epsilon \int_0^{\ol \tau_k(t)} \! \| \dot \alpha_k(s) \|^2_{H^1} \, \d s  + \epsilon \int_0^{\ol \tau_k(t)} \! \| \dot u_k(s) \|^2_{W^{1,p}} \, \d s & \leq C \exp \Big( C \tfrac{\ol \tau_k(t)}{\epsilon}   \Big) \, ,  \label{eq:estimate 2} \\
\int_0^{T} \! \| \dot \alpha_k(s) \|_{H^1} \, \d s  + \int_0^{T} \! \| \dot u_k(s) \|_{W^{1,p}} \, \d s & \leq  C  \, . \label{eq:estimate 3}
\end{align}
\end{proposition}

\begin{proof}
We only need to show the estimates on $\alpha_k(t)$, since the estimates on $u_k(t)$ simply follow from~\eqref{eq:estimate of increment}.

We start with computations which are common in the proofs of all the three inequalities in the statement. The starting point is to obtain an estimate on the time increments of $\dot \alpha_k(t)$ by testing the Euler equations at two subsequent times of the time discretisation. To do so, we fix $i \in \{2,\dots, k\}$.  The case $i = 1$ requires slightly different arguments.
By \eqref{eq:eul2} evaluated at a time $t \in (t^{i-1}_k, t^i_k)$ we get that
\[
\langle \de_\alpha \E(\alpha^i_k, u^i_k), \dot \alpha^i_k \rangle + \R(\dot \alpha^i_k;  V^{i-1}_k) + \epsilon \| \dot \alpha^i_k\|^2_{L^2} = 0 \, .
\]
On the other hand, by testing \eqref{eq:eul1} with $\beta = \dot \alpha^i_k$ at a time $t \in (t^{i-2}_k,t^{i-1}_k)$, we get
\[
\langle \de_\alpha \E(\alpha^{i-1}_k, u^{i-1}_k), \dot \alpha^i_k \rangle + \R(\dot \alpha^i_k; V^{i-2}_k) + \epsilon \langle \dot \alpha^{i-1}_k,\dot \alpha^i_k \rangle_{L^2} \geq 0 \, .
\]
Subtracting the second inequality from the first one, we infer that
\[
\langle \de_\alpha \E(\alpha^i_k, u^i_k) - \de_\alpha \E(\alpha^{i-1}_k, u^{i-1}_k), \dot \alpha^i_k \rangle + \R(\dot \alpha^i_k; V^{i-1}_k) - \R(\dot \alpha^i_k; V^{i-2}_k) + \epsilon \langle \dot \alpha^i_k - \dot \alpha^{i-1}_k , \dot \alpha^i_k \rangle_{L^2} \leq 0 \, ,
\]
namely,
\begin{equation} \label{eq:05121748}
\begin{split}
\epsilon \langle & \dot \alpha^i_k - \dot \alpha^{i-1}_k , \dot \alpha^i_k \rangle_{L^2} + \integral{\Omega}{\big( \nabla \alpha^i_k - \nabla \alpha^{i-1}_k \big) \cdot \nabla \dot \alpha^i_k}{\d x} \\
&\leq \frac{1}{2} \!\integral{\Omega}{\Big[ \mu'(\alpha^{i-1}_k)\big| \nabla u^{i-1}_k \big|^2 - \mu'(\alpha^i_k)\big| \nabla u^i_k \big|^2 \Big] \dot \alpha^i_k}{\d x} + \integral{\Omega}{\Big[ f(V^{i-1}_k) - f(V^{i-2}_k)\Big] \dot \alpha^i_k}{\d x} \\
& \leq \frac{1}{2} \!\integral{\Omega}{ \mu'(\alpha^{i-1}_k)\Big[ \big| \nabla u^{i-1}_k \big|^2 - \big| \nabla u^i_k \big|^2 \Big] \dot \alpha^i_k}{\d x} +  \frac{1}{2} \!\integral{\Omega}{\Big[ \mu'(\alpha^{i-1}_k) - \mu'(\alpha^i_k) \Big] \big| \nabla u^i_k \big|^2 \dot \alpha^i_k}{\d x} \\
& \quad + \|f'\|_{L^\infty}\integral{\Omega}{\big| V^{i-1}_k - V^{i-2}_k \big| | \dot \alpha^i_k|}{\d x} \\
& \leq \frac{1}{2}\|\mu'\|_{L^\infty} \! \integral{\Omega}{\big| \nabla u^i_k + \nabla u^i_k \big| \big| \nabla u^i_k - \nabla u^{i-1}_k \big| |\dot \alpha^i_k|}{\d x} + \frac{1}{2}\|\mu''\|_{L^\infty} \!  \integral{\Omega}{|\alpha^i_k - \alpha^{i-1}_k| \big| \nabla u^i_k \big|^2 |\dot \alpha^i_k|}{\d x}\\
& \quad + \|f'\|_{L^\infty}\integral{\Omega}{\big| \zeta^{i-1}_k - \zeta^{i-2}_k \big| | \dot \alpha^i_k|}{\d x} \\
& \leq C \tau \Big[ \ \big\| \nabla u^i_k + \nabla u^{i-1}_k \big\|_{L^2}  \big\| \nabla \dot u^i_k\big\|_{L^p} \|\dot \alpha^i_k\|_{L^{q_1}} + \big\| \nabla u^i_k  \big\|^2_{L^p} \|\dot \alpha^i_k\|^2_{L^{q_1}} \\& \hspace{5em} + \big\|\nabla u_k^{i-1} \big\|_{L^p}  \|\dot{\alpha}_k^{i-1}\|_{L^{q_1}} \|\dot \alpha^i_k\|_{L^{q_1}} + \big\| \nabla \dot u^{i-1}_k  \big\|_{L^p} \|\dot \alpha^i_k\|_{L^{q_1}} \ \Big] \, .
\end{split}
\end{equation}
In the last inequality we have chosen $p \in (2,\tilde p \, )$ and $q_1 \in (2,\infty)$ such that $\frac{1}{p} + \frac{1}{q_1} = \frac{1}{2}$, and we have employed the identity
\begin{equation*}
\zeta^{i-1}_k - \zeta^{i-2}_k = [g(\alpha_k^{i-1})-g(\alpha_k^{i-2})]\, \nabla u_k^{i-1} + g(\alpha_k^{i-2})[ \nabla u_k^{i-1} - \nabla u_k^{i-2} ]
\end{equation*}
that gives
\begin{equation}\label{0106180946}
\integral{\Omega}{\big| \zeta^{i-1}_k - \zeta^{i-2}_k \big| | \dot \alpha^i_k|}{\d x} \leq \tau \bigg(\|g'\|_{L^\infty}  \big\|\nabla u_k^{i-1} \big\|_{L^p}  \|\dot{\alpha}_k^{i-1}\|_{L^{q_1}} \|\dot \alpha^i_k\|_{L^{q_1}} + \|g\|_{L^\infty}  \big\| \nabla \dot u^{i-1}_k  \big\|_{L^p} \|\dot \alpha^i_k\|_{L^{q_1}} \bigg)\,.
\end{equation}
We remark that, taking $\zeta_k^i:= g(\alpha_k^i) |\nabla u_k^i|^\theta$, with $\theta\in [1,p)$ we could also get the conclusion in \eqref{eq:05121748} with $q'_1 \geq q_1$ such that  $\frac{\theta}{p}+ \frac{1}{q'_1}=\frac{1}{2}$, in place of $q_1$. Indeed
\begin{equation*}
\zeta_k^{i-1}-\zeta_k^{i-2} = [g(\alpha_k^{i-1})-g(\alpha_k^{i-2})]\, |\nabla u_k^{i-1}|^\theta + g(\alpha_k^{i-2})[ |\nabla u_k^{i-1}|^\theta - |\nabla u_k^{i-2}|^\theta ]\,,
\end{equation*}
and since, by the Mean Value Theorem,
\begin{equation*}
\big| |\nabla u_k^{i-1}|^\theta - |\nabla u_k^{i-2}|^\theta \big| \leq \theta (|\nabla u_k^{i-1}|+ |\nabla u_k^{i-2}|)^{\theta-1} \big(|\nabla u_k^{i-1}| - |\nabla u_k^{i-2}|\big)\,,
\end{equation*}
we have that
\begin{equation}\label{0106180947}
\begin{split}
\int \limits_\Omega|\zeta_k^{i-1}-\zeta_k^{i-2}| |\dot{\alpha}_k^i| \d x \leq \tau &\bigg(\|g'\|_{L^\infty}  \big\|\nabla u_k^{i-1} \big\|_{L^p}  \|\dot{\alpha}_k^{i-1}\|_{L^{q'_1}} \|\dot \alpha^i_k\|_{L^{q'_1}} \\ &+ \|g\|_{L^\infty} \theta \big(\|\nabla u_k^{i-1}\|_{L^p}+ \|\nabla u_k^{i-2}\|_{L^p}\| \big) \big\| \nabla \dot u^{i-1}_k  \big\|_{L^p} \|\dot \alpha^i_k\|_{L^{q'_1}} \bigg)\,.
\end{split}
\end{equation}
 Using the fact that $\langle \dot \alpha^i_k - \dot \alpha^{i-1}_k , \dot \alpha^i_k \rangle_{L^2} \geq \| \dot \alpha^i_k \|_{L^2} \Big( \| \dot \alpha^i_k \|_{L^2} - \| \dot \alpha^{i-1}_k \|_{L^2} \Big)$ and by Lemma~\ref{lemma:increment of u} we infer that
\begin{equation} \label{eq:08111609}
\begin{split}
\epsilon & \| \dot \alpha^i_k \|_{L^2} \Big( \| \dot \alpha^i_k \|_{L^2} - \| \dot \alpha^{i-1}_k \|_{L^2} \Big) + \tau \|\nabla \dot \alpha^i_k\|^2_{L^2} \\
& \leq  C \tau \Big[ \| \dot \alpha^i_k \|_{L^{q_1}} \| \dot \alpha^i_k \|_{L^{q_2}}  + \| \dot \alpha^i_k \|^2_{L^{q_1}} + \| \dot \alpha^i_k \|_{L^{q_1}} \| \dot \alpha^{i-1}_k \|_{L^{q_2}} + \| \dot \alpha^i_k \|_{L^{q_1}} \big( \| \dot w^i_k \|_{W^{1,\tilde p}} + \| \dot w^{i-1}_k \|_{W^{1,\tilde p}} \big) \Big]\\
& \leq c_1 \tau \Big[ \| \dot \alpha^i_k \|^2_{L^r} + \| \dot \alpha^i_k \|_{L^r}\| \dot \alpha^{i-1}_k \|_{L^r} +  \| \dot w^i_k \|^2_{W^{1,\tilde p}} +  \| \dot w^{i-1}_k \|^2_{W^{1,\tilde p}} \Big]  \, ,
\end{split}
\end{equation}
where $q_2 := \frac{p \tilde p}{\tilde p - p} \in (2,\infty)$ and $r = \max\{ q'_1,  q_2\} \in (2,\infty)$. We labelled the constant in the last inequality with $c_1$ in order to keep track of it in the sequel. By the compact embedding $H^1(\Omega) \compact L^r(\Omega)$ (notice that $\Omega \subset \RR^2$), we have that for every $\delta > 0$ there exists a constant $C_\delta > 0$ such that for every $\beta \in H^1(\Omega)$
\begin{equation} \label{eq:interpolation inequality}
\| \beta \|^2_{L^r}  \leq \delta \| \nabla \beta \|^2_{L^2} + C_\delta \| \beta \|^2_{L^1} \leq \delta \| \nabla \beta \|^2_{L^2} + C_\delta \| \beta \|_{L^1} \| \beta \|_{L^2} \, .
\end{equation}
Adding $c_1 \, \tau \| \dot \alpha^i_k \|^2_{L^r}  - c_1 \, \tau \| \dot \alpha^i_k \|_{L^r}\| \dot \alpha^{i-1}_k \|_{L^r} + \frac{\tau}{2} \|\dot \alpha^i_k\|^2_{L^2}$ to both sides of \eqref{eq:08111609}, choosing $\delta$ suitably small in the previous inequality  applied to $\beta = \dot \alpha^i_k$,  and multiplying by $\frac{2 \tau}{\epsilon}$ we have that
\begin{equation} \label{eq:08111721}
\begin{split}
2 & \| \dot \alpha^i_k \|_{L^2} \Big( \| \dot \alpha^i_k \|_{L^2} - \| \dot \alpha^{i-1}_k \|_{L^2} \Big) +  2 c_1 \tfrac{\tau}{\epsilon}  \| \dot \alpha^i_k \|_{L^r} \Big( \| \dot \alpha^i_k \|_{L^r} - \| \dot \alpha^{i-1}_k \|_{L^r}  \Big) + \tfrac{\tau}{\epsilon} \|\dot \alpha^i_k\|^2_{H^1}  \\
&\leq c_2 \tfrac{\tau}{\epsilon} \Big( \| \dot w^i_k \|^2_{W^{1,\tilde p}} +  \| \dot w^{i-1}_k \|^2_{W^{1,\tilde p}} \Big) + 2 c_2\tfrac{\tau}{\epsilon} \| \dot \alpha^i_k \|_{L^1} \| \dot \alpha^i_k \|_{L^2}   \, .
\end{split}
\end{equation}

Let us set
\begin{align*}
A_i &:= \Big[ \| \dot \alpha^i_k \|^2_{L^2} + c_1 \tfrac{\tau}{\epsilon} \| \dot \alpha^i_k \|^2_{L^r} \Big]^{\! \frac{1}{2}} \, , && B_i := \sqrt{\tfrac{\tau}{2\epsilon}}\|\dot \alpha^i_k\|_{H^1} \, , \\
 C_i &:= \sqrt{c_2 \tfrac{\tau}{\epsilon}} \Big[ \| \dot w^i_k \|^2_{W^{1,\tilde p}} +  \| \dot w^{i-1}_k \|^2_{W^{1,\tilde p}} \Big]^{\! \frac{1}{2}} \, , && D_i := c_2 \tfrac{\tau}{\epsilon} \| \dot \alpha^i_k \|_{L^1} \, .
\end{align*}
The quantities above are actually defined for every $i = 1, \dots, k$. When $i=1$, we define $C_1:= \sqrt{c_2 \tfrac{\tau}{\epsilon}}\| \dot w^i_k \|_{W^{1,\tilde p}}$. Denoting by $a_i := \big(\| \dot \alpha^i_k \|_{L^2} , \sqrt{c_1\tfrac{\tau}{\epsilon}} \| \dot \alpha^i_k \|_{L^r} \big) $, we get that
\begin{equation} \label{eq:09111538}
\begin{split}
\| \dot \alpha^i_k \|_{L^2} & \Big( \| \dot \alpha^i_k \|_{L^2} - \| \dot \alpha^{i-1}_k \|_{L^2} \Big) +  c_1\tfrac{ \tau}{\epsilon} \| \dot \alpha^i_k \|_{L^r} \Big( \| \dot \alpha^i_k \|_{L^r} - \| \dot \alpha^{i-1}_k \|_{L^r}  \Big) = a_i \cdot (a_i - a_{i-1}) \\
& \geq |a_i| \big( |a_i| - |a_{i-1}|) = A_i (A_i - A_{i-1}) \, .
\end{split}
\end{equation}
Since $\tau \leq \epsilon$ we have that
\begin{equation} \label{eq:09111539}
\begin{split}
\tfrac{\tau}{2\epsilon}\|\dot \alpha^i_k\|^2_{H^1} &= \tfrac{\tau}{4\epsilon} \|\dot \alpha^i_k\|^2_{H^1} + \tfrac{\tau}{4\epsilon} \|\dot \alpha^i_k\|^2_{H^1} \geq \tfrac{\tau}{4\epsilon} \|\dot \alpha^i_k\|^2_{L^2} + \tfrac{C\tau}{4\epsilon} \|\dot \alpha^i_k\|^2_{L^r} \geq \tfrac{\tau}{4\epsilon} \|\dot \alpha^i_k\|^2_{L^2} + \tfrac{C\tau^2}{4\epsilon^2} \|\dot \alpha^i_k\|^2_{L^r} \\
& \geq \tfrac{C \tau}{\epsilon} \Big[ \|\dot \alpha^i_k\|^2_{L^2} + c_1 \tfrac{\tau}{\epsilon} \|\dot \alpha^i_k\|^2_{L^r} \Big] = 2c_3 \tfrac{\tau}{\epsilon}A_i^2   \, .
\end{split}
\end{equation}
Collecting \eqref{eq:08111721}--\eqref{eq:09111539} and setting $\gamma := c_3 \tfrac{\tau}{\epsilon}$,  we obtain that
\begin{equation} \label{eq:before discrete Gronwall}
2 A_i(A_i - A_{i-1}) + 2 \gamma A_i^2 + B_i^2  \leq C_i^2 + 2 A_i D_i \, ,
\end{equation}
for every $i = 2, \dots, k$.

\underline{Proof of estimate \eqref{eq:estimate 1}}. Here we prove a slightly stronger inequality with an additional term on the left-hand side. Specifically, we show that
\begin{equation} \label{eq:stronger estimate}
\epsilon \Big[ \| \dot \alpha_k(t) \|_{L^2}^2 + c_1 \tfrac{\tau}{\epsilon} \| \dot \alpha_k(t) \|_{L^r} \Big] \leq C \exp\Big(\tfrac{C}{\epsilon} \ol \tau_k(t) \Big) \, .
\end{equation}
  By the inequalities $2A_i(A_i - A_{i-1}) \geq A_i^2 - A_{i-1}^2$ and $D_i \leq C \frac{\tau}{\epsilon} A_i$, from \eqref{eq:before discrete Gronwall} we get in particular that
\[
 A_i^2 -  A_{i-1}^2 + B_i^2 \leq C_i^2 + C \tfrac{\tau}{\epsilon} A_i^2 \, ,
\]
for $i = 2,\dots, k$. We fix $h \in \{ 2, \dots, k\}$ and we sum the inequality above for $i = 2, \dots, h$, deducing that
\begin{equation} \label{eq:estimate on Ah}
\epsilon A_h^2 - \epsilon A_1^2 + \sum_{i=2}^h B_i^2 \leq \epsilon \sum_{i=2}^h C_i^2 + C \sum_{i=2}^h  \tau A_i^2  \, .
\end{equation}

We claim that
\begin{equation} \label{eq:estimate on A1}
\epsilon A_1^2 \leq C \Big[\epsilon C_1^2 + \tau A_1^2 +  \tfrac{1}{\epsilon}   \Big]  \, .
\end{equation}
Once \eqref{eq:estimate on A1} is proven, summing \eqref{eq:estimate on Ah} and \eqref{eq:estimate on A1}, by the initial assumption on $w$~\eqref{ass:w} we conclude that
\begin{equation} \label{eq:22031225}
\epsilon A_h^2  \leq C \Big[ \tfrac{1}{\epsilon} + \sum_{i=1}^h \tau \big( \| \dot w^i_k \|^2_{W^{1,\tilde p}} +  \| \dot w^{i-1}_k \|^2_{W^{1,\tilde p}} \big)   +  \sum_{i=1}^h \tau A_i^2 \Big] \leq C \Big[ 1 + \tfrac{1}{\epsilon} + \sum_{i=1}^h \tau A_i^2 \Big]
\end{equation}
for every $h = 1, \dots, k$. By a discrete Gronwall inequality on $\epsilon A_h^2$ we deduce that
\begin{equation} \label{eq:after discrete gronwall}
\epsilon A_h^2 \leq C \Big( 1 + \tfrac{1}{\epsilon} \Big) \exp \Big( C {\tfrac{t^h_k}{\epsilon} } \Big)
\end{equation}
for every $h = 1, \dots, k$. Multiplying by $\epsilon$ and taking the square root, we get
\begin{equation} \label{eq:10011313}
\epsilon A_h \leq C  \exp \Big( C {\tfrac{t^h_k}{\epsilon}} \Big)
\end{equation}
and thus \eqref{eq:stronger estimate}.

It remains to prove~\eqref{eq:estimate on A1}. Adding and subtracting $\de_\alpha \E(\alpha_0, u_0)$ to \eqref{eq:eul2} evaluated at time $t \in (0, t^1_k)$, we deduce that
\[
\langle \de_\alpha \E(\alpha^1_k, u^1_k) - \de_\alpha \E(\alpha_0, u_0), \dot \alpha^1_k \rangle + \R(\dot \alpha^1_k, V_0) + \langle \de_\alpha \E(\alpha_0, u_0), \dot \alpha^1_k \rangle + \epsilon \| \dot \alpha^1_k \|^2_{L^2}  = 0 \, .
\]
With computations similar to those previously done in \eqref{eq:05121748}--\eqref{eq:08111609} and using the assumption $\de_\alpha \E(\alpha_0, u_0) \in L^2(\Omega)$, we infer that
\[
\begin{split}
\epsilon  \| \dot \alpha^1_k \|^2_{L^2} + \tau \| \nabla \dot \alpha^1_k \|_{L^2}^2 &\leq C \tau \Big[ \big\| \nabla u^1_k + \nabla u_0 \|_{L^2} \big\|\nabla \dot u^1_k \|_{L^p} \| \dot \alpha^1_k\|_{L^{q_1}} + \big\| \nabla u^1_k \big\|^2_{L^p}  \| \dot \alpha^1_k\|^2_{L^{q_1}} \Big] \\
 & \quad + \| f \|_{L^\infty} \| \dot \alpha^1_k\|_{L^1} + \|\de_\alpha \E(\alpha_0,u_0) \|_{L^2} \|\dot \alpha^1_k\|_{L^2} \\
 & \leq C \tau \Big[ \| \dot w^1_k \|^2_{W^{1,\tilde p}} + \| \dot \alpha^1_k\|^2_{L^{r} }\Big] + \tfrac{\epsilon}{2} \| \dot \alpha^1_k\|^2_{L^2}  + \tfrac{C}{\epsilon}\, .
\end{split}
\]
Using inequality \eqref{eq:interpolation inequality} as above, it is not difficult to see that
\begin{equation} \label{eq:estimate on alpha1}
\epsilon \| \dot \alpha^1_k \|^2_{L^2} + \tau \|\dot \alpha^1_k \|^2_{H^1} \leq C \tau \Big[ \| \dot w^1_k \|^2_{W^{1,\tilde p}}  + \| \dot \alpha^1_k\|^2_{L^1 } + \tfrac{1}{\epsilon} \Big] \, ,
\end{equation}
which in turn implies \eqref{eq:estimate on A1}.

\underline{Proof of estimate \eqref{eq:estimate 2}}. Inequalities \eqref{eq:estimate on Ah} and \eqref{eq:estimate on A1} imply in particular that
\begin{equation} \label{eq:10011158}
 \sum_{i=2}^h B_i^2 \leq \epsilon \sum_{i=2}^h C_i^2 + C \sum_{i=2}^h  \tau A_i^2 + C \Big[\epsilon C_1^2 + \tau A_1^2 +  \tfrac{1}{\epsilon}   \Big] \, .
\end{equation}
From~\eqref{eq:estimate on alpha1} we deduce that
\begin{equation} \label{eq:10011159}
\epsilon B_1^2 \leq C \Big[ \epsilon C_1^2 + \tau A_1^2 + \tfrac{1}{\epsilon} \Big]
\end{equation}
Let us fix $h  \in \{ 1, \dots, k\}$. Summing \eqref{eq:10011158} and \eqref{eq:10011159}, by~\eqref{ass:w} we obtain that
\[
\epsilon \sum_{i = 1}^h B_i^2 \leq C \Big[ \tfrac{1}{\epsilon} + \sum_{i=1}^h \tau \big( \| \dot w^i_k \|^2_{W^{1,\tilde p}} +  \| \dot w^{i-1}_k \|^2_{W^{1,\tilde p}} \big) + \sum_{i=1}^h \tau A_i^2 \Big]  \leq C \Big[1 + \tfrac{1}{\epsilon} + \sum_{i=1}^h \tau A_i^2  \Big]
\]
and thus, multiplying by $\epsilon$ and using \eqref{eq:after discrete gronwall},
\[
\epsilon \sum_{i=1}^h \tau \| \dot \alpha^i_k \|^2_{H^1} \leq C \exp \Big( C \tfrac{t^h_k}{\epsilon}  \Big) \, .
\]
In the equality above we have integrated the exponential function in time and we have used the fact that $\tau << \epsilon$. This concludes the proof of \eqref{eq:estimate 2}.

\underline{Proof of estimate~\eqref{eq:estimate 3}}. By the discrete Gronwall estimate proved in \cite[Lemma 4.1]{KRZ13a} we deduce that for every $h = 2,\dots, k$
\begin{equation} \label{eq:10111036}
\begin{split}
\Big(\sum_{i=2}^h (1+\gamma)^{2(i-h)-1} B_i^2 \Big)^{\! \frac{1}{2}} & \leq \Big( (1+\gamma)^{-2h} A_1^2 + \sum_{i=2}^h (1+\gamma)^{2(i-h)-1} C_i^2 \Big)^{\! \frac{1}{2}} + \sqrt{2} \sum_{i=2}^h (1+\gamma)^{i-k-1}D_i  \\
& \leq \bigg[ 2(1+\gamma)^{-2h} A_1^2 + 2\sum_{i=2}^h (1+\gamma)^{2(i-h)-1} C_i^2  + 4 \Big( \sum_{i=2}^h (1+\gamma)^{i-k-1}D_i  \Big)^{\!2} \bigg]^{\! \frac{1}{2}} \\
& \leq  \sqrt{2}(1+\gamma)^{-h} A_1 + 2\sum_{i=2}^h (1+\gamma)^{2(i-h)-1} C_i^2 + 1 + 2  \sum_{i=2}^h (1+\gamma)^{i-k-1}D_i
\end{split}
\end{equation}

Using the estimate 
\[
\gamma \sum_{i=2}^h (1+\gamma)^{2(i-h)-1} = \gamma (1+\gamma)^{-2h -1} \frac{(1+\gamma)^4 - (1+\gamma)^{2h+2}}{1-(1+\gamma)^2} = \frac{1+\gamma}{2+\gamma} \big[ 1 - (1+\gamma)^{-2h+2} \big] \leq 1 \, ,
\]
by the Cauchy-Schwarz inequality we estimate the left-hand side of \eqref{eq:10111036} by
\[
\gamma \sum_{i=2}^h (1+\gamma)^{2(i-h)-1} \|\dot \alpha^i_k \|_{H^1} =  \sum_{i=2}^h \big( \gamma (1+\gamma)^{2(i-h)-1} \big)^{\!\frac{1}{2}}\big( \gamma (1+\gamma)^{2(i-h)-1} \big)^{\!\frac{1}{2}} \|\dot \alpha^i_k \|_{H^1} \leq C \Big(\sum_{i=2}^h (1+\gamma)^{2(i-h)-1} B_i^2 \Big)^{\! \frac{1}{2}} ,
\] 
for $h =2, \dots, k$. Hence \eqref{eq:10111036} reads
\begin{equation} \label{eq:10111118}
\begin{split}
\gamma \sum_{i=2}^h & (1+\gamma)^{2(i-h)-1} \|\dot \alpha^i_k \|_{H^1} \\
&\leq C \Big[  1 + (1+\gamma)^{-h} A_1  + \gamma \sum_{i=2}^h (1+\gamma)^{2(i-h)-1} \big( \| \dot w^i_k \|^2_{W^{1,\tilde p}} +  \| \dot w^{i-1}_k \|^2_{W^{1,\tilde p}} \big)  +  \gamma \sum_{i=2}^h (1+\gamma)^{i-k-1}\| \dot \alpha^i_k \|_{L^1}   \Big] \, ,
\end{split}
\end{equation}
for $h = 2,\dots, k$. We multiply both sides of \eqref{eq:10111118} by $\tau$ and we sum over $h = 2,\dots, k$. Using the expression of the partial sums of the geometric series, it is possible to show that
\begin{equation} \label{eq:10011307}
\sum_{i=2}^k \tau \|\dot \alpha^i_k \|_{H^1} \leq C \Big[ 1 + \epsilon A_1 + \sum_{i=2}^k \tau \big( \| \dot w^i_k \|^2_{W^{1,\tilde p}} +  \| \dot w^{i-1}_k \|^2_{W^{1,\tilde p}} \big) +   \sum_{i=2}^k \tau \| \dot \alpha^i_k \|_{L^1}   \Big]  \, .
\end{equation}
We refer to \cite[Proposition 4.3]{KRZ13a} or \cite[Proposition 3.8]{Cri16} for more details about the computations mentioned above.

Multiplying \eqref{eq:estimate on alpha1} by $\tau$, taking the square root and using the fact that $\tau << \epsilon$, we infer that
\[
\tau \|\dot \alpha^1_k \|_{H^1} \leq C \tau \Big[ \| \dot w^1_k \|_{W^{1,\tilde p}}  + \| \dot \alpha^1_k\|_{L^1 } + 1 \Big] \, .
\]
Adding this last inequality to \eqref{eq:10011307} we obtain that
\[
\sum_{i=1}^k \tau \|\dot \alpha^i_k \|_{H^1} \leq C \Big[ 1 + \epsilon A_1 + \sum_{i=1}^k \tau \big( \| \dot w^i_k \|^2_{W^{1,\tilde p}} +  \| \dot w^{i-1}_k \|^2_{W^{1,\tilde p}} \big) +   \sum_{i=1}^k \tau \| \dot \alpha^i_k \|_{L^1}   \Big]  \, .
\]
To conclude the proof of \eqref{eq:estimate 3}, we observe that: $\epsilon A_1 \leq C$ by \eqref{eq:10011313} evaluated for $h=1$; the second sum is bounded by a constant by the initial assumption on $w$ \eqref{ass:w}; the third sum is actually a telescopic sum, namely
\[
 \sum_{i=1}^k \tau \| \dot \alpha^i_k \|_{L^1} = \integral{\Omega}{\big( \alpha_0 - \alpha^k_k \big) }{\d x} \leq |\Omega| \, .
\]
\end{proof}

In order to obtain the energy dissipation balance for the evolution $(\alpha_k, u_k)$, in Proposition~\ref{prop:discrete balance}, we integrate in time the energy evaluated on these affine interpolations. We are allowed to do so because they are absolutely continuous (actually $H^1$) in time. Since we also employ the Euler equation \eqref{eq:eul2} of Lemma~\ref{lemma:euler}, that contains also the piecewise constant interpolations, we have to estimate the difference of the piecewise affine and constant interpolations. This is done in the following remark.
\begin{remark}
For every $t \in [0,T]$
\[
\| \alpha_k(t) - \ol \alpha_k(t) \|_{H^1} = \Big\| \int_{t}^{\ol t_k(t)} \dot \alpha_k(s) \d s \Big\|_{H^1}  \leq  \int_{t}^{\ol t_k(t)} \| \dot \alpha_k(s) \|_{H^1} \d s \leq \tau^\frac{1}{2} \| \alpha_k \|_{H^1(0,T; H^1(\Omega))}
\]
and therefore, by~\eqref{eq:estimate 2},
\begin{subequations}\label{3105181043}
\begin{equation} \label{eq:affine and piecewise same alpha ol}
\| \alpha_k - \ol \alpha_k \|_{L^\infty(0,T;H^1(\Omega))} \leq C_\epsilon \tau^\frac{1}{2} \, .
\end{equation}
Similarly, we have
\begin{align}
\| \alpha_k - \ul \alpha_k \|_{L^\infty(0,T;H^1(\Omega))} & \leq C_\epsilon \tau^\frac{1}{2} \, , \label{eq:affine and piecewise same alpha ul} \\
\| u_k - \ol u_k \|_{L^\infty(0,T;W^{1,p}(\Omega))} & \leq C_\epsilon \tau^\frac{1}{2} \, , \quad \text{for } p \in [2, \tilde p) \, , \label{eq:affine and piecewise same u ol}\\
\| u_k - \ul u_k \|_{L^\infty(0,T;W^{1,p}(\Omega))} & \leq C_\epsilon \tau^\frac{1}{2} \, , \quad \text{for } p \in [2, \tilde p) \, . \label{eq:affine and piecewise same u ul}
\end{align}
\end{subequations}
\end{remark}

\subsection*{Discrete energy-dissipation balance}  Here we obtain the energy-dissipation balance, by employing the Euler condition \eqref{eq:eul2}, correcting with the piecewise affine interpolations in place of the piecewise constant ones.

\begin{proposition}[Discrete energy-dissipation balance] \label{prop:discrete balance}
\begin{equation} \label{eq:discrete balance}
\begin{split}
\E(\alpha_k(T), u_k(T)) & + \int_0^T \! \R(\dot \alpha_k(t); \V_k(t)) \, \d t + \epsilon \int_0^T \! \|\dot \alpha_k(t)\|^2_{L^2} \d t  \\
 & = \E(\alpha_0, u_0) + \int_0^T \! \langle \mu(\ol \alpha_k(t)) \nabla \ol u_k(t) , \nabla \dot w_k(t) \rangle_{L^2}  \, \d t + R_k \, ,
\end{split}
\end{equation}
where $R_k \to 0$ as  $k \to +\infty$.
\end{proposition}
\begin{proof}
By \eqref{eq:estimate 2}--\eqref{eq:estimate 3}, the piecewise affine interpolations $\alpha_k(t)$ and $u_k(t)$ are absolutely continuous in $t$. As a consequence, $t \mapsto \E(\alpha_k(t), u_k(t))$ is absolutely continuous and
\begin{equation} \label{eq:tderivative of E}
\begin{split}
\frac{\d}{\d t} \Big[ \E(\alpha_k(t), u_k(t)) \Big] & = \langle \de_\alpha \E(\alpha_k(t), u_k(t)) , \dot \alpha_k(t) \rangle + \langle \de_u \E(\alpha_k(t), u_k(t)) , \dot u_k(t) \rangle  \\
& = \langle \de_\alpha \E(\ol \alpha_k(t), \ol u_k(t)) , \dot \alpha_k(t) \rangle + \langle \de_u \E(\ol \alpha_k(t), \ol u_k(t)) , \dot u_k(t) \rangle  + \eta_k(t)
\end{split}
\end{equation}
for a.e.\ $t$, where
\begin{equation} \label{eq:def of eta}
\eta_k(t) := \langle \de_\alpha \E(\alpha_k(t), u_k(t)) - \de_\alpha \E(\ol \alpha_k(t), \ol u_k(t)) , \dot \alpha_k(t) \rangle + \langle \de_u \E(\alpha_k(t), u_k(t)) - \de_u \E(\ol \alpha_k(t), \ol u_k(t)) , \dot u_k(t) \rangle \, .
\end{equation}
Using $\dot u_k(t) - \dot w_k(t)$ as test function in \eqref{eq:pb of uik}, we deduce that
\[
\langle \de_u \E(\ol \alpha_k(t), \ol u_k(t)) , \dot u_k(t) \rangle = \langle \de_u \E(\ol \alpha_k(t), \ol u_k(t)) , \dot w_k(t) \rangle = \langle \mu(\ol \alpha_k(t)) \nabla \ol u_k(t) , \nabla \dot w_k(t) \rangle_{L^2}
\]
Together with the Euler equation for $\ol \alpha_k(t)$~\eqref{eq:eul2} and~\eqref{eq:tderivative of E}, this gives
\begin{equation}\label{3105182138}
\frac{\d}{\d t} \Big[ \E(\alpha_k(t), u_k(t)) \Big] = - \epsilon \|\dot \alpha_k(t)\|^2_{L^2} - \R(\dot \alpha_k(t); \V_k(t)) + \langle \mu(\ol \alpha_k(t)) \nabla \ol u_k(t) , \nabla \dot w_k(t) \rangle_{L^2}  + \eta_k(t)
\end{equation}
Integrating in time the previous equality, we obtain~\eqref{eq:discrete balance} with $R_k := \int_0^T \! \eta_k(t) \d t$.

Let us show that $R_k \to 0$. By H\"older's Inequality, by~\eqref{eq:affine and piecewise same alpha ol}, by~\eqref{eq:estimate of u}, and by \eqref{eq:estimate 3} we deduce that
\[
\begin{split}
\Big| \int_0^T \! &  \integral{\Omega}{\big( \mu'(\alpha_k(t)) - \mu'(\ol \alpha_k(t)) \big) |\nabla u_k(t)|^2 \dot \alpha_k(t)}{\d x} \, \d t \Big| \leq C \int_0^T \! \integral{\Omega}{|\alpha_k(t) - \ol \alpha_k(t)| |\nabla u_k(t) |^2 |\dot \alpha_k(t)|}{ \d x} \, \d t \\
& \leq C \int_0^T \! \|\alpha_k(t) - \ol \alpha_k(t) \|_{H^1} \| u_k(t) \|_{W^{1,p}}^2 \| \dot \alpha_k(t) \|_{H^1} \d t  \leq C_\epsilon \tau^\frac{1}{2}
\end{split}
\]
Furthermore by H\"older's Inequality, by~\eqref{eq:estimate of u}, by~\eqref{eq:affine and piecewise same u ol}, and by \eqref{eq:estimate 3} we infer that
\[
\begin{split}
\Big| \int_0^T \! & \integral{\Omega}{\mu'(\ol \alpha_k(t)) \big( |\nabla u_k(t)|^2 - |\nabla \ol u_k(t)|^2 \big) \dot \alpha_k(t) }{\d x} \, \d t \Big| \leq C \int_0^T  \! \integral{\Omega}{ | \nabla u_k(t) + \nabla \ol u_k(t) | \, | \nabla u_k(t) - \nabla \ol u_k(t) | \, | \dot \alpha_k(t) |  }{\d x} \, \d t \\
& \leq C  \int_0^T \! \|  u_k(t) +  \ol u_k(t) \|_{W^{1,p}} \| u_k(t) - \ol u_k(t) \|_{W^{1,p}} \| \dot \alpha_k(t) \|_{H^1} \d t  \leq C_\epsilon \tau^\frac{1}{2}  \, .
\end{split}
\]
Finally, by~\eqref{eq:affine and piecewise same alpha ol} and \eqref{eq:estimate 3} we get that
\[
\Big| \int_0^T \! \integral{\Omega}{\big( \nabla \alpha_k(t) - \nabla \ol \alpha_k(t) \big) \cdot \nabla \dot \alpha_k(t) }{\d x} \, \d t  \Big| \leq  \int_0^T \! \| \alpha_k(t) - \ol \alpha_k(t) \|_{H^1} \| \dot \alpha_k(t) \|_{H^1} \, \d t \leq C_\epsilon \tau^\frac{1}{2} \, .
\]
This shows that
\[
\lim_{k \to + \infty} \int_0^T \!  \langle \de_\alpha \E(\alpha_k(t), u_k(t)) - \de_\alpha \E(\ol \alpha_k(t), \ol u_k(t)) , \dot \alpha_k(t) \rangle \, \d t = 0 \, .
\]
With completely analogous computations it is not difficult to show that
\[
\lim_{k \to + \infty} \int_0^T \!  \langle \de_u \E(\alpha_k(t), u_k(t)) - \de_u \E(\ol \alpha_k(t), \ol u_k(t)) , \dot u_k(t) \rangle \, \d t = 0 \, .
\]
This concludes the proof.
\end{proof}
We observe that the energy balance \eqref{eq:discrete balance} holds for any couple of times $t_1<t_2 \in [0,T]$, as one can see arguing as in Proposition~\ref{prop:discrete balance} and integrating \eqref{3105182138} in the time interval $(t_1,t_2)$.

\section{Existence of viscous evolutions}\label{Sec:3}

In this section we pass to the limit as $k \to +\infty$ (i.e., as the time-step goes to zero). Notice that $\epsilon > 0$ is fixed in this section. The main result is the existence of \emph{viscous evolutions}, defined as follows.
Given $\alpha_\epsilon\in AC([0,T]; H^1(\Omega))$, $u_\epsilon \in AC([0,T]; H^1(\Omega))$ we define, as in \eqref{3005181820}, $\zeta_\varepsilon:=g(\alpha_\varepsilon)\nabla u_\varepsilon$ and, as in \eqref{3005181833},
\begin{equation}\label{eq:def of Veps}
V_\varepsilon(t):=\int_0^t \big| \dot{\zeta}_\varepsilon(s)\big| \d s\,,
\end{equation}
as a Bochner integral in $L^2(\Omega)$.  During the section we are in the constitutive assumptions of Section~\ref{sec:Ass}.

\begin{definition}\label{def:appreveps} We say that a function $(\alpha_\epsilon, u_\epsilon) \colon [0,T] \to H^1(\Omega) {\times} W^{1,p}(\Omega)$ is an $\epsilon$\emph{-approximate viscous evolution} if $\alpha_\epsilon\in H^1(0,T; H^1(\Omega))$, $u_\epsilon \in H^1(0,T; W^{1,p}(\Omega))$ and the following conditions are satisfied:
\begin{itemize}
\item[{\rm (ev0)$_\epsilon$}]{\em irreversibility\/}:
\begin{equation*} \label{irep}
 [0,T] \ni t \mapsto \alpha_\epsilon(t) \quad \text{is nonincreasing as a family of measurable functions on }\Omega\,,
 \end{equation*}
  that is $\alpha_\epsilon(t) \leq \alpha_\epsilon(s)$ a.e.\ in $\Omega$ for all $s \leq t$; 
 \item[{\rm (ev1)$_\epsilon$}] {\em equilibrium\/}: for every $t \in [0,T]$,  $u_\epsilon(t) \in H^1(\Omega)$ is a weak solution to the problem 
\begin{equation} \label{2202180026}
\left\{
\begin{aligned}
\div\big(\mu(\alpha_\epsilon(t)) \nabla u_\epsilon(t) \big) &= 0  && \text{in } \Omega \, , \\
u_\epsilon(t) & = w(t) && \text{on } \de_D \Omega \, .
\end{aligned}
\right.
\end{equation}
\item[{\rm (ev2)$_\epsilon$}] {\em  Karush-Kuhn-Tucker  inequality\/}: for a.e.\ $t \in (0,T)$ and for every $\beta \in H^1(\Omega)$ with $\beta \leq 0$ a.e.\ in $\Omega$ we have
\begin{equation} \label{2202180030}
\langle \de_\alpha \E(\alpha_\epsilon(t), u_\epsilon(t)), \beta \rangle +\R(\beta; V_\epsilon(t)) + \epsilon \langle \dot \alpha_\epsilon(t), \beta \rangle_{L^2} \geq 0 \, .
\end{equation}
\item[{\rm (ev3)$_\epsilon$}] {\em energy balance\/}:
\[
\E(\alpha_\epsilon(T), u_\epsilon(T)) + \int_0^T \! \R(\dot{\alpha}_\epsilon(t); V_\epsilon(t))\, \d t  + \epsilon \int_0^T \! \| \dot \alpha_\epsilon(t) \|^2_{L^2} \, \d t = \E(\alpha_0, u_0) + \int_0^T \! \langle \mu(\alpha_\epsilon(t)) \nabla u_\epsilon(t) , \nabla \dot w(t) \rangle_{L^2}  \, \d t  \, .
\]
\end{itemize}
\end{definition}

All the section is devoted to the proof of the result below.
\begin{theorem}\label{teo:exViscEvo}
Let $\tilde{p}>2$ be given by Lemma~\ref{lemma:increment of u}.
For every $\epsilon>0$  and $p< \tilde{p}$
there exists an $\epsilon$\emph{-approximate viscous evolution} $(\alpha_\epsilon,u_\epsilon)$ with $(\alpha_\epsilon(0), u_\epsilon(0))=(\alpha_0, u_0)$ and there is a constant $C>0$, independent of $\epsilon$, such that
\begin{equation} \label{2202180044}
\int_0^T \! \| \dot \alpha_\epsilon(s) \|_{H^1} \, \d s +  \int_0^{T} \! \| \dot u_\epsilon(s) \|_{W^{1,p}} \, \d s   \leq C \, .
\end{equation}
\end{theorem}

The strategy of the proof consists in showing first the existence of a weak form of $\epsilon$-approximate viscous evolution. This satisfies the conditions (ev0)$_\epsilon$, (ev1)$_\epsilon$, and the (ev2)$_\epsilon$, (ev3)$_\epsilon$ with a different expression of dissipation (Propositions~\ref{prop:e-d inequality}, \ref{prop:stweak}, and \ref{prop:balance weak}). Such a weak existence result allows us to improve, for fixed $\epsilon$, the \emph{a priori} convergences of the discrete-time evolutions (Proposition~\ref{prop:strconveps}) and to express the dissipation in terms of $V_\epsilon(t)$, the cumulation of $\zeta_\varepsilon$ (cf.\ \eqref{eq:def of Veps}), so recovering its form in Definition~\ref{def:appreveps}, by Lemma~\ref{lemma:desired convergence}. 

\subsection*{Compactness}
We start by exploiting the {\em a priori} bounds found in Proposition~\ref{prop:enhanced estimates} to deduce compactness of the discrete-time evolutions.
By~\eqref{eq:estimate 2} we find a subsequence (which we do not relabel) such that
\begin{align}
\alpha_k & \weak \alpha_\epsilon  \quad \text{weakly in } H^1(0,T;H^1(\Omega))  \, , \label{eq:weak alpha}\\
u_k & \weak u_\epsilon \quad \text{weakly in } H^1(0,T; W^{1,p}(\Omega)) \, , \text{ for } p \in [2, \tilde p \, ) \, , \label{eq:weak u}
\end{align}
as $k \to +\infty$. (Actually, we also extract a subsequence independent of $t$ such that the convergence in~\eqref{eq:convergence of f} below holds. We do not state this here for the sake of clarity in the presentation.) By the compact embeddings $H^1(\Omega) \compact L^q(\Omega)$ and $W^{1,p}(\Omega) \compact L^p(\Omega)$, by the Aubin-Lions lemma \cite{Aub63}, and by \eqref{3105181043}  we deduce that
\begin{subequations}\label{3105181105}
\begin{align}
\| \alpha_k - \alpha_\epsilon \|_{C([0,T];L^q(\Omega))}, \ \| \ol \alpha_k - \alpha_\epsilon \|_{L^\infty(0,T;L^q(\Omega))}, \ \| \ul \alpha_k - \alpha_\epsilon \|_{L^\infty(0,T;L^q(\Omega))} & \to 0 \, , \text{ for } q \in [1, \infty) \, , \label{eq:uniform convergence of alpha}\\
\| u_k - u_\epsilon \|_{C([0,T];L^p(\Omega))}, \ \| \ol u_k - u_\epsilon \|_{L^\infty(0,T;L^p(\Omega))}, \ \| \ul u_k - u_\epsilon \|_{L^\infty(0,T;L^p(\Omega))} & \to 0 \, , \text{ for } p \in [2, \tilde p \, ) \, . \label{eq:uniform convergence of u}
\end{align}
\end{subequations}
Moreover, from the inequality $\| \alpha_k \|_{L^\infty(0,T;H^1(\Omega))}  \leq C \big( 1 + \| \alpha_k \|_{W^{1,1}(0,T;H^1(\Omega))}\big)$ and by~\eqref{eq:estimate of u} and \eqref{eq:affine and piecewise same alpha ol}--\eqref{eq:affine and piecewise same u ul} we deduce that for every $t \in [0,T]$ we also have
\begin{align}
\alpha_k(t) \, , \ \ol \alpha_k(t) \, , \ \ul \alpha_k(t) & \weak \alpha_\epsilon(t)  \quad \text{weakly in } H^1(\Omega)  \, , \label{eq:pointwise weak alpha}\\
u_k(t) \, , \ \ol u_k(t) \, , \ \ul u_k(t) & \weak u_\epsilon(t) \quad \text{weakly in } W^{1,\tilde p}(\Omega) \, . \label{eq:pointwise weak u}
\end{align}
In particular,  for every $s \leq t$ we have $\alpha_\epsilon(t) \leq \alpha_\epsilon(s)$ a.e.\ in $\Omega$.  Moreover, for every $t \in [0,T]$ we have
\begin{equation} \label{eq:W1tilde norm}
\| u_\epsilon(t) \|_{W^{1, \tilde p}} \leq \liminf_{k \to +\infty} \| u_k(t) \|_{W^{1,\tilde p}} \leq C  \, .
\end{equation}
In view of the convergences \eqref{eq:weak alpha}, \eqref{eq:weak u},  by \eqref{eq:estimate 3} we get
\begin{equation} \label{eq:estimate 3 for alphaeps ueps}
\int_0^T \! \| \dot \alpha_\epsilon(s) \|_{H^1} \, \d s +  \int_0^{T} \! \| \dot u_\epsilon(s) \|_{W^{1,p}} \, \d s   \leq C \, ,
\end{equation}
 where $C$ is independent of $\epsilon$,  and then $V_\varepsilon$ is well defined as in \eqref{eq:def of Veps}.

\subsection*{Energy-dissipation balance and stability} In this subsection we pass to the limit as $k \to +\infty$ in the discrete energy-dissipation balance~\eqref{eq:discrete balance}. We start by discussing the easiest terms in the energy-dissipation balance, namely  the terms involving the energy, the viscous dissipation, and the work done by the boundary forces. The dissipation involving the fatigue term requires finer techniques and will be discussed below.

\vspace{1em}

From the pointwise convergences \eqref{eq:pointwise weak alpha}--\eqref{eq:pointwise weak u} and the lower semicontinuity of the energy $\E$ with respect to the weak convergence of $\alpha$ in $H^1(\Omega)$ and the weak convergence of $u$ in $W^{1,p}(\Omega)$ we deduce that
\begin{equation} \label{eq:E is lsc}
\E(\alpha_\epsilon(T), u_\epsilon(T) ) \leq \liminf_{k \to +\infty} \E(\alpha_k(T), u_k(T)) \, .
\end{equation}
Moreover, since $\dot \alpha_k \weak \dot \alpha_\epsilon$ weakly in $L^2(0,T; L^2(\Omega))$, we have that
\begin{equation} \label{eq:viscous diss is lsc}
\epsilon \int_0^T \! \|\dot \alpha_\epsilon(t)\|^2_{L^2} \d t \leq \liminf_{k \to +\infty} \Big(  \epsilon \int_0^T \! \|\dot \alpha_k(t)\|^2_{L^2} \d t \Big) \, .
\end{equation}
We claim that
\begin{equation} \label{eq:limit of work}
\lim_{k \to +\infty } \int_0^T \! \langle \mu(\ol \alpha_k(t)) \nabla \ol u_k(t) , \nabla \dot w_k(t) \rangle_{L^2}  \, \d t = \int_0^T \! \langle \mu(\alpha_\epsilon(t)) \nabla u_\epsilon(t) , \nabla \dot w(t) \rangle_{L^2}  \, \d t  \, .
\end{equation}
To show the convergence above, first of all we notice that $\mu(\ol \alpha_k(t)) \nabla \ol u_k(t) \weak \mu(\alpha_\epsilon(t)) \nabla u_\epsilon(t)$ weakly in $L^2(\Omega; \RR^2)$ for every $t \in [0,T]$ thanks to \eqref{eq:pointwise weak alpha}--\eqref{eq:pointwise weak u}. In addition, \eqref{eq:estimate of u} and assumption~\eqref{ass:w} imply
\[
\big| \langle \mu(\ol \alpha_k(t)) \nabla \ol u_k(t) , \nabla \dot w_k(t) \rangle_{L^2} \big| \leq C \int_0^T \! \| \ol u_k(t) \|_{H^{1}} \| \dot w_k(t) \|_{H^{1}} \d t \leq C \, .
\]
Since $\nabla \dot w_k(t) \to \nabla \dot w(t)$ strongly in $L^2(\Omega; \RR^2)$ for a.e.\ $t \in (0,T)$, by the Dominated Convergence Theorem the convergence in~\eqref{eq:limit of work} holds true.

We consider now the limit of the dissipation involving the fatigue term. We start with the following lemma, which shows that the affine interpolation of the cumulation is close to the piecewise constant interpolation.

\begin{lemma}\label{le:2202181003}
For every $k \in \NN$, $\epsilon>0$ we have that
\begin{equation}\label{eq:2202181006}
\|f(V_k)-f(\V_k)\|_{L^2(0,T;L^2(\Omega))} \leq  C \tau \Big(\|\alpha_k\|_{H^1(0,T;L^2(\Omega))}  + \|  u_k \|_{H^1(0,T;H^1(\Omega))}  \Big) \leq C_\epsilon\,\tau\,.
\end{equation}
\end{lemma}
\begin{proof}
By \eqref{3005181946} and \eqref{eq:estimate of u} we have
\begin{equation*}
\begin{split}
V_k(t)-\V_k(t)&= \frac{t-\ul t_k(t)}{\tau}\bigg(\big[g(\alpha_k(t))-g(\ul\alpha_k(t)) \big] \nabla u_k(t) + g(\ul \alpha_k(t))\big[ \nabla u_k(t) - \nabla \ul u_k(t) \big] \bigg) \,,
\end{split}
\end{equation*}
so that
\begin{equation*}
\begin{split}
|V_k(t)-\V_k(t)|& \leq \tau \Big( \|g'\|_{L^\infty} |\dot{\alpha}_k(t)| |\nabla u_k(t)| + |g(\ul \alpha_k(t))| |\nabla \dot{u}_k(t)| \Big) \\
& \leq C \tau \big( |\dot{\alpha}_k(t)| + |\nabla \dot{u}_k(t)| \big)\,.
\end{split}
\end{equation*}
Thus for any $\beta\in L^2(0,T;L^2(\Omega))$ (recall that $f$ is Lipschitz)
\[
\begin{split}
 \int_0^T \! \integral{\Omega}{|f(V_k(t)) - f(\V_k(t))| \, | \beta(t)|}{\d x}  \, \d t & \leq C \tau \int_0^T \! \integral{\Omega}{\big( |\dot{\alpha}_k(t)| + |\nabla \dot{u}_k(t)| \big) \, |\beta(t)|}{\d x}  \, \d t \\
&  \leq C \tau \int_0^T \! \Big( \|\dot{\alpha}_k(t)\|_{L^2} + \| \dot u_k(t) \|_{H^{1}} \Big) \|\beta(t)\|_{L^2} \, \d t  \\
& \leq C \tau \Big( \|\alpha_k\|_{H^1(0,T; L^2(\Omega))} + \|  u_k \|_{H^1(0,T;H^1(\Omega)) } \Big)  \|\beta \|_{L^2(0,T;L^2(\Omega))} \, .
\end{split}
\]
Recalling \eqref{eq:estimate 2}, the estimate above gives \eqref{eq:2202181006}. We notice that we arrive at the same conclusion also with $\zeta$ defined by $g(\alpha)|\nabla u|^\theta$, for $\theta\in [1,\tilde{p})$, arguing similarly to what done to pass from \eqref{0106180946} to \eqref{0106180947} in Proposition~\ref{prop:enhanced estimates}.
\end{proof}
\vspace{1em}

In the following lemma we show that a strong convergence of the discrete-time evolutions would guarantee the convergence  of the dissipation term. We stress that the a priori bounds on $u_k(t)$ found in Proposition~\ref{prop:enhanced estimates} only guarantee the weak convergence~\eqref{eq:weak u}. Therefore we are not allowed to apply Lemma~\ref{lemma:desired convergence} at the moment.

\begin{lemma} \label{lemma:desired convergence}
Assume that the following convergences for $\alpha_k$ and $u_k$ hold true:
\begin{subequations}\label{3105181042}
\begin{align}
\alpha_k & \to \alpha_\epsilon \quad  \text{strongly in } W^{1,1}(0,T; L^2(\Omega)) \,, \label{3105180858} \\
u_k & \to u_\epsilon \quad  \text{strongly in } W^{1,1}(0,T; W^{1,p}(\Omega)) \, , \quad \text{for } p \in [2, \tilde p) \, . \label{eq:strong assumption on u}
\end{align}
\end{subequations}
Then
\begin{equation} \label{eq:convergence of cumulation}
f (V_k) \to f (V_\epsilon) \quad \text{strongly in } L^2(0,T; L^2(\Omega)) \, ,
\end{equation}
and
\begin{equation} \label{eq:strong limit of dissipation}
\lim_{k \to +\infty }\int_0^T \! \R(\dot \alpha_k(t); \V_k(t)) \, \d t = \int_0^T \! \R(\dot \alpha_\epsilon(t); V_\epsilon(t)) \, \d t  \, ,
\end{equation}
where the cumulations $V_k$ and $V_\epsilon$ are defined in \eqref{3005181946} and \eqref{eq:def of Veps}, respectively.
\end{lemma}
\begin{proof}
For the proof it is convenient to introduce the function
\[
\d g (\beta , h) := \begin{cases}
\displaystyle \frac{g(\beta + h)  - g(\beta)}{h} \, , & \text{if } h \neq 0 \, , \\
g'(\beta) \, , & \text{if } h = 0 \, ,
\end{cases}
\]
for every $\beta, h \in \RR$. Observe that $g(\beta+h) = g(\beta) + h \d g(\beta, h)$ and since $g \in C^{1,1}(\RR)$
\[
|\d g(\beta, h) - g'(\beta)| \leq \| g'' \|_{L^\infty} |h| \, .
\]
Using the function $\d g$, we can write for every $s \in [0,T]$
\[
\dot \zeta_k(s) = \d g(\ul \alpha_k(s), \tau \dot \alpha_k(s)) \dot \alpha_k(s) \nabla \ol u_k(s) + g(\ul \alpha_k(s)) \nabla \dot u_k(s) \, .
\]
We now estimate $V_k-V_\varepsilon$ by employing~\eqref{3105180902} and~\eqref{eq:def of Veps}. For every $t \in [0,T]$ we have
\[\begin{split}
& \integral{\Omega}{| V_k(t;x) - V_\epsilon(t;x) | }{\d x} \leq \int_0^t \! \integral{\Omega}{\big| \dot \zeta_k(s;x)  -  \dot \zeta_\epsilon(s;x)  \big|}{\d x}  \, \d s \\
&
\leq \int_0^t \int\limits_\Omega \Big| \d g(\ul \alpha_k(s), \tau \dot \alpha_k(s)) \,  \dot{\alpha}_k(s) \,\nabla \ol u_k(s) +  g(\ul \alpha_k(s)) \,\nabla \dot{u}_k(s) - g'(\alpha_\varepsilon(s)) \,\dot{\alpha}_\varepsilon(s) \,\nabla u_\varepsilon(s) - g(\alpha_\varepsilon(s))\,\nabla \dot{u}_\varepsilon(s)  \Big| \d x \d s\\
&
\leq \int_0^t \int\limits_\Omega \bigg[ \tau \|g''\|_{L^\infty} \big|\dot \alpha_k(s)\big|^2 \big|\nabla \ol u_k(s)\big| + \|g''\|_{L^\infty} \big|\ul \alpha_k(s) - \alpha_\epsilon(s)\big| \big|\dot \alpha_k(s)\big| \big|\nabla \ol u_k(s)\big| \\
& \hspace{4em}
+\|g'\|_{L^\infty} \big| \dot \alpha_k(s) - \dot \alpha_\epsilon(s) \big| \big| \nabla \ol u_k(s) \big| + \|g'\|_{L^\infty} \big| \dot \alpha_\epsilon(s) \big| \big| \nabla \ol u_k(s) - \nabla u_\epsilon(s) \big| \\
& \hspace{4em}
+ \|g'\|_{L^\infty} \big| \ul \alpha_k(s) - \alpha_\epsilon(s) \big| \big| \nabla \dot u_k(s) \big| + \|g\|_{L^\infty} \big| \nabla \dot u_k(s) - \nabla \dot u_\epsilon(s) \big|  \bigg] \d x \d s
\\
&
\leq C \big( \| \ul \alpha_k - \alpha_\varepsilon \|_{L^\infty(0,T; L^q(\Omega))} \|\alpha_k\|_{W^{1,1}(0,T; L^q(\Omega))} + \|\alpha_k -\alpha_\varepsilon\|_{W^{1,1}(0,T; L^2(\Omega))}  \big) \| \ol u_k \|_{L^\infty(0,T; W^{1,p}(\Omega))} \\
&
\hspace{3em} + C \big( \tau + \| u_k - u_\varepsilon \|_{L^\infty(0,T;W^{1,p}(\Omega))} \big) \|\alpha_\varepsilon \|_{W^{1,1}(0,T; L^q(\Omega))}
\\
&
\hspace{3em} +  C \| \ul \alpha_k - \alpha_\varepsilon \|_{L^\infty(0,T; L^q(\Omega))} \| u_k\|_{W^{1,1}(0,T; W^{1,p}(\Omega))} + C \|   u_k  -    u_\epsilon  \|_{W^{1,1}(0,T;W^{1,p}(\Omega))} \, ,
\end{split}
\]
for  $q\in (2,\infty)$ such that $\frac{1}{q}+\frac{1}{p}<\frac{1}{2}$.

 Notice that we obtain the above inequality also if $\zeta_\varepsilon= g(\alpha_\varepsilon)| \nabla u_\varepsilon |^\theta$, with $\theta\in [1,\tilde{p})$, up to consider $q'>q$ with $\frac{1}{q'}+\frac{\theta}{p}<\frac{1}{2}$ in the estimates of $\alpha$, since
\begin{equation*}
\frac{\d}{\d t}\big| \nabla u \big|^\theta = \theta | \nabla u \big|^{\theta-2} \nabla u \cdot \nabla \dot{u}\,.
\end{equation*}

Let us now integrate in time the  inequality obtained above for $V_k-V_\varepsilon$: using \eqref{3105181044}, \eqref{eq:estimate 3}, \eqref{3105181043}, \eqref{3105181105}, and~\eqref{3105181042} we deduce that
\begin{equation*}
\|V_k-V_\varepsilon\|_{L^1(0,T;L^1(\Omega))} \to 0\,,
\end{equation*}
and then we get \eqref{eq:convergence of cumulation}, since $f$ is bounded.

Moreover, by weak convergence $\dot \alpha_k \weak \dot \alpha_\epsilon$ in $L^2(0,T;L^2(\Omega))$
\begin{equation*} \label{eq:convergence of dissipation}
\int_0^T \! \R(\dot \alpha_k(t); V_k(t)) \, \d t = - \int_0^T \! \integral{\Omega}{f(V_k(t)) \dot \alpha_k(t)}{\d x} \, \d t \to - \int_0^T \! \integral{\Omega}{f(V_\epsilon(t)) \dot \alpha_\epsilon(t)}{\d x} \, \d t = \int_0^T \! \R(\dot \alpha_\epsilon(t); V_\epsilon(t)) \, \d t
\end{equation*}
and, by \eqref{eq:2202181006},
\[
\begin{split}
\Big| \int_0^T \! & \big( \R(\dot \alpha_k(t); V_k(t)) - \R(\dot \alpha_k(t);\V_k(t)) \big) \, \d t \Big| \leq  \|f(V_k)-f(\V_k)\|_{L^2(0,T;L^2(\Omega))} \|\alpha_k \|_{H^1(0,T;L^2(\Omega))} \leq C_\epsilon \tau \to 0 \, ,
\end{split}
\]
as $k \to +\infty$. This concludes the proof.
\end{proof}
\begin{remark}\label{rem:2202181042}
Combining \eqref{eq:2202181006} and \eqref{eq:convergence of cumulation} we obtain that if \eqref{eq:strong assumption on u} holds, then
\begin{equation}\label{eq:2202181040}
f (\V_k) \to f (V_\epsilon) \quad \text{strongly in } L^2(0,T; L^2(\Omega)) \, .
\end{equation}
\end{remark}

At the moment we do not have convergence~\eqref{eq:strong assumption on u} at our disposal, and we cannot deduce that the convergence of the functions $f(\V_k(t))$ to $f(V_\epsilon(t))$. For this reason, in the following lemma we consider an additional variable $\tilde f_\epsilon(t)$ in the limit evolution, which later in the proof will turn out to be $f(V_\epsilon(t))$.

\begin{lemma}[Compactness for the cumulated variable]\label{le:compCum}
For every $\epsilon > 0$ there exist a nonincreasing function $t \mapsto \tilde f_\epsilon(t) \in L^\infty(\Omega)$ and a subsequence independent of $t$ (which we do not relabel) such that
\begin{equation} \label{eq:convergence of f}
f(\V_k(t)) \wstar \tilde f_\epsilon(t) \quad  \text{weakly* in } L^\infty(\Omega)  \, ,
\end{equation}
for every $t \in [0,T]$.
\end{lemma}
\begin{proof}
To prove the lemma we apply the generalized version of the classical Helly Theorem given in \cite[Helly Theorem]{MalDuc} in the space $\M_b(\Omega)$. For every $t \in [0,T]$, the sequence $\big(f(\V_k(t))\big)_k$ is equibounded in $L^\infty(\Omega)$, and thus is relatively compact in $\M_b(\Omega)$ with respect to the weak* convergence. Moreover, the functions $f(\V_k)$ have uniformly bounded variation in $\M_b(\Omega)$. Indeed, for $s \leq t$ we have $f(\V_k(t)) \leq f(\V_k(s))$ and thus, given a partition $0 = s_0 < \dots < s_m = T$, we get
\[
\sum_{j=1}^m \integral{\Omega}{\big| f(\V_k(s_{j})) - f(\V_k(s_{j-1})) \big| }{\d x} = \integral{\Omega}{ f(\V_k(0)) - f(\V_k(T))  }{\d x} \leq \|f\|_{L^\infty}.
\]

On the one hand, by \cite[Helly Theorem]{MalDuc} we deduce that there exists a subsequence independent of $t$ (which we do not relabel) and a function $t \mapsto \lambda_t \in \M_b(\Omega)$ such that
\begin{equation} \label{eq:05061508}
f(\V_k(t)) \L^2 \mres \Omega \wstar \lambda_t \quad \text{weakly* in } \M_b(\Omega) \, .
\end{equation}
On the other hand, for every $t \in [0,T]$ there exists a function $\tilde f_\epsilon(t) \in L^\infty(\Omega)$ and a subsequence $k_j(t)$ depending on $t$ such that
\begin{equation} \label{eq:05061509}
f(\V_{k_j(t)}(t)) \wstar \tilde f_\epsilon(t) \quad  \text{weakly* in } L^\infty(\Omega)  \, .
\end{equation}
By \eqref{eq:05061508} and \eqref{eq:05061509} we conclude that $\lambda_t = \tilde f_\epsilon(t) \L^2 \mres \Omega$ and the convergence in \eqref{eq:05061509} holds on the whole subsequence $k$ where \eqref{eq:05061508} is satisfied. Notice that $\tilde f_\epsilon(t)$ is nonincreasing in $t$.
\end{proof}

The first step is to deduce the existence of an evolution where the fatigue term $f(V_\epsilon(t))$ is in fact replaced by the term~$\tilde f_\epsilon(t)$. We first prove one inequality in the energy-dissipation balance for the continuous-time evolutions. The opposite inequality will follow automatically from the differential conditions satisfied by~$\alpha_\epsilon$, see Proposition~\ref{prop:balance weak} below.
\begin{proposition}[Energy-dissipation balance in weak form: first inequality]\label{prop:e-d inequality}  For every $\epsilon > 0$ we have
\begin{equation} \label{eq:e-d inequality}
\E(\alpha_\epsilon(T), u_\epsilon(T)) - \int_0^T \! \integral{\Omega}{\tilde f_\epsilon(t) \dot \alpha_\epsilon(t)}{\d x} \, \d t  + \epsilon \int_0^T \! \| \dot \alpha_\epsilon(t) \|^2_{L^2} \, \d t \leq \E(\alpha_0, u_0) + \int_0^T \! \langle \mu(\alpha_\epsilon(t)) \nabla u_\epsilon(t) , \nabla \dot w(t) \rangle_{L^2}  \, \d t  \, .
\end{equation}
\end{proposition}
\begin{proof}

In order to prove~\eqref{eq:e-d inequality}, we write the dissipation with the fatigue term as a supremum of finite sums which are continuous with respect to the convergence~\eqref{eq:convergence of f}. Specifically
\begin{equation} \label{eq:R as sup}
\int_0^T \!  \R(\dot \alpha_k(t);\V_k(t))  \, \d t = \sup_{0=s_0 < \dots < s_m = T} \Big\{ \sum_{j=1}^m \integral{\Omega}{f(\V_k(s_j)) (\alpha_k(s_{j-1}) - \alpha_k(s_{j}))}{\d x} \Big\} \, ,
\end{equation}
where the supremum is taken among all possible partitions $0=s_0 < \dots < s_m =T$, $m \in \NN$, of the interval $[0,T]$. The supremum is in fact attained on the partition $0=t^0_k < \dots < t^k_k = T$. To check this, let us fix a partition $0=s_0 < \dots < s_m =T$ and let us prove that
\begin{equation} \label{eq:06061343}
\begin{split}
\sum_{j=1}^m \integral{\Omega}{f(\V_k(s_j)) (\alpha_k(s_{j-1}) - \alpha_k(s_{j}))}{\d x} & \leq \sum_{i=1}^k \integral{\Omega}{f(\V_k(t^i_k)) (\alpha_k(t^{i-1}_k) - \alpha_k(t^i_k))}{\d x} \\
& = - \int_0^T \integral{\Omega}{f(\V_k(t)) \dot \alpha_k(t)}{\d x} \d t \,.
\end{split}
\end{equation}
Note that if we refine the partition $0=s_0 < \dots < s_m =T$  by including the nodes $t^0_k , \dots , t^k_k$, the dissipation increases, since the monotonicity of $f(\V_k)$ and of $\alpha_k$ yields the following triangular inequality:
\[
\integral{\Omega}{f(\V_k(r_3)) (\alpha_k(r_1) - \alpha_k(r_3))}{\d x} \leq \integral{\Omega}{f(\V_k(r_2)) (\alpha_k(r_1) - \alpha_k(r_2))}{\d x} + \integral{\Omega}{f(\V_k(r_3)) (\alpha_k(r_2) - \alpha_k(r_3))}{\d x}
\]
for $0 \leq r_1 \leq r_2 \leq r_3 \leq T$. Therefore we can assume without loss of generality that $\{t^0_k, \dots, t^k_k\} \subset \{s_0, \dots, s_m \}$. Let us now fix $i \in \{1, \dots, k\}$ and $1 \leq h_i < \ell_i \leq m$ such that $t^{i-1}_k = s_{h_i} < \dots < s_{\ell_i} = t^{i}_k$. Then the sum in in the left-hand side of~\eqref{eq:06061343} can be rearranged as
\[
\begin{split}
& \sum_{i=1}^k \sum_{j=h_i}^{\ell_i} \integral{\Omega}{f(\V_k(s_j)) (\alpha_k(s_{j-1}) - \alpha_k(s_{j}))}{\d x}  = \sum_{i=1}^k \sum_{j=h_i}^{\ell_i} \integral{\Omega}{f(V^{i-1}_k) \frac{s_{j}-s_{j-1}}{\tau} (\alpha^{i-1}_k- \alpha^{i}_k) }{\d x}   \\
& \quad = \sum_{i=1}^k \integral{\Omega}{f(V^{i-1}_k) \frac{t^i_k-t^{i-1}_k}{\tau} (\alpha^{i-1}_k- \alpha^{i}_k) }{\d x} = \sum_{i=1}^k \integral{\Omega}{f(\V_k(t^i_k)) (\alpha_k(t^{i-1}_k)- \alpha_k(t^{i}_k) )}{\d x} \, .
\end{split}
\]

Now we pass to the limit in~\eqref{eq:R as sup} as $k \to +\infty$. Let us fix a partition $0=s_0 < \dots < s_m =T$ and let us fix~$j \in \{0, \dots, m\}$. By~\eqref{eq:pointwise weak alpha} we have in particular that $\alpha_k(s_j) \to \alpha_\epsilon(s_j)$ and $\alpha_k(s_{j-1}) \to \alpha_\epsilon(s_{j-1})$ strongly in~$L^1(\Omega)$ and therefore, by~\eqref{eq:convergence of f}, we obtain that
\begin{equation} \label{eq:limit of f}
\integral{\Omega}{f(\V_k(s_{j})) (\alpha_k(s_{j-1}) - \alpha_k(s_j))}{ \d x} \to \integral{\Omega}{\tilde f_\epsilon(s_j) (\alpha_\epsilon(s_{j-1}) - \alpha_\epsilon(s_j) )}{\d x}
\end{equation}
as $k \to +\infty$.

On the other hand we have that
\begin{equation} \label{eq:int of f as sup}
\sup_{0=s_0 < \dots < s_m = T} \Big\{ \sum_{j=1}^m \integral{\Omega}{\tilde f_\epsilon(s_j) (\alpha_\epsilon(s_{j-1}) - \alpha_\epsilon(s_j) )}{\d x} \Big\} =  - \int_0^T \! \integral{\Omega}{\tilde f_\epsilon(t) \dot \alpha_\epsilon(t)}{\d x} \, \d t \, .
\end{equation}
The equality above follows from a general lemma proved in~\cite[Lemma A.1]{Cri16} regarding the integral representation of weighted variations. To check the fulfillment of the assumptions required by~\cite[Lemma A.1]{Cri16} we remark that:
\begin{itemize}
\item $\alpha_\epsilon \in AC([0,T];L^2(\Omega))$;
\item $\dot \alpha_\epsilon \leq 0$ a.e.\ in $\Omega$;
\item $\tilde f_\epsilon(t) \leq \tilde f_\epsilon(s)$ a.e.\ in $\Omega$ for $s \leq t$;
\item there exists a countable set $E \subset [0,T]$ such that $t \mapsto \tilde f_\epsilon(t)$ is continuous for every $t \in [0,T] \sm E$ with respect to strong $L^2$ topology (this follows from the monotonicity by \cite[Lemma A.2]{Cri16}).
\end{itemize}
Applying \cite[Lemma A.1]{Cri16} with $X := L^2(\Omega)$ and $F=L^2(\Omega)$, we get~\eqref{eq:int of f as sup}.

By~\eqref{eq:R as sup}--\eqref{eq:int of f as sup} we conclude that
\begin{equation} \label{eq:fatigue is lsc}
- \int_0^T \! \integral{\Omega}{\tilde f_\epsilon(t) \dot \alpha_\epsilon(t)}{\d x} \, \d t \leq \liminf_{k \to + \infty} \int_0^T \!  \R(\dot \alpha_k(t);\V_k(t))  \, \d t  \, .
\end{equation}
We conclude the proof using the inequality above together with \eqref{eq:E is lsc}--\eqref{eq:limit of work} and \eqref{eq:discrete balance}.
\end{proof}

\begin{proposition}[Stability in weak form]\label{prop:stweak} Let $\epsilon > 0$. For every $t \in [0,T]$, $u_\epsilon(t)$ is a weak solution to the problem
\begin{equation} \label{eq:pb of ueps}
\left\{
\begin{aligned}
\div\big(\mu(\alpha_\epsilon(t)) \nabla u_\epsilon(t) \big) &= 0  && \text{in } \Omega \, , \\
u_\epsilon(t) & = w(t) && \text{on } \de_D \Omega \, .
\end{aligned}
\right.
\end{equation}
For a.e.\ $t \in (0,T)$ and for every $\beta \in H^1(\Omega)$ with $\beta \leq 0$ a.e.\ in $\Omega$ we have
\begin{equation} \label{eq:viscous stability}
\langle \de_\alpha \E(\alpha_\epsilon(t), u_\epsilon(t)), \beta \rangle - \integral{\Omega}{\tilde f_\epsilon(t) \beta}{\d x} + \epsilon \langle \dot \alpha_\epsilon(t), \beta \rangle_{L^2} \geq 0 \, .
\end{equation}
\end{proposition}
\begin{proof}
To prove~\eqref{eq:pb of ueps} it is sufficient to observe that from~\eqref{eq:pb of uik} we have that $\ol u_k(t)$ is a weak solution to the problem
\begin{equation} \label{eq:pb of uk}
\left\{
\begin{aligned}
\div\big(\mu(\ol \alpha_k(t)) \nabla \ol u_k(t) \big) &= 0  && \text{in } \Omega \, , \\
\ol u_k(t) & = w_k(t) && \text{on } \de_D \Omega \,  ,
\end{aligned}
\right.
\end{equation}
and pass \eqref{eq:pb of uk} to the limit as $k \to +\infty$ using~\eqref{eq:pointwise weak u} and \eqref{eq:pointwise weak alpha}.

Let us fix $\beta \in H^1_-(\Omega)$. Integrating~\eqref{eq:eul1} in time, we get
\begin{equation} \label{eq:discrete balance integrated}
- \int_0^T \! \langle \de_\alpha \E(\ol \alpha_k(t), \ol u_k(t)), \beta \rangle \, \d t + \int_0^T \! \integral{\Omega}{f(\V_k(t)) \beta}{ \d x} \, \d t - \epsilon \int_0^T \! \langle \dot \alpha_k(t) , \beta \rangle_{L^2} \,  \d t \leq 0 \, .
\end{equation}

First of all, we claim that for every $t\in [0,T]$
\begin{equation} \label{eq:2302181101}
-  \langle \de_\alpha \E(\alpha_\epsilon (t),  u_\epsilon(t)), \beta \rangle \leq  \liminf_{k \to +\infty}  -  \langle \de_\alpha \E(\ol \alpha_k(t), \ol u_k(t)), \beta \rangle \, .
\end{equation}
Indeed,  since $\mu'(\ol \alpha_k(t)) \beta \to \mu'(\alpha_\epsilon(t)) \beta$ strongly in $L^2(\Omega)$ and by~\eqref{eq:pointwise weak u}, Ioffe Theorem \cite[Theorem 3.23]{Dac} yields 
\begin{equation} \label{eq:12011718}
- \integral{\Omega}{\mu'(\alpha_\epsilon(t)) \big| \nabla u_\epsilon (t) \big|^2 \beta}{\d x}  \leq \liminf_{k \to +\infty} \Big[ - \integral{\Omega}{\mu'(\ol \alpha_k(t)) \big| \nabla \ol u_k(t)\big|^2 \beta}{\d x} \Big] \, .
\end{equation} 
 Furthermore, by~\eqref{eq:pointwise weak alpha}
\begin{equation} \label{eq:12011720}
\integral{\Omega}{\nabla \ol \alpha_k(t) \cdot \nabla \beta}{\d x}  \to \integral{\Omega}{\nabla \alpha_\epsilon(t) \cdot \nabla \beta}{\d x} \, ,
\end{equation}
for every $t \in [0,T]$. Summing~\eqref{eq:12011718}--\eqref{eq:12011720} we obtain \eqref{eq:2302181101}.
Moreover, by convergence~\eqref{eq:convergence of f}, we have
\begin{equation}\label{eq:2302181104}
\integral{\Omega}{f(\V_k(t)) \beta}{ \d x} \to \integral{\Omega}{\tilde f_\epsilon(t) \beta}{ \d x}
\end{equation}
for every $t \in [0,T]$.

Finally, \eqref{eq:weak alpha} implies
\begin{equation} \label{eq:12011740}
\epsilon \int_{t_1}^{t_2} \! \langle \dot \alpha_k(t) , \beta \rangle_{L^2} \to \epsilon \int_{t_1}^{t_2} \! \langle \dot \alpha_\epsilon(t) , \beta \rangle_{L^2}\,,
\end{equation}
for every $0\leq t_1 \leq t_2\leq T$.

Collecting \eqref{eq:2302181101}, \eqref{eq:2302181104}, \eqref{eq:12011740}, and by~\eqref{eq:discrete balance integrated}, we infer that
\[
- \int_{t_1}^{t_2} \! \langle \de_\alpha \E(\alpha_\epsilon(t), u_\epsilon(t)), \beta \rangle \, \d t + \int_{t_1}^{t_2} \! \integral{\Omega}{\tilde f_\epsilon(t) \beta}{ \d x} \, \d t - \epsilon \int_{t_1}^{t_2} \! \langle \alpha_\epsilon(t) , \beta \rangle_{L^2} \,  \d t \leq 0
\]
for every $\beta \in H^1_-(\Omega)$ and $0\leq t_1 \leq t_2\leq T$.
By the arbitrariness of $t_1$, $t_2$, a localisation argument gives~\eqref{eq:viscous stability}.
\end{proof}

\begin{remark}
Using~\eqref{eq:pb of ueps} we can improve the convergence in~\eqref{eq:uniform convergence of u}, namely for every $\epsilon > 0$
\begin{equation} \label{eq:improved convergence on u}
 \| u_k - u_\epsilon \|_{C([0,T];W^{1,p}(\Omega))}, \ \| \ol u_k - u_\epsilon \|_{L^\infty(0,T;W^{1,p}(\Omega))}, \ \| \ul u_k - u_\epsilon \|_{L^\infty(0,T;W^{1,p}(\Omega))}  \to 0 \, , \text{ for } p \in [2, \tilde p \, ) \, .
\end{equation}
Indeed by~\eqref{eq:pb of uk} and \eqref{eq:pb of ueps} we deduce that the function $v := u_\epsilon(t) - \ol u_k(t) - w(t) + \ol w_k(t)$ is a weak solution to the problem
\[
\left\{
\begin{aligned}
\div\big(\mu(\alpha_\epsilon(t)) \nabla v) & = \ell && \text{in } \Omega \, , \\
v & = 0 && \text{on } \de_D \Omega \, ,
\end{aligned}
\right.
\]
where $\ell \in W^{-1,p}_{\de_D \Omega}(\Omega)$ is defined by $\ell := \div \big( (\mu(\ol \alpha_k(t)) - \mu(\alpha_\epsilon(t))) \nabla \ol u_k(t) \big) + \div \big( \mu(\alpha_\epsilon(t)) (\nabla \ol w_k(t) - \nabla w(t))  \big)$. By Remark~\ref{rmk:integrability}, \eqref{eq:estimate of u}, \eqref{eq:uniform convergence of alpha}, and~\eqref{ass:w} we deduce that
\[
\begin{split}
\| \ol u_k(t) - u_\epsilon(t) \|_{W^{1,p}} & \leq C \Big[ \| \ol \alpha_k(t) - \alpha_\epsilon(t) \|_{L^q} \| \ol u_k(t)\|_{W^{1,\tilde p}} + \|\alpha_\epsilon(t)\|_{L^q} \|\ol w_k(t) - w(t)\|_{W^{1, \tilde p}} \Big] \\
& \leq C \Big[ \| \ol \alpha_k(t) - \alpha_\epsilon(t) \|_{L^q} + \|\ol w_k(t) - w(t)\|_{W^{1, \tilde p}} \Big]  \to 0
\end{split}
\]
uniformly with respect to $t$, for a suitable $q \in (2, \infty)$. The convergence of $u_k$ and $\ul u_k$ follows from~\eqref{eq:affine and piecewise same u ol}--\eqref{eq:affine and piecewise same u ul}. 
\end{remark}

\begin{proposition}[Energy-dissipation balance in weak form]  \label{prop:balance weak}
For every $\epsilon > 0$ we have
\[
\E(\alpha_\epsilon(T), u_\epsilon(T)) - \int_0^T \! \integral{\Omega}{\tilde f_\epsilon(t) \dot \alpha_\epsilon(t)}{\d x} \, \d t  + \epsilon \int_0^T \! \| \dot \alpha_\epsilon(t) \|^2_{L^2} \, \d t = \E(\alpha_0, u_0) + \int_0^T \! \langle \mu(\alpha_\epsilon(t)) \nabla u_\epsilon(t) , \nabla \dot w(t) \rangle_{L^2}  \, \d t  \, .
\]
\end{proposition}

\begin{proof}
One inequality has been proven in Proposition~\ref{prop:e-d inequality}. To prove the opposite inequality, we observe that $ t \mapsto \E(\alpha_\epsilon(t), u_\epsilon(t))$ is absolutely continuous  and
\[
\begin{split}
\frac{\d}{\d t} \Big[ \E(\alpha_\epsilon(t), u_\epsilon(t)) \Big] & = \langle \de_\alpha \E(\alpha_\epsilon(t), u_\epsilon(t) ) , \dot \alpha_\epsilon(t) \rangle + \langle \de_u \E(\alpha_\epsilon(t), u_\epsilon(t) ) , \dot u_\epsilon(t) \rangle  \\
& \geq \integral{\Omega}{\tilde f_\epsilon(t) \dot \alpha_\epsilon(t)}{\d x} - \epsilon \| \dot \alpha_\epsilon(t)\|^2_{L^2} + \langle \mu(\alpha_\epsilon(t)) \nabla u_\epsilon(t) ) , \nabla \dot w(t) \rangle_{L^2}
\end{split}
\]
for a.e.\ $t \in (0,T)$, where in the last inequality we have used~\eqref{eq:viscous stability} and \eqref{eq:pb of ueps}. Integrating the previous inequality in time, we complete the proof.
\end{proof}

The energy-dissipation balance obtained in Proposition~\ref{prop:balance weak} above allows us to get the desired strong convergence~\eqref{eq:strong assumption on u}.

\begin{proposition}[Strong convergence of discrete-time evolutions]\label{prop:strconveps}
For every $\epsilon >0$ we have
\begin{subequations} \label{eq:strong convergences}
\begin{align} \label{eq:strong convergence of u}
\alpha_k & \to \alpha_\epsilon \quad \text{strongly in } W^{1,1}(0,T; L^q(\Omega)) \, ,\quad  \text{for } q \in [1, \infty)  \, , \\
u_k & \to u_\epsilon \quad  \text{strongly in } W^{1,1}(0,T; W^{1,p}(\Omega)) \, ,  \quad \text{for } p \in [2, \tilde p) \, .
\end{align}
\end{subequations}
\end{proposition}
\begin{proof}
From Proposition~\ref{prop:discrete balance} and Proposition~\ref{prop:balance weak} and using the convergence of the work term~\eqref{eq:limit of work}, we deduce that
\[
\begin{split}
\lim_{k \to +\infty} & \Big[ \E(\alpha_k(T), u_k(T)) + \int_0^T \! \R(\dot \alpha_k(t); \V_k(t)) \, \d t + \epsilon \int_0^T \! \|\dot \alpha_k(t)\|^2_{L^2} \d t  \Big]  \\
& = \E(\alpha_\epsilon(T), u_\epsilon(T)) - \int_0^T \! \integral{\Omega}{\tilde f_\epsilon(t) \dot \alpha_\epsilon(t)}{\d x} \, \d t  + \epsilon \int_0^T \! \| \dot \alpha_\epsilon(t) \|^2_{L^2} \, \d t \, .
\end{split}
\]
Notice that if $(a_k)_k$ and $(b_k)_k$ are two sequences such that $a_k + b_k \to a + b$ and $a \leq \liminf_k a_k$, $b \leq \liminf_k b_k$, then $a_k \to a$ and $b_k \to b$. Therefore, by~\eqref{eq:E is lsc}, \eqref{eq:viscous diss is lsc}, and \eqref{eq:fatigue is lsc} we obtain that
\begin{gather*}
\lim_{k \to +\infty} \int_0^T \! \|\dot \alpha_k(t)\|^2_{L^2} \d t = \int_0^T \! \| \dot \alpha_\epsilon(t) \|^2_{L^2} \, \d t \, .
\end{gather*}
As a consequence
\begin{equation} \label{eq:L2L2 convergence of alpha}
\dot \alpha_k \to \dot \alpha_\epsilon \quad \text{strongly in } L^2(0,T; L^2(\Omega)) \, .
\end{equation}
We want to deduce the strong convergence~\eqref{eq:strong convergence of u} from \eqref{eq:L2L2 convergence of alpha}. In order to do so, we shall control $\| \dot u_\epsilon(t) - \dot u_k(t) \|_{W^{1,p}}$ with $\| \dot \alpha_\epsilon(t) - \dot \alpha_k(t) \|_{L^q}$ for some $q \in (2, \infty)$, as we did in the proof of~\eqref{eq:estimate of increment}. For this reason it is necessary to slightly improve the integrability in the target space in~\eqref{eq:L2L2 convergence of alpha}. More precisely, we claim that for every $q \in [1, \infty)$
\begin{equation}  \label{eq:L1Lq convergence of alpha}
\dot \alpha_k \to \dot \alpha_\epsilon \quad \text{strongly in } L^1(0,T; L^q(\Omega))  \, .
\end{equation}
Indeed, let us fix $\theta \in (0,1)$ and $q > 2$ (the case $q \leq 2$ being already covered by \eqref{eq:L2L2 convergence of alpha}) and let us define $r > q$ in such a way that $\frac{1}{q} = \frac{\theta}{2} + \frac{1-\theta}{r}$. Using the interpolation inequality between the spaces $L^2(\Omega)$ and $L^q(\Omega)$, H\"older's Inequality, \eqref{eq:estimate 3}, \eqref{eq:estimate 3 for alphaeps ueps}, and \eqref{eq:L2L2 convergence of alpha} we obtain that
\[
\begin{split}
\int_0^T \! \|\dot \alpha_k(t) - \dot \alpha_\epsilon(t) \|_{L^q} \, \d t  & \leq \int_0^T \! \|\dot \alpha_k(t) - \dot \alpha_\epsilon(t) \|^\theta_{L^2} \|\dot \alpha_k(t) - \dot \alpha_\epsilon(t) \|^{1-\theta}_{L^r}  \, \d t \\
& \leq \Big( \int_0^T \! \|\dot \alpha_k(t) - \dot \alpha_\epsilon(t) \|_{L^2} \, \d t \Big)^{\! \theta} \Big( \int_0^T \! \|\dot \alpha_k(t) - \dot \alpha_\epsilon(t) \|_{H^1} \, \d t \Big)^{\! 1-\theta}  \\
& \leq C \Big( \int_0^T \! \|\dot \alpha_k(t) - \dot \alpha_\epsilon(t) \|_{L^2} \, \d t \Big)^{\! \theta} \to 0
\end{split}
\]
as $k \to +\infty$. This proves~\eqref{eq:L1Lq convergence of alpha}.

We are now ready to prove~\eqref{eq:strong convergence of u}.
Differentiating~\eqref{eq:pb of ueps} in time and by~\eqref{eq:pb of dot uk} we obtain that for a.e.\ $t \in (0,T)$ the function $v :=  \dot u_\epsilon(t) -  \dot u_k(t) - \dot w(t) + \dot w_k(t)$ is a weak solution to the problem
\[
\left\{
\begin{aligned}
\div\big(\mu(\alpha_\epsilon(t)) \nabla v) & = \ell && \text{in } \Omega \, , \\
v & = 0 && \text{on } \de_D \Omega \, ,
\end{aligned}
\right.
\]
where $\ell \in W^{-1,p}_{\de_D \Omega}(\Omega)$ is defined by
\[
\begin{split}
\ell & := - \div\Big( \mu'(\alpha_\epsilon(t)) \dot \alpha_\epsilon(t) \big( \nabla u_\epsilon(t) - \nabla \ol u_k(t) \big)  \Big) - \div\Big( \mu(\alpha_\epsilon(t)) \big(  \nabla \dot w(t) - \nabla \dot w_k(t) \big) \Big) \\
& \quad + \div\Big( \big( \mu(\ul \alpha_k(t)) -\mu(\alpha_\epsilon(t)) \big) \nabla \dot u_k(t) \Big) + \div\Big( \big( \tfrac{\mu(\ul \alpha_k(t) + \tau \dot \alpha_k(t)) - \mu(\ul \alpha_k(t))}{\tau} - \mu'(\alpha_\epsilon(t)) \dot \alpha_\epsilon(t) \big )  \nabla \ol u_k(t) \Big) \, .
\end{split}
\]
Observe that for a.e.\ $t \in \Omega$
\[
\begin{split}
\Big| & \tfrac{\mu(\ul \alpha_k(t) + \tau \dot \alpha_k(t)) - \mu(\ul \alpha_k(t))}{\tau} - \mu'(\alpha_\epsilon(t)) \dot \alpha_\epsilon(t) \Big| \\
& \leq \Big| \tfrac{\mu(\ul \alpha_k(t) + \tau \dot \alpha_k(t)) - \mu(\ul \alpha_k(t))}{\tau} - \mu'(\ul \alpha_k(t)) \dot \alpha_k(t) \Big| + \big| \mu'(\ul \alpha_k(t))   - \mu'(\alpha_\epsilon(t)) \big| \big| \dot \alpha_k(t) \big| + \big| \mu'(\alpha_\epsilon(t)) \big| \big| \dot \alpha_k(t) - \dot \alpha_\epsilon(t)\big| \\
& \leq C \Big[ \tau \big| \dot \alpha_k(t) \big| + \big| \ul \alpha_k(t) - \alpha_\epsilon(t) \big|  \big| \dot \alpha_k(t) \big| + \big| \dot \alpha_k(t) - \dot \alpha_\epsilon(t)\big| \Big]
\end{split}
\]
a.e.\ in $\Omega$. Therefore, by Remark~\ref{rmk:integrability}, \eqref{eq:estimate of u}, \eqref{eq:estimate of increment}, and \eqref{eq:W1tilde norm} we get that
\[
\begin{split}
\| \dot u_k(t) - \dot u_\epsilon(t) \|_{W^{1,p}} & \leq C  \Big[  \| \dot \alpha_\epsilon(t) \|_{L^q} \|   u_\epsilon(t) -  \ol u_k(t) \|_{W^{1, p_1}} + \| \dot w_k(t) - \dot w(t) \|_{W^{1,\tilde p}} \\
& \qquad + \| \ul \alpha_k(t) -  \alpha_\epsilon(t)\|_{L^q}  \|\dot u_k(t) \|_{W^{1, p_1}} +  \tau \|\dot \alpha_k(t)\|_{L^q} \|\ol u_k(t)\|_{W^{1,\tilde p}} \\
& \qquad + \|\ul \alpha_k(t) - \alpha_\epsilon(t) \|_{L^q} \| \dot \alpha_k(t) \|_{L^r} \| \ol u_k(t)\|_{W^{1,p_1}} + \| \dot \alpha_k(t) - \dot \alpha_\epsilon(t)\|_{L^q} \| \ol u_k(t)\|_{W^{1,\tilde p}}  \Big] \\
& \leq C  \Big[  \| \dot \alpha_\epsilon(t) \|_{H^1}\|   u_\epsilon(t) -  \ol u_k(t) \|_{W^{1,p_1}} + \| \dot w_k(t) - \dot w(t) \|_{W^{1,\tilde p}} \\
& \qquad + \|\ul \alpha_k(t) - \alpha_\epsilon(t) \|_{L^q} \| \dot \alpha_k(t) \|_{H^1} +  \tau \|\dot \alpha_k(t)\|_{H^1} + \| \dot  \alpha_k(t) - \dot \alpha_\epsilon(t)\|_{L^q}  \Big] \, ,
\end{split}
\]
where $q, r \in (2,\infty)$, and $p_1 \in (2, \tilde p \, )$ are suitable exponents. Integrating in time the previous inequality and by H\"older's Inequality we obtain
\[
\begin{split}
\| \dot u_k - \dot u_\epsilon \|_{L^1(0,T;W^{1,p}(\Omega))} &  \leq C \Big[ \| \alpha_\epsilon \|_{W^{1,1}(0,T;H^1(\Omega))} \| u_\epsilon -  \ol u_k \|_{L^\infty(0,T;W^{1,p_1}(\Omega))} + \| w_k -  w \|_{W^{1,1}(0,T;W^{1,\tilde p}(\Omega))} \\
& \qquad + \|\ul \alpha_k - \alpha_\epsilon \|_{L^\infty(0,T;L^q(\Omega))} \| \alpha_k \|_{W^{1,1}(0,T;H^1(\Omega))} +   \tau \|\dot \alpha_k\|_{W^{1,1}(0,T;H^1(\Omega))} \\
& \qquad + \| \dot  \alpha_k - \dot \alpha_\epsilon\|_{L^1(0,T;L^q(\Omega))} \Big] \\
& \leq C \Big[  \| u_\epsilon -  \ol u_k \|_{L^\infty(0,T;W^{1,p_1}(\Omega))} + \| w_k -  w \|_{W^{1,1}(0,T;W^{1,\tilde p}(\Omega))}\\
& \qquad + \|\ul \alpha_k - \alpha_\epsilon \|_{L^\infty(0,T;L^q(\Omega))} + \tau  + \| \dot  \alpha_k - \dot \alpha_\epsilon\|_{L^1(0,T;L^q(\Omega))} \Big] \, .
\end{split}
\]
By~\eqref{eq:improved convergence on u}, \eqref{ass:w}, \eqref{eq:uniform convergence of alpha}, and~\eqref{eq:L1Lq convergence of alpha} we conclude that the right-hand side in the inequality above converges to zero as $k \to +\infty$.
\end{proof}

\begin{proof}[Proof of Theorem~\ref{teo:exViscEvo}]
For fixed $\varepsilon>0$, Propositions~\ref{prop:stweak} and \ref{prop:balance weak} show that $(\alpha_\epsilon, u_\epsilon)$, obtained by \eqref{eq:weak alpha}, \eqref{eq:weak u} as weak limit of a sequence of discrete-time evolutions $(\alpha_k, u_k)$, satisfy the conditions of Definition~\ref{def:appreveps} in a weak sense. In fact, (ev0)$_\epsilon$, (ev1)$_\epsilon$ hold, while (ev2)$_\epsilon$, (ev3)$_\epsilon$ are satisfied with $\tilde{f}_\epsilon(t)$ in place of $f(V_\epsilon(t))$, where $\tilde{f}_\epsilon(t)$ is such that (cf.\ \eqref{eq:convergence of f})
\begin{equation*}
f(\V_k(t)) \wstar \tilde f_\epsilon(t) \quad  \text{weakly* in } L^\infty(\Omega)  \, ,
\end{equation*}
  for every $t\in [0,T]$.

Actually we find, in Proposition~\ref{prop:strconveps}, that \emph{a posteriori} we have an enhanced convergence for the displacement evolutions that guarantees the strong convergence
\begin{equation*}
f (\V_k) \to f (V_\epsilon) \quad \text{strongly in } L^2(0,T; L^2(\Omega)) \, ,
\end{equation*}
by Lemma~\eqref{lemma:desired convergence} and Remark~\ref{rem:2202181042}. We conclude that for a.e.\ $t\in (0,T)$
\[
\tilde f_\epsilon(t)= f (V_\epsilon(t)) \,,
\]
 so that (ev2)$_\epsilon$, (ev3)$_\epsilon$ are satisfied with $f(V_\epsilon(t))$ and $(\alpha_\epsilon, u_\epsilon)$ is an $\epsilon$-approximate viscous evolution. The estimate \eqref{2202180044} follows immediately from \eqref{eq:estimate 3 for alphaeps ueps}.
\end{proof}

We conclude this section by a characterisation of the energy balance for $\epsilon$-approximate viscous evolutions, that will be employed in the next section to pass to the limit as $\epsilon$ tends to 0.
We first deduce the following lemma.
\begin{lemma}\label{le:2202181112}
Let $(\alpha_\epsilon,u_\epsilon)\in H^1(0,T; H^1(\Omega))\times H^1(0,T;W^{1,p}(\Omega))$ satisfies (ev0)$_\epsilon$, (ev1)$_\epsilon$ of Definition~\ref{def:appreveps}. Then (ev3)$_\epsilon$ for $(\alpha_\epsilon,u_\epsilon)$ is equivalent to:
\begin{itemize}
\item[{\rm (ev3')$_\epsilon$}:] for a.e.\ $t\in (0,T)$
\begin{equation}\label{3006181736}
\langle \de_\alpha \E(\alpha_\epsilon(t), u_\epsilon(t) ) , \dot \alpha_\epsilon(t) \rangle + \R(\dot{\alpha}_\epsilon(t); V_\epsilon(t))  + \epsilon \| \dot \alpha_\epsilon(t) \|^2_{L^2} = 0 \,.
\end{equation}
\end{itemize}
\end{lemma}

\begin{proof}
Being $\alpha_\epsilon$, $u_\epsilon$
absolutely continuous (in time) we get that $t\mapsto \mathcal{E}(\alpha_\epsilon(t), u_\epsilon(t))$ is absolutely continuous and
\begin{equation}\label{eq:2202181125}
\begin{split}
\frac{\d}{\d t} \Big[ \E(\alpha_\epsilon(t), u_\epsilon(t)) \Big] & = \langle \de_\alpha \E(\alpha_\epsilon(t), u_\epsilon(t) ) , \dot \alpha_\epsilon(t) \rangle + \langle \de_u \E(\alpha_\epsilon(t), u_\epsilon(t) ) , \dot u_\epsilon(t) \rangle  \\
& = \langle \de_\alpha \E(\alpha_\epsilon(t), u_\epsilon(t) ) , \dot \alpha_\epsilon(t) \rangle + \langle \mu(\alpha_\epsilon(t)) \nabla u_\epsilon(t) , \nabla \dot w(t) \rangle_{L^2}\,,
\end{split}
\end{equation}
using (ev1)$_\epsilon$. Differentiating in time (ev3)$_\epsilon$ gives then the equivalence between (ev3)$_\epsilon$ and (ev3')$_\epsilon$.
\end{proof}
\begin{remark}\label{rem:2202181151}
Arguing in a similar way (cf.\ also Proposition~\ref{prop:balance weak} and \cite[Proposition~4.2]{CriLaz16}) it is not difficult to see that if $(\alpha_\epsilon,u_\epsilon)\in H^1(0,T; H^1(\Omega)\times W^{1,p}(\Omega))$ satisfies (ev0)$_\epsilon$, (ev1)$_\epsilon$, (ev2)$_\epsilon$ of Definition~\ref{def:appreveps} and
\begin{itemize}
\item[{\rm (ev3'')$_\epsilon$}:]
\[
\E(\alpha_\epsilon(T), u_\epsilon(T)) + \int_0^T \! \R(\dot{\alpha}_\epsilon(t); V_\epsilon(t))\, \d t  + \epsilon \int_0^T \! \| \dot \alpha_\epsilon(t) \|^2_{L^2} \, \d t \leq \E(\alpha_0, u_0) + \int_0^T \! \langle \mu(\alpha_\epsilon(t)) \nabla u_\epsilon(t) , \nabla \dot w(t) \rangle_{L^2}  \, \d t  \,,
\]
\end{itemize}
then $(\alpha_\epsilon,u_\epsilon)$ is an $\epsilon$-approximate viscous evolution.
\end{remark}
Let us introduce some notation in view of the characterisation of the energy balance for $\epsilon$-approximate viscous evolutions. 

For $(\alpha, u) \in W^{1,1}(0,T;H^1(\Omega){\times}H^1(\Omega))$ and $\tilde{f}\in L^2(\Omega)$ (that we regard as an element of $(H^1(\Omega))'$ with $\langle \tilde{f}, \beta \rangle=\int_\Omega \tilde f\,\beta \d x$) we define
\begin{equation}\label{eq:2202181214}
\Phi(g):=\sup_{\beta\in F} \langle -g, \beta \rangle \quad\text{ for every }g\in (H^1(\Omega))'\,,\qquad \Psi(\alpha, u, \tilde{f}):=\Phi\big(\de_\alpha \E(\alpha, u )-\tilde{f} \, \big)\,,
\end{equation}
where
\begin{equation*}
F:= \{\beta \in H^1_-(\Omega) \colon \|\beta\|_{L^2}\leq 1\}
\end{equation*}

Employing Lemma~\ref{le:2202181112} we obtain the following characterisation of the energy balance, which is invariant under time reparametrisation.
\begin{proposition}\label{prop:2202181211}
Let $(\alpha_\epsilon, u_\epsilon)$ be an $\epsilon$-approximate viscous evolution. Then with the notation above we have that
\begin{equation}\label{eq:2202181319}
\epsilon \|\dot{\alpha}_\epsilon(t)\|_{L^2}=\Psi(\alpha_\epsilon(t), u_\epsilon(t), f(V_\epsilon(t)))\,,
\end{equation}
and one may recast the energy balance (ev3)$_\epsilon$ as
 \begin{equation}\label{eq:2302181249}
 \begin{split}
\E(\alpha_\epsilon(T), u_\epsilon(T)) + \int_0^T \! \R(\dot{\alpha}_\epsilon(t); V_\epsilon(t))\, \d t  &+  \int_0^T \! \| \dot \alpha_\epsilon(t) \|_{L^2} \Psi(\alpha_\epsilon(t), u_\epsilon(t), f(V_\epsilon(t))) \, \d t \\&= \E(\alpha_0, u_0) + \int_0^T \! \langle \mu(\alpha_\epsilon(t)) \nabla u_\epsilon(t) , \nabla \dot w(t) \rangle_{L^2}  \, \d t  \, .
\end{split}
\end{equation}
\end{proposition}

\begin{proof}
By (ev2)$_\epsilon$ we get that for every $\beta \in H^1_-(\Omega)$
\begin{equation*}
\epsilon \langle \dot{\alpha}_\epsilon(t), \beta \rangle \geq \langle -\de_\alpha \E(\alpha_\epsilon(t), u_\epsilon(t)), \beta \rangle+ \R(\beta; V_\epsilon(t)) \,.
\end{equation*}
On the other hand Lemma~\ref{le:2202181112} implies that the equality above is attained for $\ol\beta=\frac{\dot{\alpha}_\epsilon(t)}{\|\dot{\alpha}_\epsilon(t)\|_{L^2}}$ and this gives \eqref{eq:2202181319}, since $\ol \beta$ is in $F$ and $\R(\beta; V_\epsilon(t))  = -\int_\Omega f(V_\epsilon(t)) \, \beta \d x$. Then \eqref{eq:2302181249} follows immediately from the energy balance (ev3)$_\epsilon$.
\end{proof}

\begin{remark}\label{rem:2302182016}
Arguing in the same way of \cite[Lemma~4.4]{CriLaz16} (see also \cite[Lemma~A.2]{Neg16}) we deduce that
\begin{equation*}
\Phi(g)=\d_2(g, G)\quad\text{ for every }g \in (H^1(\Omega))'\,,
\end{equation*}
where
\begin{equation*}
G:= \{h\in (H^1(\Omega))' \colon \langle h, \beta \rangle \geq 0 \, \text{ for every }\beta \in H^1_-(\Omega)\}\,, \qquad
\d_2(g,G):=\min\{\|h\|_{L^2} \colon h\in L^2(\Omega)\,, \, h+g \in G\}\,.
\end{equation*}
\end{remark}
\section{Vanishing viscosity limit}

This section concerns the asymptotics of the viscous evolution, whose existence has been proven in Section~\ref{Sec:3}, under the constitutive assumptions in Section~\ref{sec:Ass}, as the viscosity parameter $\varepsilon$ vanishes. We use a rescaling technique, common to many other works (see e.g.\ \cite{DMDSSol, KRZ13a, KRZ15, CriLaz16}).
Let $\{(\alpha_\epsilon, u_\epsilon)\}_{\epsilon>0}$ be a family of $\epsilon$-approximate viscous evolutions satisfying the uniform $W^{1,1}$ bound in time \eqref{2202180044}, for a given $p<\tilde{p}$, where $\tilde{p}$ is given by Lemma~\ref{lemma:increment of u}. The existence of these evolutions has been shown in Theorem~\ref{teo:exViscEvo}.
For $\epsilon>0$ and $t\in [0,T]$ we set
\begin{equation}\label{eq:2202181908}
\seo(t):=t+\int_0^t \|\dot{\alpha}_\epsilon(s)\|_{H^1} \d s + \int_0^t \|\dot{u}_\epsilon(s)\|_{W^{1,p}} \d s\,.
\end{equation}
Then $\seo$ is absolutely continuous and
\begin{equation*}
\seo(t_2)-\seo(t_1)\geq t_2-t_1 \quad \text{ for every }0\leq t_1 \leq t_2 \leq S_\epsilon:=\seo(T)\,,
\end{equation*}
in particular $\seo$ is strictly increasing and bijective on its domain. We denote by $\teo\colon [0,S_\epsilon]\to [0,T]$ the inverse of $\seo$. In view of \eqref{2202180044}, we have that $T\leq S_\epsilon<C$, for $C>0$ independent of $\epsilon$, and then, up to a subsequence, $S_\epsilon\to S$ as $\epsilon\to 0$, with $S\geq T$. We define the rescaled evolution on $[0,S_\epsilon]$ by setting
\begin{equation}\label{eq:2302181228}
\aeo(s):=\alpha_\epsilon(\teo(s))\,, \quad \ueo(s):=u_\varepsilon(	\teo(s))\,, \quad \zeta_\varepsilon^\circ(s):=\zeta_\epsilon(\teo(s))\,, \quad \Veo(s):=V_\epsilon(\teo(s)) \,.
\end{equation}
Up to extending $\teo$ with $\teo(S_\epsilon)$ in $(S_\epsilon, \ol S]$, for $\ol S:= \sup_{\varepsilon>0}S_\epsilon$ ($\epsilon$ small), we assume the rescaled functions above defined on the fixed time interval $[0,S]$.
By a change of variable we have from \eqref{eq:def of Veps} that
\begin{equation*}
\Veo(s)=\int_0^s |\dot{\zeta}_\epsilon^\circ(\sigma)|\d \sigma \quad \text{a.e.\ in }\Omega, \text{ for every }s\in [0,S]\,.
\end{equation*}
Since \eqref{eq:2202181908} gives that $\teo$ is nondecreasing and that
\begin{equation}\label{eq:2202181909}
\teo(s_2)-\teo(s_1) +\|\aeo(s_2)-\aeo(s_1)\|_{H^1} + \|\ueo(s_2)-\ueo(s_1)\|_{W^{1,p}} \leq s_2-s_1
\end{equation}
for every $0\leq s_1\leq s_2\leq S$, we deduce (cf.\ also e.g.\ \cite{DMDSSol, CriLaz16, KRZ15}) that, up to a (not relabeled) subsequence
\begin{equation}\label{eq:2202182206}
(\teo, \aeo, \ueo) \wstar (t^\circ, \alpha^\circ, u^\circ)\quad \text{weakly$^*$ in }W^{1,\infty}(0,S; [0,T]{\times}H^1(\Omega){\times} W^{1,p}(\Omega))\,,
\end{equation}
for a suitable $(t^\circ, \alpha^\circ, u^\circ)$ with
\begin{equation*}
\dot{t}^\circ(s)+\|\dot{\alpha}^\circ(s)\|_{H^1} +  \|\dot{u}^\circ(s)\|_{W^{1,p}}\leq 1 \quad \text{ for a.e.\ }s\in (0,S)\,.
\end{equation*}
In view of the equicontinuity (with respect to $\varepsilon$) of $(\aeo, \ueo)$, it follows that for every $s\in [0,S]$ and $s_\epsilon\to s$
\begin{equation}\label{eq:2202182102}
\aeo(s_\epsilon)\weak \alpha^\circ(s) \quad \text{weakly in }H^1(\Omega)\,, \qquad \ueo(s_\epsilon)\weak u^\circ(s)\quad \text{weakly in }W^{1,p}(\Omega)\,.
\end{equation}
 Moreover, we define
\begin{equation} \label{eq:wcirc}
w^\circ(s) := w(t^\circ(s)) \, , \quad \text{for every } s \in [0,S] \, .
\end{equation}

Similarly to the analogous situation in Section~\ref{Sec:3}, the weak convergences above are not enough to guarantee pointwise convergence for the cumulations of the strains, even if the cumulation of $\zeta^\circ$
\begin{equation*}
V^\circ(s)=\int_0^s \big|\dot{\zeta}^\circ(\sigma)\big|\d \sigma \quad \text{a.e.\ in }\Omega, \text{ for every }s\in [0,S]\,
\end{equation*}
is well defined as a Bochner integral in $L^2(\Omega)$.
We may only say, passing through an Helly type selection principle as in Lemma~\ref{le:compCum} that there exists $\tilde{f}^\circ \colon [0,S]\to L^\infty(\Omega)$, increasing in time for a.e.\ fixed $x \in \Omega$, such that
\begin{equation}\label{eq:convCumulEps}
f(\Veo(s))\wstar \tilde{f}^\circ(s) \quad\text{weakly$^*$ in }L^\infty(\Omega)\,, \text{ for every $s\in [0,S]$}\,.
\end{equation}
However, differently from Section~\ref{Sec:3}, in view of the loss of the viscous term in the limit evolution we are not able to improve the convergences \eqref{eq:2202182206} \emph{a posteriori}, so to express $\tilde{f}^\circ$ in terms of $V^\circ$, but we prove only an inequality, see Proposition~\ref{prop:ordineCum}.

We obtain then the following existence result for limit of rescaled $\varepsilon$-approximate viscous evolutions, that we call \emph{rescaled quasistatic viscosity evolutions}, which is the main result of the paper.
\begin{theorem}\label{teo:evoRepar}
With the notation above the function $(t^\circ, \alpha^\circ, u^\circ)\in W^{1,\infty}(0,S; [0,T]{\times}H^1(\Omega){\times} W^{1,p}(\Omega))$, defined as limit of rescaled $\epsilon$-approximate viscous evolutions in \eqref{eq:2202182206}, satisfies the following properties:
\begin{itemize}
\item[{\rm (ev0)}]{\em irreversibility\/}:
\begin{equation*} \label{ir0}
 [0,S] \ni s \mapsto \alpha^\circ(s) \quad \text{is nonincreasing as a family of measurable functions on }\Omega\,,
 \end{equation*}
  that is $\alpha^\circ(t) \leq \alpha^\circ(s)$ a.e.\ in $\Omega$ for all $s \leq t$; 
 \item[{\rm (ev1)}] {\em equilibrium\/}: for every $s \in [0,S]$,  $u^\circ(s) \in H^1(\Omega)$ is a weak solution to the problem 
\begin{equation} \label{eq:2302182006}
\left\{
\begin{aligned}
\div\big(\mu(\alpha^\circ(s)) \nabla u^\circ(s) \big) &= 0  && \text{in } \Omega \, , \\
u^\circ(s) & = w^\circ(s) && \text{on } \de_D \Omega \, .
\end{aligned}
\right.
\end{equation}
\item[{\rm (ev2)}] {\em  Karush-Kuhn-Tucker inequality \/}:  for a.e.\ $s \in (0,S)\setminus U^\circ$ and  for every $\beta \in H^1(\Omega)$ with $\beta \leq 0$ a.e.\ in $\Omega$ we have
\begin{equation} \label{eq:2302182007}
\langle \de_\alpha \E(\alpha^\circ(s), u^\circ(s)) -\tilde{f}^\circ(s), \beta \rangle \geq 0 \, ,
\end{equation}
 where $U^\circ:=\{s\in [0,S]\colon t^\circ \text{ is constant in a neighbourhood of } s \}$. 
\item[{\rm (ev3)}] {\em energy balance\/}:
\[
\begin{split}
\E(\alpha^\circ(S), u^\circ(S)) &- \int_0^S \! \langle \tilde{f}^\circ(s), \dot{\alpha}^\circ(s) \rangle\, \d s +\int_0^S \|\dot{\alpha}^\circ(s)\|_{L^2}\,\Psi(\alpha^\circ(s), u^\circ(s), \tilde{f}^\circ(s))\,\mathrm{d}s  \\&= \E(\alpha_0, u_0) + \int_0^S \! \langle \mu(\alpha^\circ(s)) \nabla u^\circ(s) , \nabla \dot w^\circ(s) \rangle  \, \d s  \, .
\end{split}
\]
\end{itemize}
Moreover, for every $s\in [0,S]$ we have that
\begin{equation}\label{0106191253}
\E(\alpha^\circ(s), u^\circ(s))=\lim_{\epsilon\to 0}\E(\aeo(s), \ueo(s))\,,
\end{equation}
\begin{equation*}\label{eq:2302182320}
- \int_0^S \! \langle \tilde{f}^\circ(s), \dot{\alpha}^\circ(s) \rangle\, \d s= -\lim_{\epsilon\to 0}  \int_0^S \! \langle f(\Veo(s)), \dot{\alpha}^\circ_\epsilon(s) \rangle\, \d s= \lim_{\epsilon\to 0} \int_0^S \! \mathcal{R}(\dot{\alpha}^\circ_\epsilon(s); \Veo(s))\, \d s\,,
\end{equation*}
and
\begin{equation*}
\int_0^S \|\dot{\alpha}^\circ(s)\|_{L^2}\,\Psi(\alpha^\circ(s), u^\circ(s), \tilde{f}^\circ(s))\,\mathrm{d}s=\lim_{\epsilon\to 0} \int_0^S  \| \dot \alpha_\epsilon^\circ(s) \|_{L^2} \Psi(\aeo(s), \ueo(s), f(\Veo(s))) \, \d s\,.
\end{equation*}
\end{theorem}

\begin{remark}\label{rem:2302182010}
The  Karush-Kuhn-Tucker inequality  (ev2) is equivalent to $\Psi(\alpha^\circ(s), u^\circ(s), \tilde{f}^\circ(s))=0$, by the definition of $\Psi$ \eqref{eq:2202181214}, so that the term in $\Psi$ in the energy balance gives a contribution only in the zones where the evolution is not stable. Notice also that $\Psi(\alpha^\circ(s), u^\circ(s), \tilde{f}^\circ(s))=\d_2(\partial_\alpha\E(\alpha^\circ(s), u^\circ(s)) - \tilde{f}^\circ(s), G)$, (cf.\ Remark~\ref{rem:2302182016}) a sort of $L^2$-distance from the (first order) stability set $G$.
\end{remark}
\begin{remark}\label{0106182044}
If $(t^\circ, \alpha^\circ, u^\circ)\in W^{1,\infty}(0,S; [0,T]{\times}H^1(\Omega){\times} W^{1,p}(\Omega))$ satisfies (ev0), (ev1), (ev2), then (ev3) is equivalent to say that for a.e.\ $s\in (0,S)$
\begin{equation}\label{0106182337}
\langle \de_\alpha \E(\alpha^\circ(s), u^\circ(s)) -\tilde{f}^\circ(s), \dot{\alpha}^\circ(s) \rangle + \|\dot{\alpha}^\circ(s)\|_{L^2}\,\Psi(\alpha^\circ(s), u^\circ(s), \tilde{f}^\circ(s)) =0\,,
\end{equation}
arguing as done for Lemma~\ref{le:2202181112}, by differentiation. Notice that, by definition of $\Psi$ \eqref{eq:2202181214}, we always have
\begin{equation}\label{2706180942}
\langle \de_\alpha \E(\alpha^\circ(s), u^\circ(s)) -\tilde{f}^\circ(s), \beta \rangle + \|\beta\|_{L^2}\,\Psi(\alpha^\circ(s), u^\circ(s), \tilde{f}^\circ(s)) \geq 0\,,
\end{equation}
for every $\beta \in H^1(\Omega)$ with $\beta \leq 0$ a.e.\ in $\Omega$.
\end{remark}
By the convergence of the energies \eqref{0106191253} we deduce the following relation between $\tilde{f}^\circ$ and $f(V^\circ)$.
\begin{proposition}\label{prop:ordineCum}
For every $s\in [0,S]$
\begin{equation}\label{eq:ordineCum}
\tilde{f}^\circ(s) \leq f(V^\circ) \quad\text{a.e.\ in }\Omega\,.
\end{equation}
\end{proposition}
\begin{proof}
In this proof we use the notion of essential variation of a time-dependent family of functions, whose definition is given in~Definition~\ref{def:essvar} in the Appendix. Here we recall that, by Proposition~\ref{prop:essvar},
\[
V_\epsilon^\circ(s) = \essvar(\zeta^\circ_\epsilon;0,s) = \esssup_{0=s_0 < \dots < s_m = s} \Big\{ \sum_{j=0}^m |\zeta_\epsilon^\circ(s_j) - \zeta_\epsilon^\circ(s_{j-1})| \Big\} \, .
\]
Hence, since $f$ is nonincreasing, we have that for every partition $0\leq s_0 < \dots < s_m \leq s$
\begin{equation}\label{eq:2302181030}
f(\Veo(s))\leq f\Big(\sum_{j=1}^m |\zeta_\varepsilon^\circ(s_j)-\zeta_\varepsilon^\circ(s_{j-1})|\Big)\quad\text{ a.e.\ in }\Omega\,.
\end{equation}
Indeed
\begin{equation*}
\sum_{j=1}^m |\zeta_\varepsilon^\circ(s_j)-\zeta_\varepsilon^\circ(s_{j-1})|  \leq \Veo(s)\,,
\end{equation*}
as functions on $\Omega$.
By \eqref{eq:2202182102} and \eqref{0106191253} we have that $\nabla \ueo(s) \to \nabla u^\circ(s)$ in $L^2(\Omega)$ for every $s\in [0,S]$, so that
\begin{equation*}
\zeta_\varepsilon^\circ(s) \to \zeta^\circ(s)\quad\text{in }L^2(\Omega) \,,
\end{equation*}
and then
\begin{equation}\label{0106181258}
\sum_{j=1}^m |\zeta_\varepsilon^\circ(s_j)-\zeta_\varepsilon^\circ(s_{j-1})|  \to \sum_{j=1}^m |\zeta^\circ(s_j)-\zeta^\circ(s_{j-1})| \quad\text{ in }L^2(\Omega)
\end{equation}
as $\epsilon\to 0$ for every fixed partition.
Testing \eqref{eq:convCumulEps} with characteristic functions of any Borel set $B\subset \Omega$ and employing \eqref{eq:2302181030}, \eqref{0106181258}, we can pass to $\epsilon\to 0$  and obtain
\begin{equation*}
\int\limits_B \tilde{f}^\circ(s) \d x \leq \int \limits_B f\Big(\sum_{j=1}^m |\zeta^\circ(s_j)-\zeta^\circ(s_{j-1})|\Big) \d x\,,
\end{equation*}
that gives, since $B\subset \Omega$ Borel is arbitrary,
 \begin{equation*}
  \tilde{f}^\circ(s)  \leq  f\Big(\sum_{j=1}^m |\zeta^\circ(s_j)-\zeta^\circ(s_{j-1})|\Big) \,.
 \end{equation*}
 By the arbitrariness of the partition, and since $f$ is nonincreasing, this implies
 \begin{equation*}
  \tilde{f}^\circ(s)  \leq  f\big(\mathrm{essVar}(\zeta^\circ;0,s)\big) \,,
 \end{equation*}
and \eqref{eq:ordineCum} follows because $V^\circ(s)=\mathrm{essVar}(\zeta^\circ;0,s)$ by Proposition~\ref{prop:essvar}.
\end{proof}
\begin{remark}\label{rem:2302181959}

By Proposition~\ref{prop:ordineCum} we have that the  Karush-Kuhn-Tucker inequality  (ev2) holds also for $f(V^\circ(s))$ in place of $\tilde{f}^\circ(s)$, that is for a.e.\ $s \in (0,S)\setminus U^\circ$ and $\beta \in H^1(\Omega)$, $\beta \leq 0$, we have
\begin{equation*}
\langle \de_\alpha \E(\alpha^\circ(s), u^\circ(s)) -f(V^\circ(s)), \beta \rangle \geq 0 \, .
\end{equation*}
However we can guarantee only the inequality
\[
\begin{split}
\E(\alpha^\circ(S), u^\circ(S)) &- \int_0^S \! \langle f(V^\circ(s)), \dot{\alpha}^\circ(s) \rangle\, \d s +\int_0^s \|\dot{\alpha}^\circ(s)\|_{L^2}\,\Psi(\alpha^\circ(s), u^\circ(s), \tilde{f}^\circ(s))\,\mathrm{d}s  \\&\geq \E(\alpha_0, u_0) + \int_0^S \! \langle \mu(\alpha^\circ(s)) \nabla u^\circ(s) , \nabla \dot w^\circ(s) \rangle  \, \d s  \, .
\end{split}
\]
in place of (ev3), if we consider $f(V^\circ(s))$ instead of $\tilde{f}^\circ$.
\end{remark}

\begin{proof}[Proof of Theorem~\ref{teo:evoRepar}]
Since in general $t^\circ$ is not invertible, we consider its left and right inverse, defined by
\begin{equation*}
\begin{split}
s^\circ_-(t)&:=\sup\{s\in[0,S]\colon t^\circ(s)<t\} \quad \text{ for }t\in (0,T]\,, s^\circ_-(0):=0\,,\\
s^\circ_+(t)&:=\inf\{s\in[0,S]\colon t^\circ(s)>t\} \quad \text{ for }t\in [0,T)\,,  s^\circ_+(T):=S\,.
\end{split}
\end{equation*}
For every $t\in [0,T]$ we have that  $t^\circ(s^\circ_-(t))=t=t^\circ(s^\circ_-(t))$ and
\begin{equation}\label{eq:2302181638}
s^\circ_-(t)\leq \liminf_{\epsilon\to 0} \seo(t)\leq \limsup_{\epsilon\to 0} \seo(t) \leq s^\circ_+(t)\,,
\end{equation}
while $s^\circ_-(t^\circ(s))\leq s\leq s^\circ_+(t^\circ(s))$ for every $s\in [0,S]$. The set
\begin{equation}\label{eq:defS0}
S^\circ:=\{t\in [0,T]\colon s^\circ_-(t)<s^\circ_+(t)\}
\end{equation}
is at most countable, and
\begin{equation}\label{eq:defU0}
U^\circ=\bigcup_{t\in S^\circ}(s^\circ_-(t), s^\circ_+(t))\,.
\end{equation}

Arguing as done in Proposition~\ref{prop:stweak}, by \eqref{eq:2202182102}, we pass (ev1)$_\epsilon$ to the limit and obtain (ev1), while (ev0) is immediate from the pointwise convergence of $\aeo(s)$ to $\alpha^\circ(s)$ for every $s \in [0,S]$.

\noindent \textbf{Proof of (ev2).} It is enough to show that $A^\circ \subset U^\circ$, where
\begin{equation}\label{eq:defA0}
A^\circ:=\{s\in[0,S]\colon \Psi(\alpha^\circ(s), u^\circ(s),\tilde{f}^\circ(s))>0\}\,.
\end{equation}
Arguing as in the proof Proposition~\ref{prop:stweak} to obtain \eqref{eq:2302181101} and \eqref{eq:2302181104}, we deduce that
for every $\beta\in H^1_-(\Omega)$ and every $s \in [0,S]$
\begin{equation*}
\langle-\de_\alpha \E(\alpha^\circ(s), u^\circ(s)) + \tilde{f}^\circ(s), \beta \rangle \leq \liminf_{\epsilon\to 0} \langle-\de_\alpha \E(\aeo(s), \ueo(s)) + f(\Veo(s)), \beta \rangle\,,
\end{equation*}
so that, passing to the supremum for $\beta \in H^1_-(\Omega)$,
\begin{equation}\label{eq:2302181132}
\Psi(\alpha^\circ(s), u^\circ(s),\tilde{f}^\circ(s)) \leq \liminf_{\epsilon\to 0} \Psi(\aeo(s), \ueo(s),f(\Veo(s)))\,.
\end{equation}
By \eqref{eq:derivative of E} and the convergences \eqref{eq:2202182102} we have that the map $s\mapsto \langle-\de_\alpha \E(\alpha^\circ(s), u^\circ(s)), \beta \rangle$ is continuous for every $\beta\in H^1_-(\Omega)$. Also $s \mapsto \langle \tilde{f}^\circ(s), \beta \rangle$ is continuous: indeed
\begin{equation*}
|\langle f(\Veo(s_2))- f(\Veo(s_1)), \beta \rangle|\leq\|f'\|_{L^\infty}  \int \limits_\Omega |\beta| \int_{s_1}^{s_2} | \dot{\zeta}^\circ_\epsilon(\sigma) | \d \sigma \d x \leq  \|f'\|_{\infty} \Big(\int_{s_1}^{s_2}\|\dot{\zeta}^\circ_\epsilon(\sigma)\|_{L^2} \d \sigma \Big) \|\beta\|_{L^2} \leq C(s_2-s_1) \|\beta\|_{L^2}\,,
\end{equation*}
since
\begin{equation*}
\|\dot{\zeta}_\varepsilon^\circ(\sigma)\|_{L^2}= \| g'(\aeo(\sigma))\,\dot{\alpha}_\varepsilon^\circ(\sigma) \, \nabla \ueo(\sigma) + g(\aeo(\sigma))\, \nabla \dot{u}_\varepsilon^\circ(\sigma)\|_{L^2} \leq C \big(\|\dot{\alpha}_\varepsilon^\circ(\sigma)\|_{H^1} + \|\nabla \dot{u}_\varepsilon^\circ(\sigma)\|_{L^2}\big) \leq C\,,
\end{equation*}
and we pass to the limit as $\epsilon\to 0$ to get
\begin{equation*}
|\langle \tilde f^\circ(s_2)- \tilde f^\circ (s_1), \beta \rangle| \leq C(s_2-s_1) \|\beta\|_2\,.
\end{equation*}
Therefore $s \mapsto \langle-\de_\alpha \E(\alpha^\circ(s), u^\circ(s))+ \tilde{f}^\circ(s), \beta \rangle$ is continuous, and
\begin{equation}\label{eq:2302181232}
s\mapsto \Psi(\alpha^\circ(s), u^\circ(s),\tilde{f}^\circ(s)) \quad\text{ is lower semicontinuous.}
\end{equation}
In particular, $A^\circ$ is an open set. We now follow closely the argument in \cite[Theorem~5.4, proof of (ev3) therein]{CriLaz16}, to say that
\begin{equation}\label{lim''}
\limsup_{\varepsilon\rightarrow 0}\,\dot{t}_\varepsilon^\circ(s)>0\quad \text{for a.e.\ $s\in (0,S)\setminus D^\circ$}\,,
\end{equation}
where $D^\circ:=\{s\in (0,S)\,\colon\,\dot{t}^\circ(s)=0\}$.

Arguing by contradiction, there exists a measurable set $B\subset (0,S)\setminus D^\circ$ with positive measure such that
\[ \lim_{\varepsilon\rightarrow 0}\dot{t}_\varepsilon^\circ(s)=0\quad \text{for every $s\in B$}\,,
\]
$\teo$ being nondecreasing.
Since the functions $\teo$ are 1-Lipschitz, the Dominated Convergence Theorem implies that
\[\lim_{\varepsilon\rightarrow 0}\int_B\dot{t}_\varepsilon^\circ(s)\,\mathrm{d}s=0\,.\]
On the other hand, from $\teo \weak t^\circ$ weakly$^*$ in $W^{1,\infty}$ (see \eqref{eq:2202182102})
\[\lim_{\varepsilon\rightarrow 0}\int_B\dot{t}_\varepsilon^\circ(s)\,\mathrm{d}s=\int_B\dot{t}^\circ(s)\,\mathrm{d}s\,,\]
and this contradicts
\[\int_B\dot{t}^\circ(s)\,\mathrm{d}s>0\,,\]
that follows from the definition of $D^\circ$.

By \eqref{eq:2302181132} and \eqref{eq:2202181319} evaluated in $t=\teo(s)$ (cf.\ \eqref{eq:2302181228}) we deduce
\begin{equation*}
\begin{split}
0\leq \Psi(\alpha^\circ(s),u^\circ(s), \tilde{f}^\circ(s))\leq \liminf_{\varepsilon\rightarrow 0}\Psi(\aeo(s),\ueo(s), f(\Veo(s)))=\liminf_{\varepsilon\rightarrow 0} \varepsilon\|\dot{\alpha}_\varepsilon(\teo(s))\|_{L^2}=\liminf_{\varepsilon\rightarrow 0} \varepsilon\frac{\|\dot{\alpha}_\varepsilon^\circ(s)\|_{L^2}}{\dot{t}_\varepsilon^\circ(s)}=0
\end{split}
\end{equation*}
for a.e.\ $s\in (0,S)\setminus D^\circ$.
This implies that $\dot{t}^\circ(s)=0$ for a.e.\ $s\in A^\circ$. Being $A^\circ$ open by \eqref{eq:2302181232}, every $s\in A^\circ$ has an open neighborhood where $\dot{t}^\circ=0$; then $A^\circ \subset U^\circ$, because $t^\circ$ is Lipschitz.

\noindent \textbf{Proof of (ev3).} Looking at the version of the energy balance \eqref{eq:2302181249} proven in Proposition~\ref{prop:2202181211}, this is invariant under time reparametrisation. Then, by the change of variables $t=\teo(s)$  (in the left hand side)  we get 
\begin{equation}\label{eq:2302181257}
 \begin{split}
\E(\aeo(S), \ueo(S)) + \int_0^S \! \R(\dot \alpha^\circ_\epsilon(s); \Veo(s))\, \d s  &+ \int_0^S \! \| \dot \alpha_\epsilon^\circ(s) \|_{L^2} \Psi(\aeo(s), \ueo(s), f(\Veo(s))) \, \d s \\&= \E(\alpha_0, u_0) +  \int_0^T \! \langle \mu(\alpha_\epsilon(t)) \nabla u_\epsilon(t) , \nabla \dot w(t) \rangle_{L^2}  \, \d t  \, .
\end{split}
\end{equation} 
Arguing as done in Proposition~\ref{prop:balance weak} to deduce \eqref{eq:fatigue is lsc}, we obtain
\begin{equation*}
- \int_0^S \! \langle\tilde{f}^\circ(s), \dot \alpha^\circ(s)\rangle \, \d s=  \sup_{0=s_0 < \dots < s_m = S} \Big\{ \sum_{j=1}^m \langle\tilde{f}^\circ(s_j), \alpha^\circ(s_{j-1}) - \alpha^\circ(s_j) \rangle \Big\}
\end{equation*}
and then, since \eqref{eq:2202182102} and \eqref{eq:convCumulEps} give
\[
\langle\tilde{f}^\circ(s_j), \alpha^\circ(s_{j-1}) - \alpha^\circ(s_j) \rangle=\lim_{\varepsilon \to 0} \langle f(\Veo(s_j)), \aeo(s_{j-1}) - \aeo(s_j) \rangle\]
 for any $s_{j-1}$, $s_j$, we deduce
\begin{equation}\label{eq:2302181450}
- \int_0^S \! \langle\tilde{f}^\circ(s), \dot \alpha^\circ(s)\rangle \, \d s\leq \liminf_{\epsilon\to 0} \int_0^S \! \R(\dot \alpha^\circ_\epsilon(s); \Veo(s))\, \d s\,,
\end{equation}
recalling \eqref{eq:int of f as sup}.
Moreover, we claim that
\begin{equation}\label{eq:2302181459}
\int_{A^\circ} \|\dot{\alpha}^\circ(s)\|_{L^2}\,\Psi(\alpha^\circ(s), u^\circ(s), \tilde{f}^\circ(s))\,\mathrm{d}s\leq\liminf_{\varepsilon\rightarrow 0}\int_{A^\circ}  \| \dot \alpha_\epsilon^\circ(s) \|_{L^2} \Psi(\aeo(s), \ueo(s), f(\Veo(s))) \, \d s\,.
\end{equation}
Indeed, for every compact set $C\subset A^\circ$ and every continuous function $\psi\colon C \to [0,+\infty)$ such that
\begin{equation*}
\Psi(\alpha^\circ(s), u^\circ(s), \tilde{f}^\circ(s)) >\psi(s)\quad\text{for every $s\in C$}\,,
\end{equation*}
by the compactness of $C$ and \eqref{eq:2302181132},
for $\varepsilon$ sufficiently small we get
\begin{equation*}
\Psi(\aeo(s), \ueo(s), f(\Veo(s)))>\psi(s)\quad\text{for every $s\in C$}\,.
\end{equation*}
Then, by approximating the semicontinuous function $s\mapsto \Psi(\alpha^\circ(s), u^\circ(s), \tilde{f}^\circ(s))$ from below by continuous functions, in order to prove \eqref{eq:2302181459} it is sufficient to show
\begin{equation*}\label{618sol}
\int_C\|\dot{\alpha}^\circ(s)\|_{L^2}\,\psi(s)\,\mathrm{d}s\leq \liminf_{\varepsilon\rightarrow 0} \int_C\|\dot{\alpha}_\varepsilon^\circ(s)\|_{L^2}\,\psi(s)\,\mathrm{d}s
\end{equation*}
for every compact $C\subset A^\circ$ and every continuous function $\psi\colon C \to [0,+\infty)$.
This is done as in \cite[Theorem~5.4]{CriLaz16} or \cite[Lemma~6.4]{DMDSSol}, using a localisation argument and the fact that for every $\varphi \in C_c(\Omega)$ with $\|\varphi\|_{L^2}=1$ the functions $s\mapsto\langle\varphi,\dot{\alpha}_\varepsilon^\circ(s)\rangle$ are equi-Lipschitz on $[0,S]$ and converge to $s\mapsto\langle\varphi,\dot{\alpha}^\circ(s)\rangle$ for every $s$.

By \eqref{eq:2302181450}, \eqref{eq:2302181459}, and the semicontinuity of the internal energy (cf.\ \eqref{eq:E is lsc}) we obtain the lower semicontinuity of the left hand side of the energy balance \eqref{eq:2302181257}.

Let us now study the limit with respect to $\epsilon$ of the right hand side of \eqref{eq:2302181257}.
Since for every $t\in [0,T]\setminus S^\circ$ it holds that $s^\circ_-(t)=\lim_{\epsilon\to 0} \seo(t)$ (see \eqref{eq:2302181638}), then
\begin{equation*}
\alpha_\epsilon(t)\weak \alpha^\circ(s^\circ_-(t))\quad\text{ in }H^1(\Omega)\,,\qquad u_\epsilon(t)\weak u^\circ(s^\circ_-(t)) \quad\text{ in }W^{1,p}(\Omega)\,,
\end{equation*}
and
\begin{equation*}
\int_0^T  \langle \mu(\alpha^\circ(s^\circ_-(t))) \nabla u^\circ(s^\circ_-(t)), \nabla \dot w(t) \rangle_{L^2} \d t =\lim_{\epsilon\to 0}\int_0^T \! \langle \mu(\alpha_\epsilon(t)) \nabla u_\epsilon(t) , \nabla \dot w(t) \rangle_{L^2}  \, \d t
\end{equation*}
by the Dominated Convergence Theorem. On the other hand,  recalling~\eqref{eq:wcirc}, 
\begin{equation*}
\begin{split}
\int_0^T  \langle \mu(\alpha^\circ(s^\circ_-(t))) \nabla u^\circ(s^\circ_-(t)), \nabla \dot w(t) \rangle_{L^2} \d t &=\int_0^S  \langle \mu(\alpha^\circ(s^\circ_-(t^\circ(s)))) \nabla u^\circ(s^\circ_-((t^\circ(s))), \nabla \dot w(t^\circ(s)) \, \dot{t}^\circ(s) \rangle_{L^2} \d s \\&= \int_0^S  \langle \mu(\alpha^\circ(s)) \nabla u^\circ(s), \nabla \dot w^\circ(s) \rangle_{L^2} \d s\,,
\end{split}
\end{equation*}
since $\dot{t}^\circ(s)=0$ for a.e.\ $s\in U^\circ$ and $s^\circ_-((t^\circ(s))=s$ for a.e.\ $s\in (0,S)\setminus U^\circ$.
Therefore the right hand side of \eqref{eq:2302181257} passes to the limit and we conclude the energy  inequality
\[
\begin{split}
\E(\alpha^\circ(S), u^\circ(S)) &- \int_0^S \! \langle \tilde{f}^\circ(s), \dot{\alpha}^\circ(s) \rangle\, \d s +\int_0^S \|\dot{\alpha}^\circ(s)\|_{L^2}\,\Psi(\alpha^\circ(s), u^\circ(s), \tilde{f}^\circ(s))\,\mathrm{d}s  \\& \leq \E(\alpha_0, u_0) + \int_0^S \! \langle \mu(\alpha^\circ(s)) \nabla u^\circ(s) , \nabla \dot w^\circ(s) \rangle  \, \d s  \, .
\end{split}
\] 
To prove the converse inequality, we differentiate with respect to the time variable the energy (which is absolutely continuous, since $\alpha^\circ$, $u^\circ$ are Lipschitz). We obtain that for a.e.\ $s \in (0,S)$
\begin{equation}\label{2610181048}
\begin{split}
\frac{\d}{\d s} \Big[ \E(\alpha^\circ(s), u^\circ(s)) \Big] & = \langle \de_\alpha \E(\alpha^\circ(s), u^\circ(s) ) , \dot \alpha^\circ(s) \rangle + \langle \de_u \E(\alpha^\circ(s), u^\circ(s) ) , \dot u^\circ(s)\rangle  \\
& = \langle \de_\alpha \E(\alpha^\circ(s), u^\circ(s) ) , \dot \alpha^\circ(s) \rangle + \langle \mu(\alpha^\circ(s)) \nabla u^\circ(s), \nabla \dot w^\circ(s) \rangle_{L^2}\,, \\
& \geq \langle \tilde{f}^\circ(s), \dot{\alpha}^\circ(s) \rangle - \|\dot \alpha^\circ(s)\|_{L^2} \Psi(\alpha^\circ(s),u^\circ(s),\tilde f^\circ(s))+ \langle \mu(\alpha^\circ(s)) \nabla u^\circ(s), \nabla \dot w^\circ(s) \rangle_{L^2}\,,
\end{split}
\end{equation}
employing (ev1), evaluated in $\beta=\dot{\alpha}^\circ(s)$ (which is in $H^1_-(\Omega)$ for a.e.\ $s$), in the second equality and~\eqref{2706180942} in the inequality above.
We deduce the energy balance (ev3) by integrating \eqref{2610181048} in $(0,S)$
As a byproduct, we also obtain that \eqref{eq:2302181450} and \eqref{eq:2302181459} hold true as limits as $\varepsilon \to 0$. 
\end{proof}

We conclude by showing some properties of an evolution $(t^\circ, \alpha^\circ, u^\circ)$, obtained as limit of rescaled $\varepsilon$-approximate viscous evolution, in the spirit of e.g.\ \cite{DMDSSol, KRZ13a, CriLaz16}.
We are in particular interested in its description in the time subset $U^\circ \subset [0,S]$, where it is not rate independent: if $\alpha^\circ$ remains constant in $(s_1,s_2)\subset U^\circ$, then all the evolution is trivial in $(s_1,s_2)$ (Remark~\ref{rem:Le6.1}); on the other hand, if $\dot{\alpha}^\circ > 0$ in space in a time interval, up to a further time reparametrisation we have that the system is governed by an equation satisfied in the transition between the initial and final configurations: this equation (see \eqref{0106182321} and \eqref{3006181741} in Proposition~\ref{prop:0106182313}) corresponds formally to consider 1 as viscosity parameter in \eqref{3006181736} in Lemma~\ref{le:2202181112}, governing the $\varepsilon$-approximate viscous evolutions.

\begin{remark}\label{rem:Le6.1}
If $\dot{\alpha}^\circ(s)=0$ for every $s\in (s_1,s_2)\subset U^\circ$, then $\dot{u}^\circ(s)=0$ for every $s\in (s_1,s_2)\subset U^\circ$. Indeed, by definition of $U^\circ$, it follows that $t^\circ(s)=t^\circ(s_1)$ and $w^\circ(s)=w^\circ(s_1)$ for every $s\in (s_1,s_2)$, and then $u^\circ(s)=u^\circ(s_1)$, the unique solution of
\begin{equation*}
\min_{u=w^\circ(s_1) \text{ on }\partial_D \Omega} \int \limits_\Omega \mu(\alpha^\circ(s_1)) |\nabla u|^2 \d x\,,
\end{equation*}
by (ev1) in Theorem~\ref{teo:evoRepar}.
\end{remark}

\begin{proposition}\label{prop:0106182313}
Let $(s_1,s_2)\subset A^\circ$ (with $A^\circ$ defined in \eqref{eq:defA0}) containing no subintervals where $\dot{\alpha}^\circ(s)=0$ in $\Omega$ for a.e.\ $s$, and let for every $s\in (s_1,s_2)$
\begin{equation*}
\varrho(s):=\int_{\frac{s_1+s_2}{2}}^s \frac{\|\dot{\alpha}^\circ(\sigma)\|_{L^2}}{\Psi(\alpha^\circ(\sigma), u^\circ(\sigma), \tilde{f}^\circ(\sigma)) } \d \sigma\,.
\end{equation*}
Then $\varrho$ is locally Lipschitz and stricly increasing, the function
\begin{equation*}
\alpha^\sharp(r):=\alpha^\circ(\varrho^{-1}(r)) \quad \text{for }r\in \varrho((s_1,s_2))
\end{equation*}
has bounded variation and is continuous into $H^1(\Omega)$, and
\begin{equation}\label{0106182321}
\|\dot{\alpha}^\circ(\varrho^{-1}(r))\|_{L^2}^2 \Big[\big\langle\de_\alpha \E(\alpha^\sharp(r), u^\sharp(r)) -\tilde{f}^\sharp(r), \,\dot{\alpha}^\sharp(r) \big\rangle + \| \dot{\alpha}^\sharp(r)\|_{L^2}^2   \Big] = 0\,,
\end{equation}
for every $r\in \varrho^{-1}((s_1,s_2))$,
where $u^\sharp(r):=u^\circ(\varrho^{-1}(r))$, $\tilde{f}^\sharp(r):=\tilde{f}^\circ(\varrho^{-1}(r))$.

If, moreover, $\dot{\alpha}^\circ(s)$ is not 0 for every $s\in (s_1,s_2)$ and $\|\dot{\alpha}^\circ(s)\|_{L^2} > \delta_K$ for every $K\compact (s_1,s_2)$, then $\varrho$ is locally bi-Lipschitz, $\alpha^\sharp$ is locally Lipschitz, and
\begin{equation}\label{3006181741}
\big\langle\de_\alpha \E(\alpha^\sharp(r), u^\sharp(r)) -\tilde{f}^\sharp(r), \, \dot{\alpha}^\sharp(r)\big \rangle + \| \dot{\alpha}^\sharp(r)\|^2_{L^2} = 0 \,.
\end{equation}
\end{proposition}
\begin{proof}
By \eqref{eq:defA0} and \eqref{eq:2302181232}, for any $K\compact A^\circ$ we get that $\Psi(\alpha^\circ(\sigma), u^\circ(\sigma), \tilde{f}^\circ(\sigma)) \geq \delta_K > 0$ for $\sigma\in K$. Then $\varrho$ is locally Lipschitz on $(s_1,s_2)$, and it is strictly increasing by the assumptions that in no subintervals of $(s_1,s_2)$ we have $\dot{\alpha}^\circ(s)=0$ in $\Omega$. This gives also $\varrho^{-1}$ continuous and strictly increasing, so that $\alpha^\sharp$ is continuous and with bounded variation, since $\alpha^\circ$ is Lipschitz.

The change of variables $s=\varrho^{-1}(r)$ in \eqref{0106182337} gives
\begin{equation*}
\begin{split}
  \big\langle\de_\alpha \E(\alpha^\sharp(r), u^\sharp(r)) -\tilde{f}^\sharp(r), \dot{\alpha}^\circ(\varrho^{-1}(r)) \big\rangle   + \Psi(\alpha^\sharp(r), u^\sharp(r), \tilde{f}^\sharp(r)) \, \|\dot{\alpha}^\circ(\varrho^{-1}(r))\|_{L^2} = 0\,,
\end{split}
\end{equation*}
that is
\begin{equation}\label{0206180002}
\begin{split}
\|\dot{\alpha}^\circ(\varrho^{-1}(r))\|_{L^2}  \big\langle\de_\alpha \E(\alpha^\sharp(r), u^\sharp(r)) -\tilde{f}^\sharp(r), \dot{\alpha}^\circ(\varrho^{-1}(r)) \big\rangle   + \Psi(\alpha^\sharp(r), u^\sharp(r), \tilde{f}^\sharp(r)) \, \|\dot{\alpha}^\circ(\varrho^{-1}(r))\|_{L^2}^2 = 0\,,
\end{split}
\end{equation}
for every $r\in \varrho^{-1}((s_1,s_2))$ and every $\beta\in H^1_-(\Omega)$.
Now $\alpha^\sharp$ is weakly differentiable in $H^1(\Omega)$ at a.e.\ $r\in \varrho^{-1}((s_1,s_2))$, and we have the chain rule
\begin{equation}\label{0206180952}
\dot{\alpha}^\sharp(r)=\dot{\alpha}^\circ(\varrho^{-1}(r)) \frac{\d}{\d t}\varrho^{-1}(r) = \dot{\alpha}^\circ(\varrho^{-1}(r)) \frac{\Psi(\alpha^\sharp(r), u^\sharp(r), \tilde{f}^\sharp(r)) }{\|\dot{\alpha}^\circ(\varrho^{-1}(r))\|_{L^2}} \quad \text{a.e.\ in }\Omega
\end{equation}
for a.e.\ $r$ such that $\|\dot{\alpha}^\circ(\varrho^{-1}(r))\|_{L^2}>0$. Then \eqref{0206180002} and \eqref{0206180952} imply
\begin{equation*}
\begin{split}
\|\dot{\alpha}^\circ(\varrho^{-1}(r))\|_{L^2} \Big[ \big\langle \de_\alpha \E(\alpha^\sharp(r), u^\sharp(r)) -\tilde{f}^\sharp(r), \dot{\alpha}^\circ(\varrho^{-1}(r)) \big\rangle  + \langle \dot{\alpha}^\sharp(r), \dot{\alpha}^\circ(\varrho^{-1}(r)) \rangle_{L^2} \Big] = 0\,,
\end{split}
\end{equation*}
Recalling that  $\Psi(\alpha^\sharp(r), u^\sharp(r), \tilde{f}^\sharp(r))>0$ for every $r\in \varrho^{-1}(s_1,s_2)$, the two previous inequalities imply~\eqref{0106182321}.
At this stage, \eqref{3006181741} follows easily since $\|\dot{\alpha}^\circ(\varrho^{-1}(r))\|_{L^2}^2>0$ for every $r \in \varrho^{-1}((s_1,s_2))$.
\end{proof}

As usual in an analysis based on time rescaling, one could see that in the original, \emph{faster}, time variable $t\in [0,T]$, the evolution is rate-independent outside an at most countable number of jump times, which is a subset of $S^\circ$ introduced in \eqref{eq:defS0}. In order to describe the evolution of the system during these jump, one has to employ the description given by Remark~\ref{rem:Le6.1} and Proposition~\ref{prop:0106182313}.
Here we do not perform directly this analysis, based on inverse rescaling in time, since it would be very similar to that in e.g.\ \cite[Section~5]{Sol2} and \cite[Proposition~6.7]{CriLaz16}, to which we refer the interested reader.

\begin{appendices}

\section{Auxiliary results}\label{sec:app}

\subsection*{The essential variation} In this appendix $X$ denotes a measure space. We do not label the measure on $X$ and the notions of $L^p$ space and of a.e.-equivalence refer to the measure on $X$. Moreover we fix $n \geq 1$.

We define here the notion of essential variation, namely the variation for a time-dependent family of measurable functions, in the sense of a.e.\ inequality.

\begin{definition} \label{def:essvar}
Let us consider a function $t \mapsto \zeta(t)$, with $\zeta(t) \colon X \to \RR^n$. Let $0 \leq s \leq t \leq T$. The {\em essential variation} of $\zeta$ in the interval $[s,t]$ is the function $\essvar(\zeta;s,t) \colon X \to [0,+\infty]$ defined by
\[
\essvar(\zeta;s,t) := \esssup_{s=s_0 < \dots < s_m = T} \Big\{ \sum_{j=0}^m |\zeta(s_j) - \zeta(s_{j-1})| \Big\} \, ,
\]
where the essential supremum is taken over all partitions $0=s_0 < \dots < s_m = t$, $m \in \NN$.
\end{definition}

\begin{remark} \label{rmk:essvar is additive}
For every $t_1 \leq t_2 \leq t_3$ we have
\[
\essvar(\zeta;t_1,t_3) = \essvar(\zeta;t_1,t_2) + \essvar(\zeta;t_2,t_3) \, \quad \text{a.e.\ in } X \, .
\]
\end{remark}

For completeness, we recall here the definition of the essential supremum of a family of measurable functions.
\begin{definition}\label{def:esssup}
Let $(v_i)_{i\in I}$ be a family of measurable functions from $X$ to $[-\infty,\infty]$. Let $\ol v \colon X \to [-\infty,\infty]$ be a measurable function such that
\begin{itemize}
\item[(i)] $\ol v \geq v_i$ a.e.\ in $X$, for every $i\in I$;
\item[(ii)] if $v \colon X \to [-\infty,\infty]$ is a measurable function such that $v \geq v_i$ a.e.\ in $X$, for every $i\in I$, then $v \geq \ol v$ a.e.\ in $X$.
\end{itemize}
The functions $\ol v$ is called an \emph{essential supremum} of the family $(v_i)_{i\in I}$. In fact, there exists a unique (up to a.e.\ equivalence) essential supremum $\ol v$ of the family $(v_i)_{i \in I}$. We denote it by $ \displaystyle \esssup_{i \in I} v_i := \ol v$.
\end{definition}

In the next proposition we provide an explicit formula for the essential variation of a function $\zeta$ that is absolutely continuous in time. A quick survey about the notion and the main properties of the Bochner integral can be found in the appendix of~\cite{Bre}; for a more detailed treatment of the subject we refer to~\cite{Dun-Sch}.

\begin{proposition} \label{prop:essvar}
Let $p \in [1, \infty)$, and let $\zeta \in AC([0,T];L^p(X;\RR^n))$. Then
\[
\essvar(\zeta;0,t)(x) = \int_0^t \! |\dot \zeta(r;x) | \d r \, , \quad \text{for a.e.\ } x \in X \, ,
\]
where the integral in the right-hand side is a Bochner integral in $L^p(X)$.
\end{proposition}
\begin{proof}
We start by claiming that
\begin{equation} \label{eq:claim on essvar}
\essvar(\zeta;0,\cdot) \in AC([0,T];L^p(X)) \quad \text{and} \quad \frac{\d}{\d t} \essvar(\zeta;0,t) = |\dot \zeta(t)| \text{ in } L^p(X) \, .
\end{equation}
To prove the claim, let us fix $s \leq t$ and a partition $s=s_0 < \dots < s_m = t$. By the absolute continuity of $\zeta$ we obtain that
\begin{equation} \label{eq:12061101}
\sum_{j=1}^m |\zeta(s_j) - \zeta(s_{j-1})| = \sum_{j=1}^m \Big| \int_{s_{j-1}}^{s_j} \! \dot \zeta(r)  \, \d r \Big|  \leq  \sum_{j=1}^m \int_{s_{j-1}}^{s_j} \! | \dot \zeta(r) |  \, \d r = \int_{s}^{t} \! | \dot \zeta(r) |  \, \d r \,
\end{equation}
a.e.\ in $X$, where the last integral is a Bochner integral in $L^p(X)$. Note that the second inequality in~\eqref{eq:12061101} can be proven, e.g., with an approximation argument via step functions. Taking the essential supremum in~\eqref{eq:12061101}, by Remark~\ref{rmk:essvar is additive} we deduce that
\begin{equation} \label{eq:bound on essvar}
\essvar(\zeta;0,t) - \essvar(\zeta;0,s)\leq \int_{s}^{t} \! | \dot \zeta(r) |  \, \d r \quad \text{ a.e.\ in } X \, .
\end{equation}
Inequality~\eqref{eq:bound on essvar} computed for $s = 0$ yields, in particular, that $\essvar(\zeta; 0, t) \in L^p(X)$ for every $t \in [0,T]$. Moreover, it shows that $\essvar(\zeta;0,\cdot) \in AC([0,T];L^p(X))$. By~\eqref{eq:bound on essvar} and by Lebesgue's Differentiation Theorem for vector-valued functions \cite[p.\ 217]{Dun-Sch} we get that
\[
\frac{\d}{\d t} \essvar(\zeta;0,t) = \lim_{s\to t^-} \frac{\essvar(\zeta;0,t) - \essvar(\zeta;0,s)}{t-s}  \leq | \dot \zeta(t) |
\]
if $t$ is a differentiability point for $\essvar(\zeta;0,\cdot)$ and it is a Lebesgue point for $|\dot \zeta|$, the limit being taken with respect to the $L^p$-norm.

On the other hand, $s < t$ is a particular partition of the interval $[s,t]$, therefore
\[
\frac{|\zeta(t) - \zeta(s)|}{t-s} \leq \frac{\essvar(\zeta; 0, t) - \essvar(\zeta; 0, s)}{t-s} \, \quad \text{a.e.\ in } X \, .
\]
Taking the limit as $s \to t^-$ with respect to the $L^p$-norm of both sides, we obtain
\[
| \dot \zeta(t) |  \leq \frac{\d}{\d t} \essvar(\zeta;0,t)  \, \quad \text{a.e.\ in } X \, ,
\]
if $t$ is a differentiability point for $\essvar(\zeta;0,\cdot)$ and $\zeta$. This proves that $\frac{\d}{\d t} \essvar(\zeta;0,t) =  | \dot \zeta(t) |$.

Finally, since $\essvar(\zeta;0,\cdot) \in AC([0,T];L^p(X))$, we conclude that
\[
\essvar(\zeta;0,t) = \int_0^t \! \frac{\d}{\d t} \essvar(\zeta;0,r) \, \d r = \int_0^t \! | \dot \zeta(t) | \, \d r \, \quad \text{a.e.\ in } X \, .
\]

\end{proof}

\end{appendices}

\bigskip
\noindent {\bf Acknowledgements.}
The authors wish to thank Adriana Garroni for several interesting discussions and fruitful advices.

 Roberto Alessi has been supported by the MATHTECH-CNR-INdAM project and the MIUR-DAAD Joint Mobility Program: ``Variational approach to fatigue phenomena with phase-field models: modeling, numerics and experiments''.
Vito Crismale has been supported by a public grant as part of the \emph{Investissement d'avenir} project, reference ANR-11-LABX-0056-LMH, LabEx LMH, and acknowledge the financial support from the Laboratory Ypatia and the CMAP. He is currently funded by the Marie Sk\l odowska-Curie Standard European Fellowship No.\ 793018.
Gianluca Orlando has been supported by the Alexander von Humboldt Foundation.

Vito Crismale and Gianluca Orlando acknowledge the
kind hospitality of the Department of Mathematics of Sapienza University of Rome, where part of this research was developed.

\bigskip
\bibliographystyle{siam}

\end{document}